\documentclass{amsart}

\usepackage{graphicx}
\usepackage[hidelinks]{hyperref}
\usepackage{caption}
\usepackage{subcaption}
\usepackage[table,dvipsnames]{xcolor}
\usepackage{cancel}
\usepackage{verbatim} 
\usepackage[normalem]{ulem} 
\usepackage{makecell} 
\usepackage{adjustbox} 

\theoremstyle{definition}
\newtheorem{theorem}{Theorem}[section]
\makeatletter
    \def\@endtheorem{\hfill$\Diamond$\endtrivlist\@endpefalse}
\makeatother
\newtheorem{definition}[theorem]{Definition}
\makeatletter
    \def\@endtheorem{\endtrivlist\@endpefalse}
\makeatother

\newtheorem{lemma}[theorem]{Lemma}
\newtheorem{remark}[theorem]{Remark}
\newtheorem{example}[theorem]{Example}

\usepackage{algorithm}
\usepackage[indLines=false]{algpseudocodex}

\usepackage{makerobust} 

\makeatletter
\newcommand{\leqnomode}{\tagsleft@true}
\newcommand{\reqnomode}{\tagsleft@false}
\makeatother

\newcommand\hcancel[2][black]{\setbox0=\hbox{$#2$}%
\rlap{\raisebox{.45\ht0}{\textcolor{#1}{\rule{\wd0}{1pt}}}}#2}

\newcommand{\FastZigguDir}{\textsc{FastZigguWithDirections}}
\newcommand{\FastZiggu}{\textsc{FastZiggu}}

\newcommand{\Z}{\textbf{Z}}
\newcommand{\Core}{\textsf{core}}
\newcommand{\reflect}{\mathsf{reflect}}
\newcommand{\inv}{\mathsf{inv}}
\newcommand{\first}{\mathsf{first}}
\newcommand{\last}{\mathsf{last}}
\newcommand{\ruler}{\mathsf{ruler}}
\newcommand{\seqCore}{\mathsf{core}}
\newcommand{\set}[1]{\left\{#1\right\}}
\newcommand{\OEIS}[1]{\textsc{Oeis} #1}

\newcommand{\rulerBinary}[1]{\mathsf{rulerBRGC}_{#1}}
\newcommand{\srulerBinary}[1]{\mathsf{rulerBRGC}_{#1}{\pm}}
\newcommand{\rulerQuat}[1]{\Delta\QUAT{#1}}
\newcommand{\srulerQuat}[1]{\Delta^{\pm}{\QUAT{#1}}}
\newcommand{\rulerLong}[1]{\Delta{\LONG{#1}}}
\newcommand{\srulerLong}[1]{\Delta^{\pm}{\LONG{#1}}}
\newcommand{\rulerShort}[2][]{\Delta{\SHORT{#2}^{#1}}}
\newcommand{\srulerShort}[2][]{\Delta^{\pm}{\SHORT{#2}^{#1}}}
\newcommand{\reverse}[1]{{#1}^{R}}
\newcommand{\complement}[1]{\overline{#1}}
\newcommand{\complementMax}[1]{\overline{\mathsf{max}}\left(#1\right)}
\newcommand{\complementMid}[1]{\overline{\mathsf{middle}}\left(#1\right)}

\newcommand{\BRGC}[1]{\mathsf{BRGC}(#1)}
\newcommand{\lexBinary}[1]{\mathsf{lex}(#1)}

\newcommand{\nextGC}[1]{{\mathsf{next}\mathbb Q}(#1)}
\newcommand{\nextLong}[1]{\mathsf{next}{\mathbb V}(#1)}
\newcommand{\nextShort}[1]{\mathsf{next}{\mathbb S}(#1)}
\newcommand{\Q}{\mathbb Q}
\newcommand{\V}{\mathbb V}
\renewcommand{\S}{\mathbb S}
\newcommand{\QUAT}[1]{\mathbb{Q}_{#1}}
\newcommand{\LONG}[1]{\mathbb{V}_{#1}} 
\newcommand{\SHORT}[1]{\mathbb{S}_{#1}}
\newcommand{\cdotsSquish}{\,{\cdot}{\cdot}{\cdot}\,}
\newcommand{\rank}{\mathsf{rank}}
\newcommand{\rankBRGC}[1]{\mathsf{rank}\mathsf{BRGC}(#1)}
\newcommand{\rankQuat}[1]{\mathsf{rank}\Q(#1)}
\newcommand{\rankLong}[1]{\mathsf{rank}\V(#1)}
\newcommand{\rankShort}[1]{\mathsf{rank}\S(#1)}

\newcommand{\Aaron}[1]{\textcolor{purple}{#1}}
\newcommand{\MG}[1]{\textcolor{blue}{#1}}

\newcommand{\Nurikabe}[1]{\mathcal{N}_{#1}}

\newcommand{\digita}[1]{\textbf{\textcolor{red}{#1}}}
\newcommand{\digitb}[1]{\textbf{\textcolor{orange}{#1}}}
\newcommand{\digitc}[1]{\textbf{\textcolor{Dandelion}{#1}}}
\newcommand{\digitd}[1]{\textbf{\textcolor{ForestGreen}{#1}}}
\newcommand{\digite}[1]{\textbf{\textcolor{blue}{#1}}}

\newcommand{\Visit}{\textbf{visit }}

\usepackage{tikz}
\usepackage{tkz-euclide}
\usetikzlibrary{decorations.pathreplacing}
\usetikzlibrary{fit}
\usepackage{placeins}
\usepackage{array}

\newcommand{\col}[2]{
    \tikz[scale=.15]{
    \tkzDefPoint(0,0){A} \tkzDefPoint(1,0){B}
    \tkzDefSquare(A,B)
    \tkzDrawPolygon[black](A,B,tkzFirstPointResult,%
    tkzSecondPointResult)
    \tkzFillPolygon[#1](A,B,tkzFirstPointResult,%
    tkzSecondPointResult)
    \tkzDefSquare(B,A)
    \tkzDrawPolygon[black](B,A,tkzFirstPointResult,%
    tkzSecondPointResult)
    \tkzFillPolygon[#2](B,A,tkzFirstPointResult,%
    tkzSecondPointResult)
}\hspace{0.085cm}}
\newcommand{\colbb}{\col{black}{black}}
\newcommand{\colwb}{\col{white}{black}}
\newcommand{\colbw}{\col{black}{white}}
\newcommand{\colww}{\col{white}{white}}

\newcommand{\xgap}{1}
\newcommand{\lw}{0.5}
\newcommand{\rad}{0.3}
\newcommand{\mazetwo}[2]{%
\begin{tikzpicture}[scale=0.07]
\draw[color=black, line width=\lw] (0,0) -- (0,1) -- (3,1) -- (3,2) -- (0,2) -- (0,3) -- (3,3);
\filldraw[draw=black, fill=white]  (3-#2,3-#1) circle (\rad);
\end{tikzpicture}
}

\newcommand{\mazethree}[3]{%
\begin{tikzpicture}[scale=0.07]
\draw[color=black, line width=\lw] (0,0) -- (0,1) -- (3,1) -- (3,2) -- (0,2) -- (0,3) -- (3,3);
\filldraw[draw=black, fill=white]  (3-#2,3-#1) circle (\rad);
\draw[color=black, line width=\lw] (3+\xgap,0) -- (3+\xgap,1) -- (3+\xgap+3,1) -- (3+\xgap+3,2) -- (3+\xgap,2) -- (3+\xgap,3) -- (3+\xgap+3,3);
\filldraw[draw=black, fill=white] (3+\xgap+3-#3,3-#2) circle (\rad);
\end{tikzpicture}}

\newcommand{\mazechiral}[3]{
\begin{tikzpicture}[scale=0.07]
\draw[color=black, line width=\lw] (0,0) -- (1,0) -- (1,3) -- (2,3) -- (2,0) -- (3,0) -- (3,3);
\filldraw[draw=black, fill=white]  (3-#1,3-#2) circle (\rad);
\draw[color=black, line width=\lw] (3+\xgap,0) -- (3+\xgap,1) -- (3+\xgap+3,1) -- (3+\xgap+3,2) -- (3+\xgap,2) -- (3+\xgap,3) -- (3+\xgap+3,3);
\filldraw[draw=black, fill=white] (3+\xgap+3-#3,3-#2) circle (\rad);
\end{tikzpicture}}

\newcommand{\mazefour}[4]{%
\begin{tikzpicture}[scale=0.07]
\draw[color=black, line width=\lw] (0,0) -- (0,1) -- (3,1) -- (3,2) -- (0,2) -- (0,3) -- (3,3);
\filldraw[draw=black, fill=white]  (3-#2,3-#1) circle (\rad);
\draw[color=black, line width=\lw] (3+\xgap,0) -- (3+\xgap,1) -- (3+\xgap+3,1) -- (3+\xgap+3,2) -- (3+\xgap,2) -- (3+\xgap,3) -- (3+\xgap+3,3);
\filldraw[draw=black, fill=white] (3+\xgap+3-#3,3-#2) circle (\rad);
\draw[color=black, line width=\lw] (6+2*\xgap,0) -- (6+2*\xgap,1) -- (6+2*\xgap+3,1) -- (6+2*\xgap+3,2) -- (6+2*\xgap,2) -- (6+2*\xgap,3) -- (6+2*\xgap+3,3);
\filldraw[draw=black, fill=white] (6+2*\xgap+3-#4,3-#3) circle (\rad);
\end{tikzpicture}}

\newcommand{\dial}[1]{%
\pgfmathsetmacro{\a}{20*#1-30}%
\rotatebox[origin=c]{\a}{$\rightarrow$}}

\MakeRobust\mazetwo
\MakeRobust\mazethree
\MakeRobust\mazechiral
\MakeRobust\mazefour
\MakeRobust\dial

\author[M. Goertz]{Madeleine Goertz}
\address{Department of Mathematics, California Polytechnic State University, San Luis Obispo, California, United States}
\email{\href{mailto:mgoertz@calpoly.edu}{mgoertz@calpoly.edu}}
\thanks{The authors would like to thank Williams College and Steven Miller for organizing the SMALL 2024 REU. The authors received financial support from NSF Grant DMS-2241623.}

\author[A. Williams]{Aaron Williams}
\address{Department of Computer Science, Williams College, Williamstown, Massachusetts, United States}
\email{\href{mailto:aaron.williams@williams.edu}{aaron.williams@williams.edu}}

\subjclass[2020]{05A05}

\keywords{Puzzles, Loopless Algorithm, Ziggu, Zigguflat, Nurikabe}

\title[The Ziggu family of puzzles]{The Ziggu family of puzzles: \\ An Algorithmic Solution \& a Combinatorial Interpretation with Nurikabe}
\title{The curious properties of the Quaternary Gray Code and how it can be used to solve Ziggu Puzzles}
\title[The Quaternary Gray Code and Ziggu Puzzles]{The Quaternary Gray Code and how it can be used to solve Ziggurat and other Ziggu Puzzles}

\begin{document}

\begin{abstract}

We investigate solutions to the new ``Ziggu'' family of exponential puzzles.
These puzzles have $p$ pieces that form $m$ mazes. 
We encode the state of each puzzle as an quaternary number (i.e., base $4$) with $n=m+1$ digits, where each digit gives the horizontal or vertical position in one maze. 
For example, the commercial version of \emph{Zigguflat} has $p=6$ pieces connected into $m=4$ mazes and its state requires $n=5$ digits to describe.
We show that the number of states on a shortest solution is $6 \cdot 2^n - 3n - 5$ (\OEIS{A101946}).
There is only one solution of this length, and it is generated from the start configuration by a simple algorithm: 
make the leftmost modification that doesn't undo the previous modification.
Replacing ``leftmost'' with ``rightmost'' instead generates the unique longest solution that visits all $(3^{n+1} - 1)/2$ states (\OEIS{A003462}).
In this way, Ziggu puzzles can be viewed as $4$-ary, $3$-ary, or $2$-ary puzzles based on how the number of state encodings, valid states, or minimum states grow with each additional maze.

Classic Gray code puzzles (e.g., Spin-Out) provide natural and illuminating comparisons.
Gray code puzzles with $p$ pieces have $2^p$ (\OEIS{A000079}) or $\lfloor \frac{2}{3} \cdot 2^p \rfloor$ (\OEIS{A000975} \cite{stockmeyer2017exploration}) states on their unique solution, and at most one modification doesn't undo the previous modification.
The states visited in a Gray code puzzle solution follow the well-known binary reflected Gray code.
We show that Ziggu puzzles instead follow the quaternary reflected Gray code.
More specifically, the shortest and longest solutions are both sublists of this order, meaning that some numbers are skipped over but the relative order of the remaining numbers does not change.
The sequence of piece changes in a Gray code puzzle solution is the binary ruler sequence.
This sequence has a familiar ``metronome'' property wherein every second change involves the first piece: 
$1,2,1,3,1,2,1,4,\ldots$ (\OEIS{A001511}).
In contrast, the shortest Ziggu solution has a ``stair-climbing'' property wherein successive changes differ by at most one:
$1,1,1,2,1,1,1,2,1,1,1,2,3,2,1,\ldots$.

These results show how to solve Ziggu puzzles from the start configuration.
To help solve the puzzle from an arbitrary configuration we provide $O(n)$-time ranking, comparison, and successor algorithms, which give the state's position along a solution, the relative order of two states, and the next state, respectively.
While Gray code puzzles have simpler recursive descriptions and successor rules, the Ziggu puzzle has a much simpler loopless algorithm to generate its shortest solution than the Gray code puzzles do. 
The two families are also intimately related as they have the same comparison function.
Finally, we enrich the literature on sequence \OEIS{A101946} by providing a bijection between Ziggu states and valid $2$-by-$n$ Nurikabe grids.
\end{abstract}

\maketitle

\newpage
\setcounter{tocdepth}{1} 
\tableofcontents

\section{Introduction}
\label{sec:intro}

\begin{figure}
    \centering
    \includegraphics[width=1.0\linewidth]{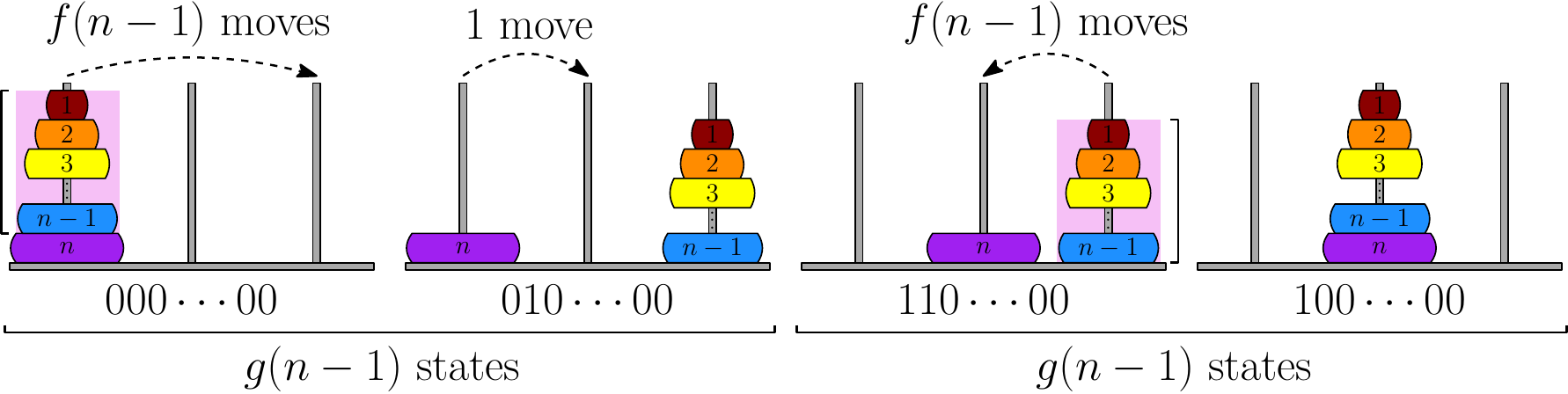}
    \caption{The \emph{Towers of Hanoi} starts with all discs on the left peg, and ends with all discs on the middle peg.
    A single move consists of picking up the top disc on one peg and putting it down on another peg, with the condition that a larger disc can never be placed on top of a smaller disc.
    It is solved by recursively moving the top $n-1$ discs to third peg, moving disc $n$ to the middle peg, then recursively moving the smallest $n-1$ discs to the middle peg.
    The total number of moves is $2^n-1 $ by the recurrence $f(n) = 2f(n-1) + 1$ with $f(1)=1$, while the total number of states is $2^n$ by the recurrence $g(n) = 2 \cdot g(n-1)$ with $g(1)=2$.}
    \label{fig:recurrenceTowers}
\end{figure}

The mathematics of puzzle solving is often categorized as recreational, but it also provides foundational educational experiences.
For example, many students in Computer Science and Mathematics first learn about recursion, induction, and exponential growth through the \emph{Towers of Hanoi} puzzle \cite{hinz2018tower}.
More specifically, the standard approach of solving the puzzle involves recursively moving the smallest $n-1$ discs, followed by moving disc $n$ once, followed by recursively moving the smallest $n-1$ discs again, as illustrated in Figure \ref{fig:recurrenceTowers}.
The number of moves in this solution is given by the recurrence in \eqref{eq:recurrenceTowers}.
\begin{equation} \label{eq:recurrenceTowers}
    T(n) = 2 \cdot T(n-1) + 1 \text{ with } T(1) = 1.
\end{equation}
Students can then use induction to prove that $2^n-1$ is a closed formula for \eqref{eq:recurrenceTowers}.
This makes it an \emph{exponential puzzle} since the shortest solution grows exponentially with each successive piece.
More specifically, it is a \emph{$2$-ary puzzle} as it grows like~$\theta(2^n)$ (where ``big-theta'' is the equality analogue of ``big-O's'' less than or equal).

While these conclusions often mark the end of a student's investigation of the \emph{Towers of Hanoi}, the puzzle can lead to much deeper discoveries in combinatorics and computation.
For example, students may observe that the puzzle has $3^n$ states, but only $2^n$ of them are used in the shortest solution.
There is a natural, but not immediately obvious, approach to encoding these states as an $n$-bit binary string.
This leads to the well-known order of binary strings known as the \emph{binary reflected Gray code} \cite{Fra1953} and its rich generalizations \cite{Mut2023}.
Similarly, the sequence of discs that move in the shortest solution follow the \emph{binary ruler sequence} (\OEIS{A001511}) or more precisely the \emph{paper folding sequence} (\OEIS{A164677}).
Given an arbitrary state of the puzzle, students may wonder if they can directly determine the next state (i.e., without any additional context), or how many steps it will take to reach the solution.
These questions lead to the development of the successor rules and ranking algorithms for the Gray code and their coverage in \emph{The Art of Computer Programming} \cite{knuth2011art}.
In other words, the development of \eqref{eq:recurrenceTowers} and its closed formula are just the first lessons provided by \emph{Towers of Hanoi}.
Furthermore, a wide variety of other exponential puzzles have been constructed over the years that follow a similar recurrence or associated sequences.
This family of puzzles is sometimes known as \emph{Gray code puzzles}.

This paper considers the new ``Ziggu'' family of exponential puzzles and various facets involved in solving them.
As we will see, the shortest solution follows a similar but more complicated recurrence.
\begin{equation} \label{eq:recurrenceZiggu}
    Z(n) = 3 \cdot Z(n-1) - 2 \cdot Z(n-2) + 3 \text{ with } Z(1) = 4 \text{ and } Z(2) = 13.
\end{equation}
Roughly speaking, solving a puzzle with $m=n-1$ Ziggu mazes requires solving its first $m-1$ Ziggu mazes three times, however, the solution to its first $m-2$ Ziggu mazes can be skipped in the latter two sub-solutions; see Figure \ref{fig:mazeMoves_short}.
The recurrence in \eqref{eq:recurrenceZiggu} has a closed form of $Z(n) = 6 \cdot 2^n - 3n - 5$ (\OEIS{A101946}), whose dominant $2^n$ term makes it a \emph{$2$-ary exponential puzzle} in terms of its shortest solution.
Just like the \emph{Towers of Hanoi}, we'll see that the puzzles can also be viewed as $3$-ary puzzles as they have $(3^{n+1} - 1) / 2$ (\OEIS{A003462}) total states.
In fact, we'll also view them as $4$-ary puzzles due to an intimate connection with the quaternary (base-$4$) reflected Gray code.

The remainder of this introductory section discusses the Ziggu and Gray code families of puzzles.
Then we highlight our new results.

\subsection{Ziggu Puzzles: Ziggurat, Zigguflat, Zigguhooked ...}
\label{sec:intro_Ziggu}

\begin{figure}
    \centering
    \begin{subfigure}{0.24\textwidth}
        \includegraphics[width=0.95\textwidth]{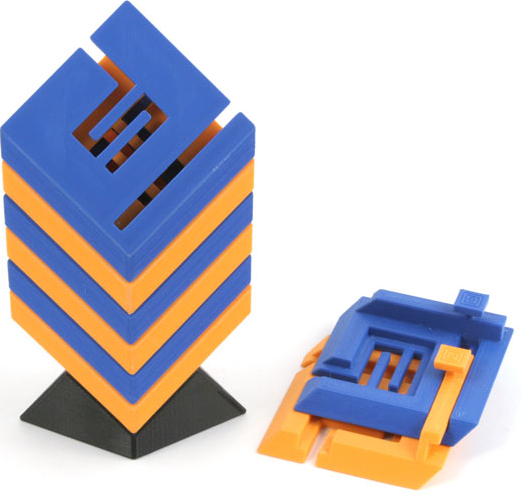}
        \caption{Ziggurat}
        \label{fig:ziggus_ziggurat}
    \end{subfigure}
    \begin{subfigure}{0.24\textwidth}
        \rotatebox{180}{\includegraphics[width=0.95\textwidth]{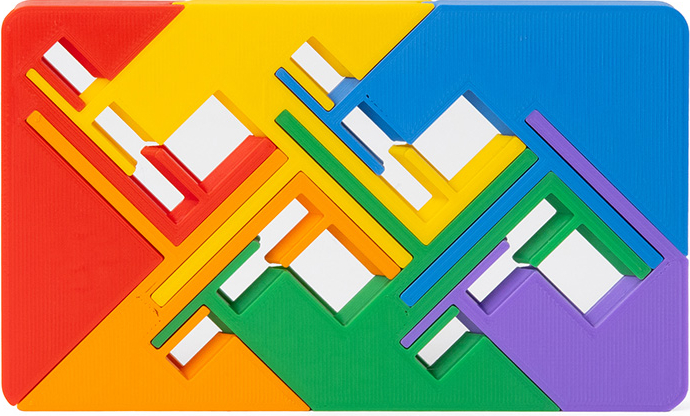}}
        \caption{Zigguflat}
        \label{fig:ziggus_zigguflat}
    \end{subfigure}
    \begin{subfigure}{0.24\textwidth}
        \includegraphics[width=0.95\textwidth]{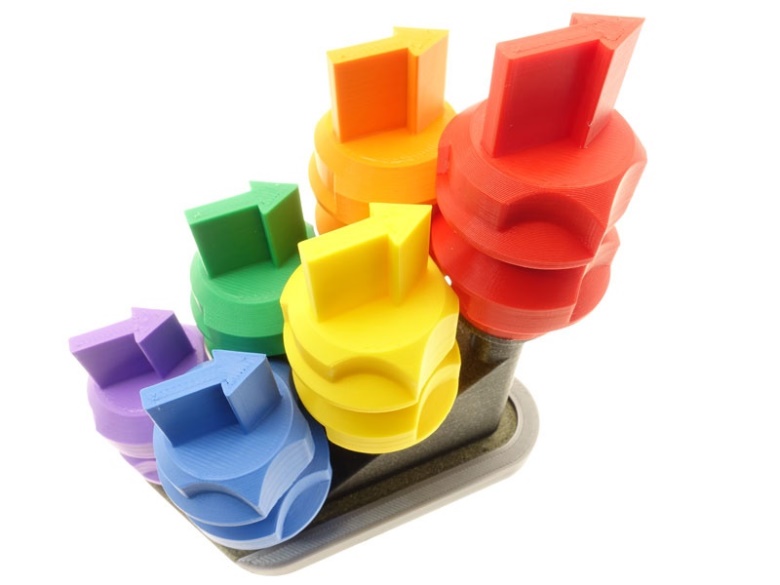}
        \caption{Ziggutrees}
        \label{fig:ziggus_ziggutrees}
    \end{subfigure}
    \begin{subfigure}{0.24\textwidth}
        \rotatebox{180}{\includegraphics[width=0.95\textwidth]{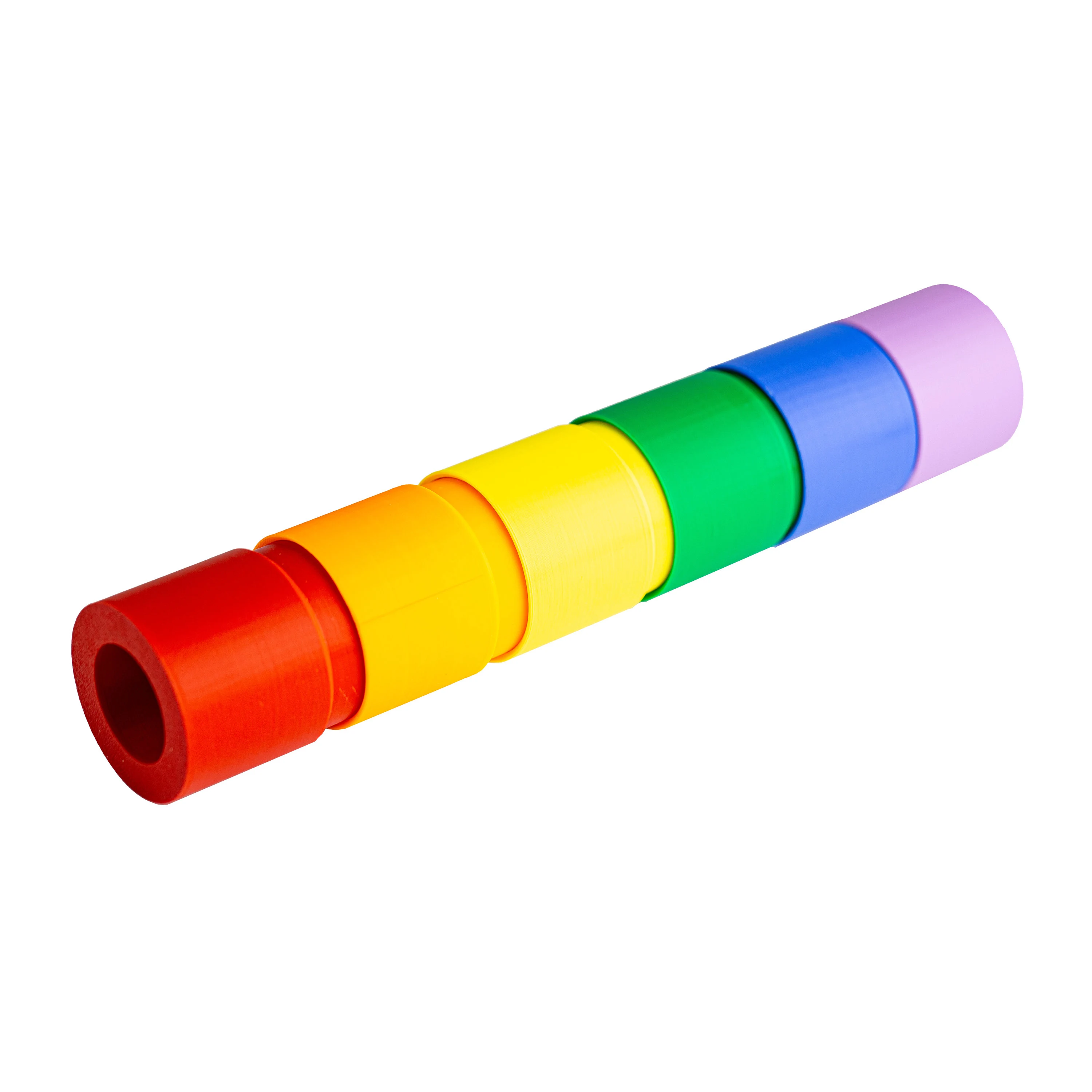}}
        \caption{Ziggustretch}
        \label{fig:ziggus_ziggustretch}
    \end{subfigure}
    \begin{subfigure}{0.24\textwidth}
        \reflectbox{\includegraphics[width=0.95\textwidth]{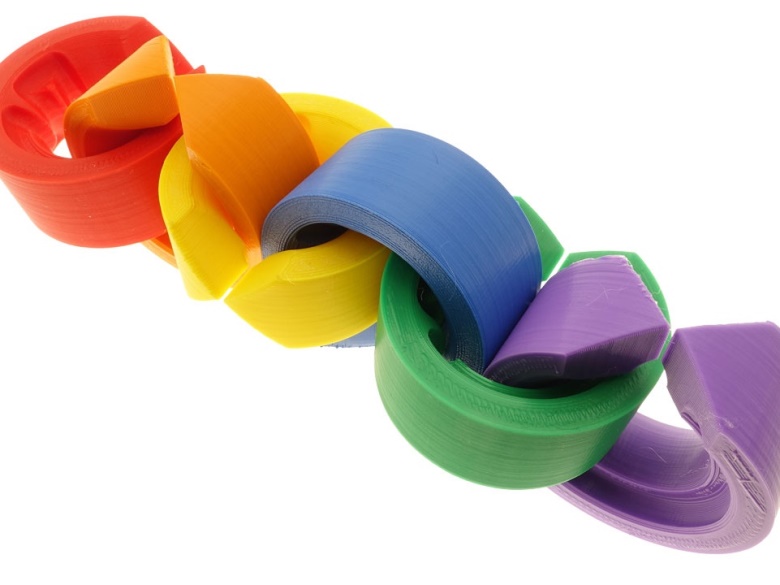}}
        \caption{Zigguchain}
        \label{fig:ziggus_zigguchain}
    \end{subfigure}
    \begin{subfigure}{0.24\textwidth}
        \includegraphics[width=0.95\textwidth]{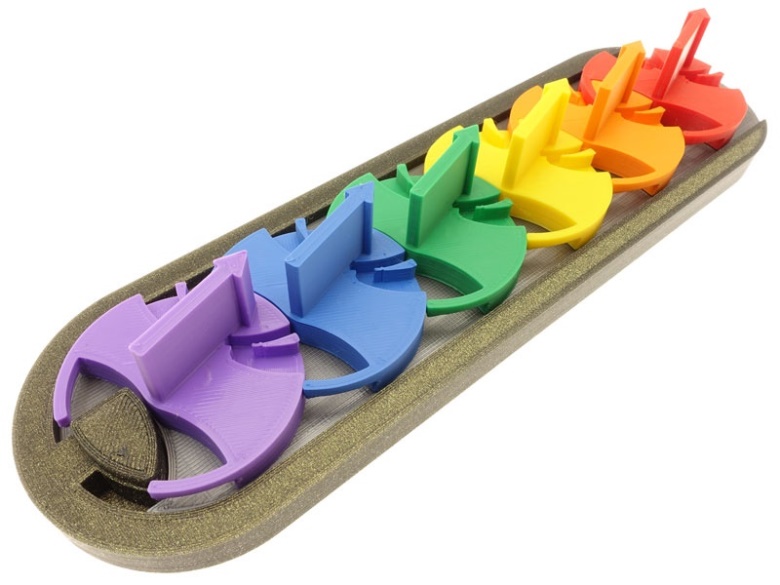}
        \caption{Zigguhooked}
        \label{fig:ziggus_zigguhooked}
    \end{subfigure}
    \begin{subfigure}{0.24\textwidth}
        \reflectbox{\includegraphics[width=0.95\textwidth]{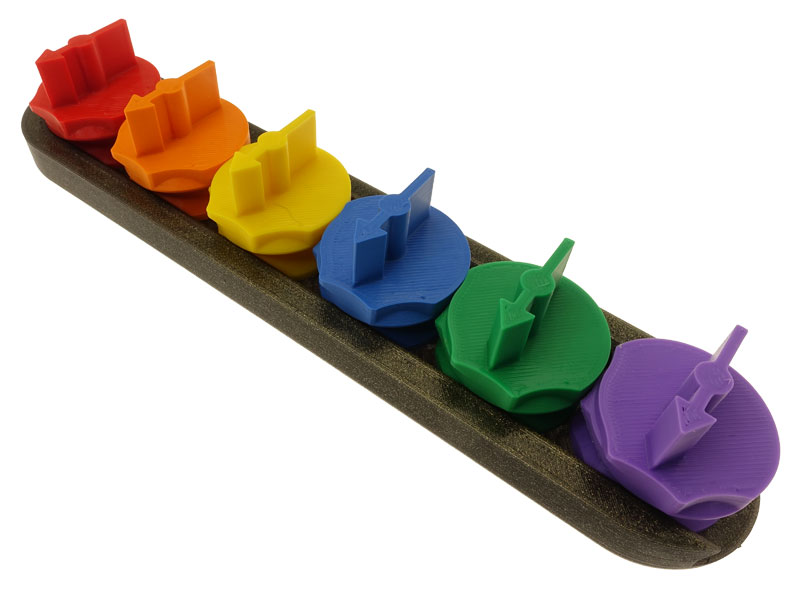}}
        \caption{Ziggutwist}
        \label{fig:ziggus_ziggutwist}
    \end{subfigure}
    \begin{subfigure}{0.24\textwidth}
        \reflectbox{\includegraphics[width=0.95\textwidth]{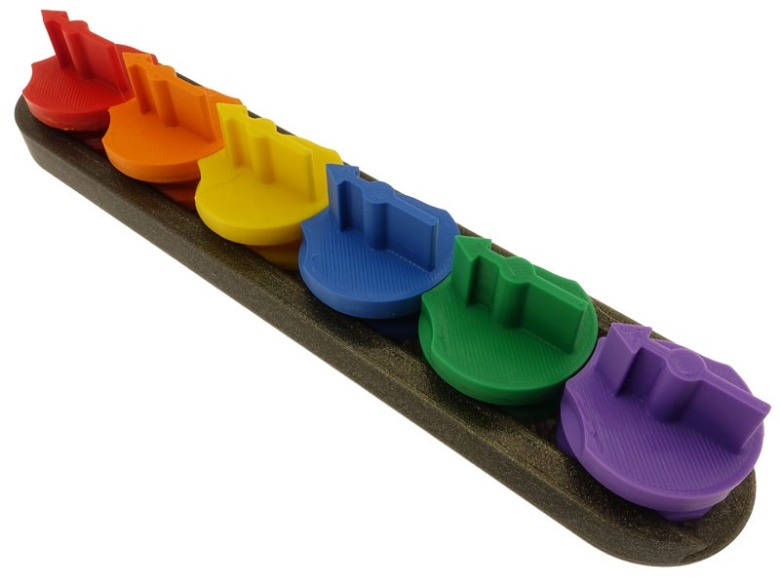}}
        \caption{Ziggutwist Too}
        \label{fig:ziggus_ziggutwist_too}
    \end{subfigure}
    \caption{The Ziggu family of puzzles.
    Each puzzle has $p=6$ pieces in alternating or rainbow colors.
    Some puzzles also have a black stand or frame which does not count as a piece.
    }
    \label{fig:ziggus}
%
\vspace{1.5em}
\small
\begin{tabular}{|c|c|c|c|c|c|c|c|c|} 
\hline
\rowcolor{lightgray}
puzzles & type & \makecell{pieces \\[-0.1em] $p$} & \makecell{mazes \\[-0.1em] $m$} & \makecell{digits \\[-0.1em] $n$} & \makecell{shortest \\[-0.1em] solution} & \makecell{longest \\[-0.1em] solution} \\ 
\hline
\makecell{Ziggurat, Zigguflat, \\ Ziggutrees, Ziggustretch} &
4-neighbor & 6 & 4 & 5 & 172 & 364 \\ \hline
\makecell{Zigguchain, Zigguhooked, \\ Ziggutwist, Ziggutwist Too} &
2-neighbor & 6 & 5 & 6 & 361 & 1093\\
\hline
\end{tabular}
\captionof{table}{There are two types of commercially available Ziggu puzzles.
Each piece either interacts with the previous and next pieces ($2$-neighbor) or the previous two and next two pieces ($4$-neighbor).
This affects how many mazes that must be solved, the number of quaternary digits we use to represent an individual state, and the number of states in their shortest and longest solutions.
}
\label{tab:ziggus}
%
\vspace{1.5em}
\small
\begin{tabular}{|c|c|c|c|l|c|} 
\hline
\rowcolor{lightgray}
strings & set & cardinality & growth & \multicolumn{1}{c}{values} & sequence \\ 
\hline
encoding & $\QUAT{n}$ & $4^n$ & $\theta(4^n)$  & $4, 16, 64, 256,1024,\mathbf{4096},\ldots$ & \OEIS{A000302} \\
valid state & $\LONG{n}$ & $(3^{n+1} - 1)/2$ & $\theta(3^n)$ & $4, 13, 40, 121, 364, \mathbf{1093},\ldots$ & \OEIS{A003462} \\
ziggu & $\SHORT{n}$ & $6 \cdot 2^n-3n-5$ & $\theta(2^n)$ & $4, 13, 34, 79, 172, \mathbf{361}, \ldots$ & \OEIS{A101946} \\
\hline
\end{tabular}
\captionof{table}{Integer sequences related to Ziggu puzzles encoded by $n$ digits.
For example, Zigguhooked with $p=6$ dials is encoded by $n=6$ digits, so there are $4^6 = 4096$ state encodings.
Only $1093$ of the encodings represent valid states and $361$ are on a shortest solution.
In particular, $132222 \in \QUAT{6}$ is an encoding, but it is not a valid state as the dials cannot be physically rotated as \dial{1}\dial{3}\dial{2}\dial{2}\dial{2}\dial{2}.
So this encoding only contributes to the top row's quantity as $132222 \notin \LONG{6}$ and $132222 \notin \SHORT{6}$. 
The growths of these quantities also explain how Ziggu puzzles can be viewed as quaternary ($4$-ary), ternary ($3$-ary), or binary ($2$-ary) puzzles.
}
\label{tab:sequences}
\end{figure}

The \emph{Ziggurat} puzzle was created by Eitan Cher and Bram Cohen and submitted to the \emph{2021 Puzzle Design Competition} \cite{zotero-1627}.
The puzzle is notable for several reasons.
It is an \emph{exponential puzzle} since it consists of some number of (nearly) identical 
pieces and the number of steps to required to solve the puzzle grows exponentially in some base.
Unlike previous puzzles of this type it is a \emph{burr} meaning that the goal is to fully disassemble (or reassemble) it.
It is also \emph{frameless} meaning that the puzzle consists only of its pieces without any other structure holding it together%
\footnote{Note that the black base in Figure \ref{fig:ziggus_ziggurat} is simply a stand and is not part of the puzzle.}. 
The name of the puzzle is based on how the pieces are initially stacked into a 3-dimensional pyramidal configuration resembling a ziggurat temple. The underlying maze mechanic underpinning the puzzle is discussed in Section \ref{sec:mazes}.

Oskar van Deventer implemented the same underlying maze structure in a variety of creative ways in the design of subsequent exponential puzzles including \emph{Zigguflat}, \emph{Ziggutwist}, \emph{Ziggutwist Too}, \emph{Zigguchain}, \emph{Zigguhooked}, \emph{Ziggutrees}, and \emph{Ziggustretch}.
All of these puzzles are seen in Figure \ref{fig:ziggus} with $p=6$ pieces.
In some of these puzzles the mazes can be clearly seen (e.g., Zigguflat) while in others they are partially obscured (e.g., Ziggutwist) or fully obscured (e.g., Ziggustretch). 
Ziggu puzzles are typically color-coded so that the red piece is the only one that can be modified in the starting configuration.
Right-handed solvers often hold the puzzle with this piece on the right.
We'll assume this orientation in our discussions, and we'll also refer to the red piece as the first piece, the orange piece as the second piece, and so on.
(In \emph{Ziggurat} the bottom piece is the first piece.)

While each Ziggu puzzle is unique, they fall into two primary categories, summarized in Table \ref{tab:ziggus}. Ziggurat, Zigguflat, Ziggutrees, and Ziggustretch (Figures \ref{fig:ziggus_ziggurat}-\ref{fig:ziggus_ziggustretch}) are all \textit{4-neighbor} puzzles, meaning that each piece touches at most four other pieces. In particular, Figure \ref{fig:ziggus_zigguflat} shows that the yellow piece touches the orange and red pieces in one direction, as well as the blue and green pieces in the other direction. On the other hand, Zigguchain, Zigguhooked, Ziggutwist, and Ziggutwist Too are all \textit{2-neighbor} puzzles, as each piece touches up to two other pieces. This is easily seen in Figures \ref{fig:ziggus_zigguchain}-\ref{fig:ziggus_ziggutwist_too}, where the pieces are arranged linearly and each piece touches its left and right neighbor. 

In the remainder of this paper, we consider Zigguflat and Zigguhooked as examples of 4- and 2-neighbor puzzles. These have the advantage of using distinctly different implementations of the Ziggu maze while still being planar, which makes the mechanics visible. As noted in Table \ref{tab:ziggus}, while both Zigguflat and Zigguhooked consist of $p=6$ pieces, Zigguhooked will take longer to solve since it has $5$ mazes, while Zigguflat has only $4$ mazes.

Despite these seemingly large differences, all of the Ziggu puzzles are identical in an abstract sense.
More specifically, all Ziggu puzzles with $m$ mazes have identical state spaces and solutions.

\subsection{Gray Code Puzzles: Chinese Rings, Towers of Hanoi, Spin-Out ...}
\label{sec:intro_Gray}

\begin{figure}
    \centering
    \begin{subfigure}{0.24\textwidth}
        \centering
        \includegraphics[width=0.95\textwidth]{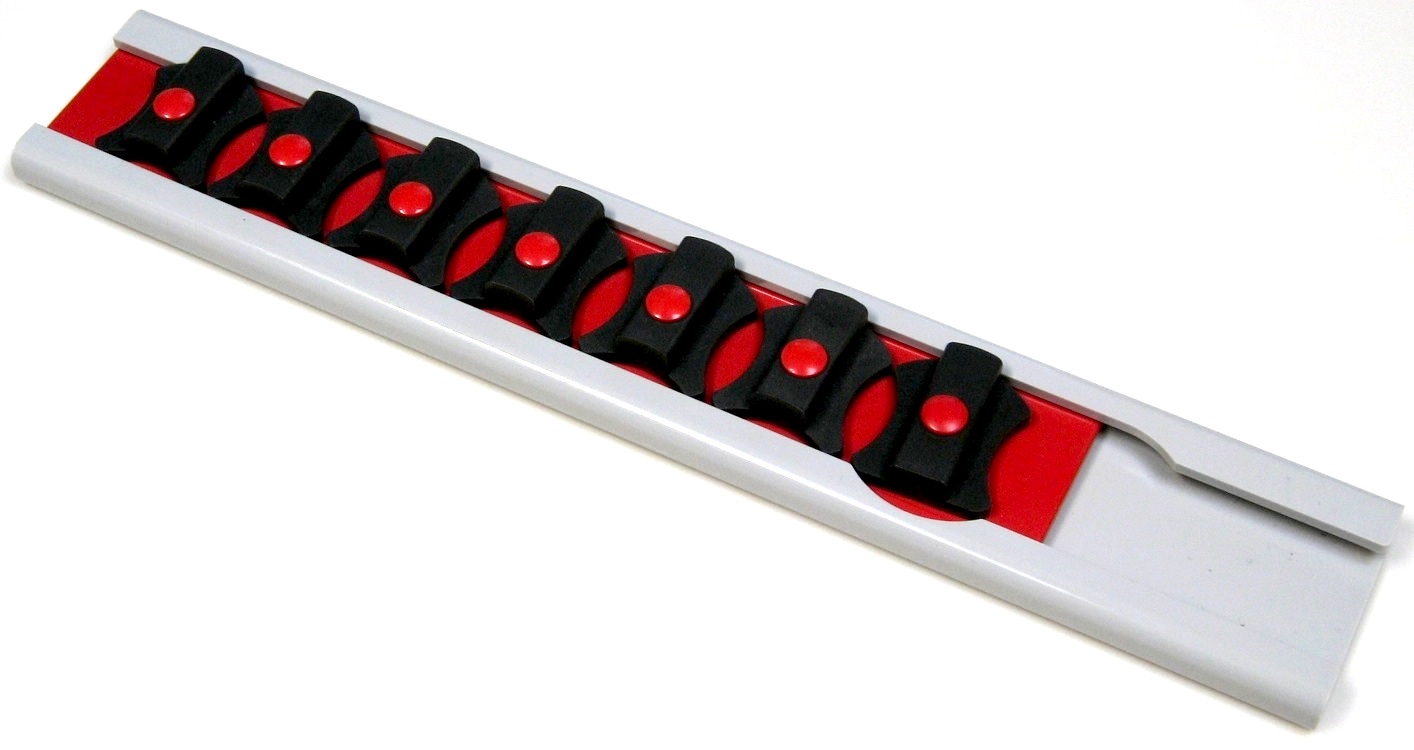}
        \caption{Spin-Out \cite{Sto}}
        \label{fig:grays_spinout}
    \end{subfigure}
    \hfill
    \begin{subfigure}{0.24\textwidth}
        \centering
        \includegraphics[width=0.95\textwidth]{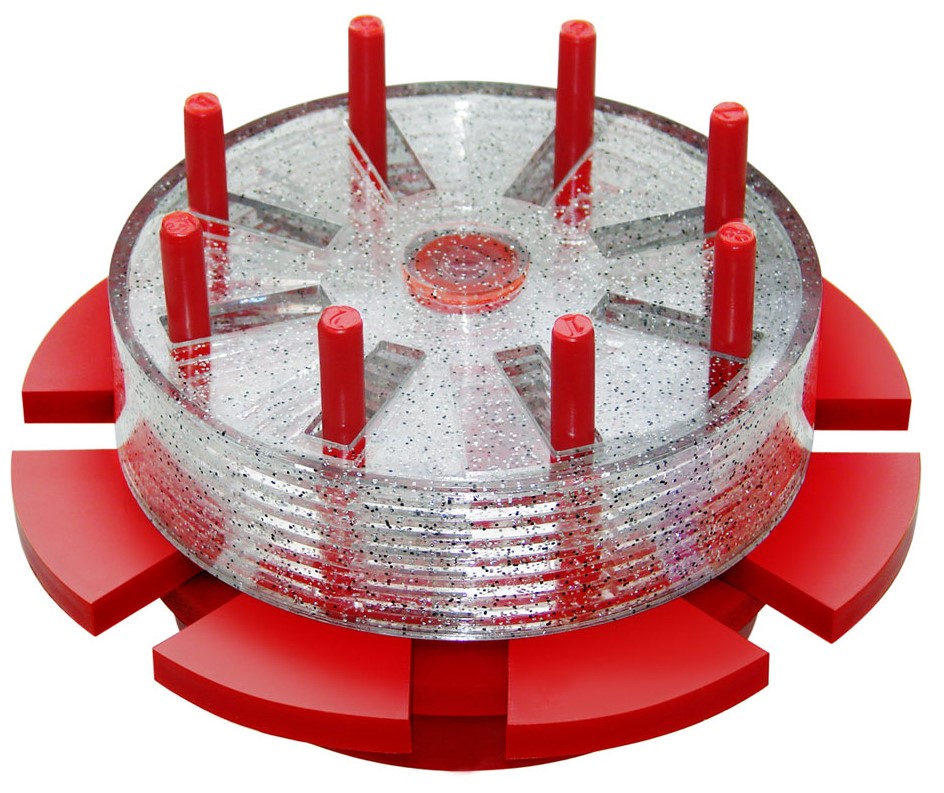}
        \caption{The Brain}
        \label{fig:grays_brain}
    \end{subfigure}
    \hfill
    \begin{subfigure}{0.24\textwidth}
        \centering
        \includegraphics[width=0.85\textwidth]{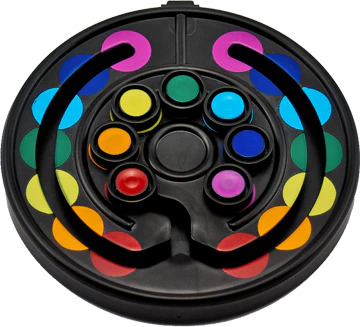}
        \caption{Rudenko's Disk}
        \label{fig:grays_rudenko}
    \end{subfigure}
    \hfill
    \begin{subfigure}{0.24\textwidth}
        \centering
        \includegraphics[width=0.85\textwidth]{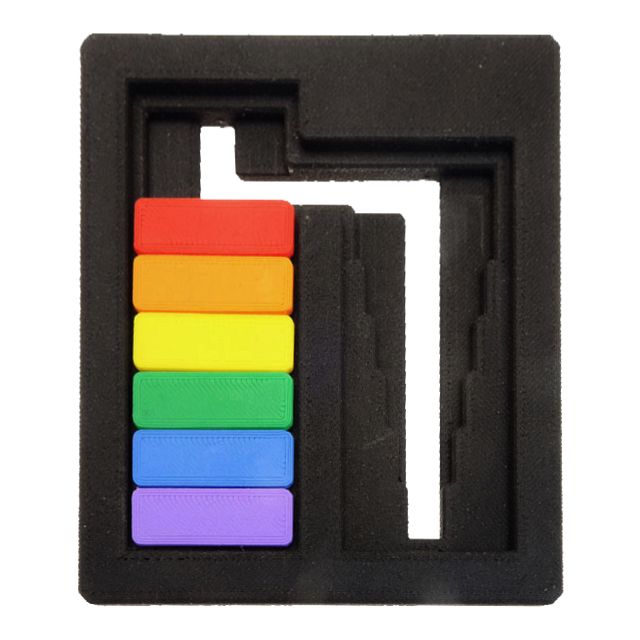}
        \caption{Panex Jr.}
        \label{fig:grays_panexjr}
    \end{subfigure}
    \begin{subfigure}{0.24\textwidth}
        \centering
        \includegraphics[width=0.95\textwidth]{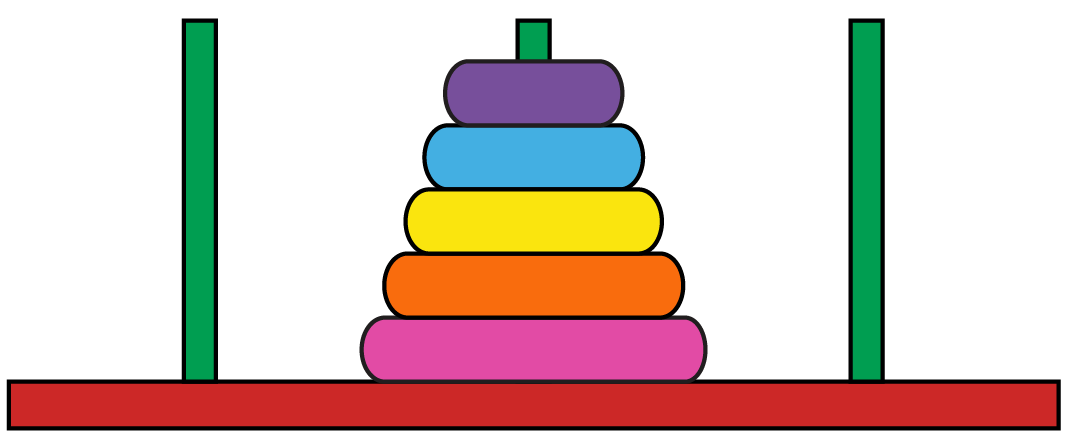}
        \caption{Towers of Hanoi}
        \label{fig:grays_hanoi}
    \end{subfigure}
    \hfill
    \begin{subfigure}{0.24\textwidth}
        \centering
        \includegraphics[width=0.95\textwidth]{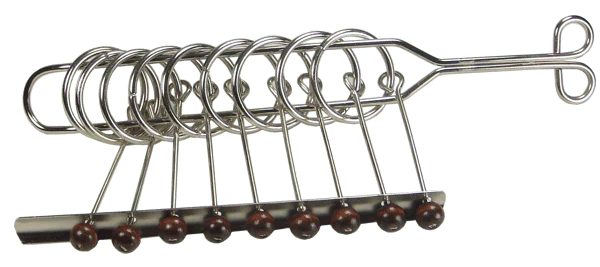}
        \caption{Baguenaudier}
        \label{fig:grays_baguenaudier}
    \end{subfigure}
    \hfill
    \begin{subfigure}{0.255\textwidth}
        \centering
        \includegraphics[width=0.95\textwidth]{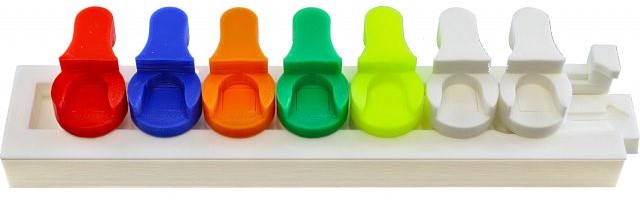}
        \caption{{\scriptsize Crazy Elephant Dance}}
        \label{fig:grays_elephant}
    \end{subfigure}
    \hfill
    \begin{subfigure}{0.24\textwidth}
        \centering
        \includegraphics[width=0.95\textwidth]{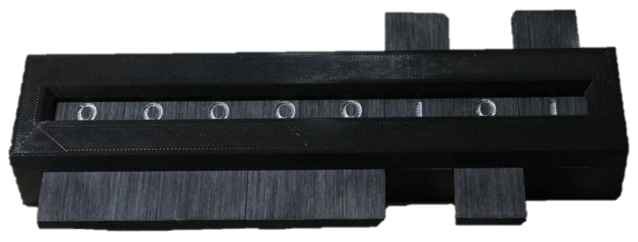}
        \caption{B-Bar}
        \label{fig:grays_bbar}
    \end{subfigure}
    \caption{The Gray code family of puzzles.}
    \label{fig:GrayCodePuzzles}
\vspace{1.5em}
\small
\begin{tabular}{|c|c|c|c|c|c|c|c|} 
\hline
\rowcolor{lightgray}
puzzle & base & \Gape[0pt][2pt]{\makecell{pieces \\ $n$}} & \Gape[0pt][2pt]{\makecell{start \\ state}} & \Gape[0pt][2pt]{\makecell{target \\ state}} & \Gape[0pt][2pt]{\makecell{shortest \\ solution}} \\ \hline
Spin-Out & binary & $7$ & $1111111$ & $0000000$ & $86$ \\ \hline
The Brain & binary & $8$ & $00000000$ & $11111111$ & $171$  \\ \hline
B-Bar & binary & $8$ & $00000000$ & $10000000$ & $256$  \\ \hline
Baguenaudier & binary & varies & $1^n$ & $0^n$ & \OEIS{A005578}  \\ \hline
Towers of Hanoi & \makecell{binary \\ ternary} & varies & $0^n$ & $10^{n-1}$ & $2^n$  \\ \hline
\end{tabular}
\captionof{table}{The number of pieces $n$ in various Gray code puzzles, along with the number of states in the unique (shortest) solutions. 
Towers of Hanoi is ternary ($3$-ary) in terms of the total number of states, but only binary ($2$-ary) in terms of its shortest solution.}
\label{tab:Grays}
\vspace{1.5em}
    \makebox[\textwidth][c]{\includegraphics[width=1.2\textwidth]{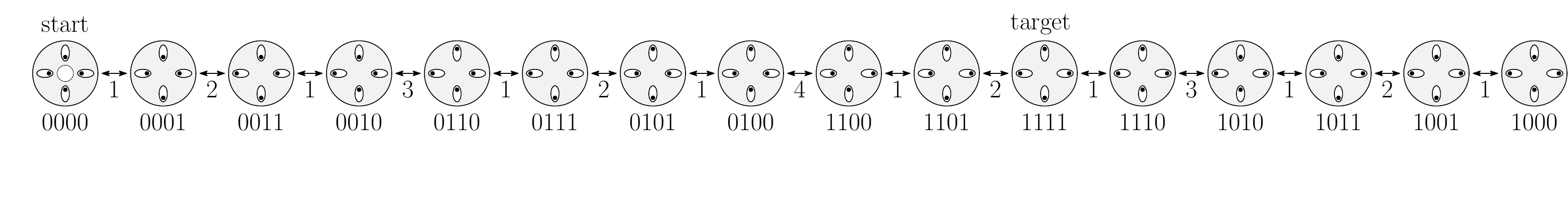}}
    \centering
    \caption{
    Most Gray code puzzles with $n$ pieces have the same state space.
    It is a path between the first and last strings in the binary reflected Gray code.
    This is shown above for \emph{The Brain} with $n=4$ wedges.
    Each bit in the state indicates whether a given wedge is pushed in or pulled out.
    The edge labels in the graph provide the wedge to toggle (i.e., the bit to complement in $b_4 b_3 b_2 b_1$) to move to a neighboring state.
    Note that the change sequence from $0000$ to $1111$ is $1,2,1,3,\ldots$ while the change sequence from $1111$ to $0000$ is $2,1,4,1,\ldots$.
    So it is important to differentiate between the start and target states (as in Table \ref{tab:Grays}).
    In particular, \emph{Spin-Out} and \emph{The Brain} puzzles with $n$ pieces have different solutions of the same length; the latter puzzle is easier as there is no choice for the first move.
    }
    \label{fig:statesBrain4}
\end{figure}

The most famous family of exponential puzzles is the \emph{Gray code puzzle} family, which are related to a well-known pattern of binary numbers known as the \emph{binary reflected Gray code (BRGC)} or simply Gray code.
A sampling appears in Figure \ref{fig:GrayCodePuzzles}.

The family consists of vintage puzzles including \emph{Spin-Out} (see \cite{pruhs1993spin} and \cite{cull2013tale}) and \emph{The Brain}.
In these puzzles the underlying state space is simply a path labeled by the Gray code.
The puzzles start at state $1^n$ and are solved when reaching state $0^n$, where $n$ is the number of pieces and exponentiation denotes repetition.
The state space of \emph{The Brain} with $n=4$ pieces appears in Figure \ref{fig:statesBrain4}.

Most famously, the family contains historical puzzles like the \emph{Towers of Hanoi} where the Gray code provides the shortest solution of length $2^n$ among its $3^n$ states, and \emph{Baguenaudier} (or ``Chinese Rings'') \cite{rougetet2024binaire} where the state space is the Gray code except that some pairs of moves can be made simultaneously.

The Gray code puzzle family also includes more recent puzzles like \emph{Rudenko's Disc} where the Gray code provides a non-shortest solution, and \emph{B-Bar} where the only solution is the Gray code.
More specifically, the solution to B-Bar is the unique Hamilton path in the Gray code state space.
More broadly, Hamilton paths have been used as the basis for other sequential puzzles, including Hamilton's \emph{Icosian game} \cite{pegg2009icosian} \cite{mutze2023hamilton}.
Ziggu puzzles have a more complicated state space (see Figure \ref{fig:statesZiggu}) but this state also has a unique Hamilton path.

Close relatives of the Gray code puzzle family include \emph{Crazy Elephant Dance} (see \cite{cooper2019chinese}) which follows a ternary (i.e., base-$3$) Gray code, and \emph{Panex Jr} which is a restricted version of \emph{Rudenko's Disc} that cannot be solved using the Gray code.

\subsection{New Results}
\label{sec:intro_new}

Those who have held a Ziggu puzzle may have stumbled upon the following approach to exploring its configurations:
modify the leftmost piece that does not undo the most recent modification.
Pleasantly, this approach will eventually solve the puzzle.
Moreover, one can argue that it does so in as few moves as possible.
However, it isn't immediately clear how many moves are made during the solution.
We prove that exactly $6 \cdot 2^n-3n-5$ states are on the shortest solution, and we refer to these states as \emph{Ziggu states} or strings.
Furthermore, if one replaces ``leftmost'' with ``rightmost'' then the exploratory approach again solves the puzzle.
Moreover, it does so in as many moves as possible.
In fact, this longest solution visits every one of the $(3^{n+1} - 1)/2$ states, and it is the only solution to do so.
Both of these formulae assume that state is encoded as one of the $4^n$ quaternary (i.e., base $4$) strings with $n$ digits, as explained later on.
These three formulae are summarized in Table \ref{tab:sequences} and they explain why Ziggu puzzles can be considered to be $4$-ary, $3$-ary, or $2$-ary exponential puzzles.


These results imply that Gray code puzzles are simpler than Ziggu puzzles.
This is because Gray code puzzles can be solved using a similar exploratory approach:
modify the only piece that does not undo the most recent modification.
In Gray code puzzles there is never a choice of which piece to modify (so long as the previous modification is known) and this approach gives the only solution.
(Most Gray code puzzles have some additional nuance that is not captured by these statements.)

\begin{figure}[ht]
    \begin{subfigure}{1.0\textwidth}
        \includegraphics[width=1.0\textwidth]{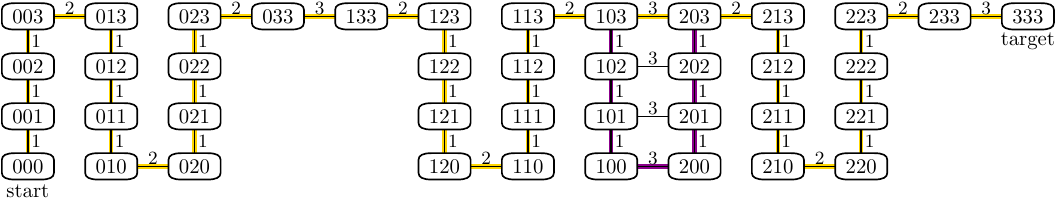}
        \caption{State graph for Ziggu puzzles encoded with $n=3$ digits.}
        \label{fig:statesZiggu_3}
    \end{subfigure}
    \begin{subfigure}{1.0\textwidth}
        \vspace{1.5em}
        \hspace{-1.5em}
        \includegraphics[width=1.11\textwidth]{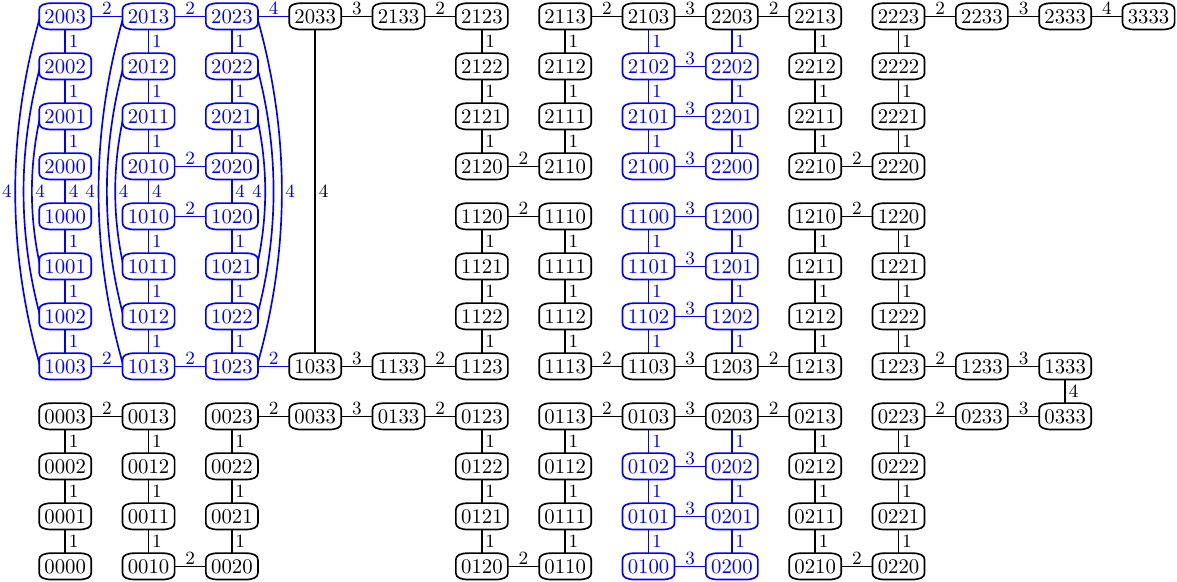}
        \caption{State graph for Ziggu puzzles encoded with $n=4$ digits.}
        \label{fig:statesZiggu_4}
    \end{subfigure}
    \caption{
    Ziggu puzzles with $n$ digits have the same state space.
    In (a) the gold edges form the unique shortest solution while the additional purple edges form the unique longest solution.
    In (b) the blue vertices appear on the longest solution but not the shortest solution.
    The lengths of the shortest and longest paths are (a) $34$ and $40$ and (b) $79$ and $81$, respectively, as in Table~\ref{tab:sequences}. 
    Note that (a) is a subgraph of (b) (with $0$ prefixed to each state label).
    }
    \label{fig:statesZiggu}
\end{figure}

We show that the connection between Ziggu puzzles and Gray code puzzles runs even deeper.
In Gray code puzzles each state can be represented by an $n$-bit binary string, and the next state is obtained by flipping the bit according to the binary reflected Gray code.
Along the same lines, the longest and shortest solutions to Ziggu puzzles can be understood using the quaternary (i.e., base $4$) reflected Gray code.
More specifically, the longest solution is obtained by skipping over invalid encodings in the quaternary reflected Gray code, and the shortest solution is obtained by additionally skipping over non-Ziggu states in the longest solution.
This is illustrated by the purple edges in the state space graph in Figure \ref{fig:statesZiggu}, and the gaps in Table \ref{tab:orders}. 
We also investigate the sequence of piece changes (rather than states) made during the longest and shortest solutions.
As we will see, there is again similarity to the analogous sequence for Gray code puzzles, which is known as the binary ruler sequence.

These results show how to solve Ziggu puzzles from the start configuration.
To help solve the puzzle from an arbitrary configuration we provide $O(n)$-time ranking, comparison, and successor algorithms, which give the state's position along a solution, the relative order of two states, and the next state, respectively.
Finally, we enrich the literature on sequence \OEIS{A101946} by providing a bijection between Ziggu states and valid $2$-by-$n$ Nurikabe grids. 

The title of this paper pays homage to an article by Martin Gardner \cite{gardner1972curious} (also see \cite{gardner2020knotted}) which discusses the Gray code and several of the aforementioned puzzles.

\begin{table}
\makeatletter
\setlength{\tabcolsep}{1.75pt}
\makeatother
\footnotesize
\centerline{\begin{tabular}{cccc|cccc|cccc|cc}
$\rank\Q$ & $\QUAT{3}$ & $\rulerQuat{3}$ & $\srulerQuat{3}$ & $\rank\V$ & $\LONG{3}$ & $\rulerLong{3}$ & $\srulerLong{3}$ & $\rank\S$ & $\SHORT{3}$ & $\rulerShort{3}$ & $\srulerShort{3} $& maze & dials\\
\hline
$ 0$ & $000$ & $1$ & $ 1$ & $ 0$ & $000$ & $1$ & $ 1$ & $ 0$ & $000$ & $1$ & $ 1$ & \mazethree{0}{0}{0} & \dial{0}\dial{0}\dial{0} \\[-0.2em]
$ 1$ & $001$ & $1$ & $ 1$ & $ 1$ & $001$ & $1$ & $ 1$ & $ 1$ & $001$ & $1$ & $ 1$ & \mazethree{0}{0}{1} & \dial{0}\dial{0}\dial{1} \\[-0.2em]
$ 2$ & $002$ & $1$ & $ 1$ & $ 2$ & $002$ & $1$ & $ 1$ & $ 2$ & $002$ & $1$ & $ 1$ & \mazethree{0}{0}{2} & \dial{0}\dial{0}\dial{2} \\[-0.2em]
$ 3$ & $003$ & $2$ & $ 2$ & $ 3$ & $003$ & $2$ & $ 2$ & $ 3$ & $003$ & $2$ & $ 2$ & \mazethree{0}{0}{3} & \dial{0}\dial{0}\dial{3} \\[-0.2em]
$ 4$ & $013$ & $1$ & $-1$ & $ 4$ & $013$ & $1$ & $-1$ & $ 4$ & $013$ & $1$ & $-1$ & \mazethree{0}{1}{3} & \dial{0}\dial{1}\dial{3} \\[-0.2em]
$ 5$ & $012$ & $1$ & $-1$ & $ 5$ & $012$ & $1$ & $-1$ & $ 5$ & $012$ & $1$ & $-1$ & \mazethree{0}{1}{2} & \dial{0}\dial{1}\dial{2} \\[-0.2em]
$ 6$ & $011$ & $1$ & $-1$ & $ 6$ & $011$ & $1$ & $-1$ & $ 6$ & $011$ & $1$ & $-1$ & \mazethree{0}{1}{1} & \dial{0}\dial{1}\dial{1} \\[-0.2em]
$ 7$ & $010$ & $2$ & $ 2$ & $ 7$ & $010$ & $2$ & $ 2$ & $ 7$ & $010$ & $2$ & $ 2$ & \mazethree{0}{1}{0} & \dial{0}\dial{1}\dial{0} \\[-0.2em]
$ 8$ & $020$ & $1$ & $ 1$ & $ 8$ & $020$ & $1$ & $ 1$ & $ 8$ & $020$ & $1$ & $ 1$ & \mazethree{0}{2}{0} & \dial{0}\dial{2}\dial{0} \\[-0.2em]
$ 9$ & $021$ & $1$ & $ 1$ & $ 9$ & $021$ & $1$ & $ 1$ & $ 9$ & $021$ & $1$ & $ 1$ & \mazethree{0}{2}{1} & \dial{0}\dial{2}\dial{1} \\[-0.2em]
$10$ & $022$ & $1$ & $ 1$ & $10$ & $022$ & $1$ & $ 1$ & $10$ & $022$ & $1$ & $ 1$ & \mazethree{0}{2}{2} & \dial{0}\dial{2}\dial{2} \\[-0.2em]
$11$ & $023$ & $2$ & $ 2$ & $11$ & $023$ & $2$ & $ 2$ & $11$ & $023$ & $2$ & $ 2$ & \mazethree{0}{2}{3} & \dial{0}\dial{2}\dial{3} \\[-0.2em]
$12$ & $033$ & $1$ & $-1$ & $12$ & $033$ & $3$ & $ 3$ & $12$ & $033$ & $3$ & $ 3$ & \mazethree{0}{3}{3} & \dial{0}\dial{3}\dial{3} \\[-0.2em]
$13$ & $032$ & $1$ & $-1$ &&&&&&&&&&\\[-0.2em]
$14$ & $031$ & $1$ & $-1$ &&&&&&&&&&\\[-0.2em]
$15$ & $030$ & $3$ & $ 3$ &&&&&&&&&&\\[-0.2em]
$16$ & $130$ & $1$ & $ 1$ &&&&&&&&&&\\[-0.2em]
$17$ & $131$ & $1$ & $ 1$ &&&&&&&&&&\\[-0.2em]
$18$ & $132$ & $1$ & $ 1$ &&&&&&&&&&\\[-0.2em]
$19$ & $133$ & $2$ & $-2$ & $13$ & $133$ & $2$ & $-2$ & $13$ & $133$ & $2$ & $-2$ & \mazethree{1}{3}{3} & \dial{1}\dial{3}\dial{3} \\[-0.2em]
$20$ & $123$ & $1$ & $-1$ & $14$ & $123$ & $1$ & $-1$ & $14$ & $123$ & $1$ & $-1$ & \mazethree{1}{2}{3} & \dial{1}\dial{2}\dial{3} \\[-0.2em]
$21$ & $122$ & $1$ & $-1$ & $15$ & $122$ & $1$ & $-1$ & $15$ & $122$ & $1$ & $-1$ & \mazethree{1}{2}{2} & \dial{1}\dial{2}\dial{2} \\[-0.2em]
$22$ & $121$ & $1$ & $-1$ & $16$ & $121$ & $1$ & $-1$ & $16$ & $121$ & $1$ & $-1$ & \mazethree{1}{2}{1} & \dial{1}\dial{2}\dial{1} \\[-0.2em]
$23$ & $120$ & $2$ & $-2$ & $17$ & $120$ & $2$ & $-2$ & $17$ & $120$ & $2$ & $-2$ & \mazethree{1}{2}{0} & \dial{1}\dial{2}\dial{0} \\[-0.2em]
$24$ & $110$ & $1$ & $ 1$ & $18$ & $110$ & $1$ & $ 1$ & $18$ & $110$ & $1$ & $ 1$ & \mazethree{1}{1}{0} & \dial{1}\dial{1}\dial{0} \\[-0.2em]
$25$ & $111$ & $1$ & $ 1$ & $19$ & $111$ & $1$ & $ 1$ & $19$ & $111$ & $1$ & $ 1$ & \mazethree{1}{1}{1} & \dial{1}\dial{1}\dial{1} \\[-0.2em]
$26$ & $112$ & $1$ & $ 1$ & $20$ & $112$ & $1$ & $ 1$ & $20$ & $112$ & $1$ & $ 1$ & \mazethree{1}{1}{2} & \dial{1}\dial{1}\dial{2} \\[-0.2em]
$27$ & $113$ & $2$ & $-2$ & $21$ & $113$ & $2$ & $-2$ & $21$ & $113$ & $2$ & $-2$ & \mazethree{1}{1}{3} & \dial{1}\dial{1}\dial{3} \\[-0.2em]
$28$ & $103$ & $1$ & $-1$ & $22$ & $103$ & $1$ & $-1$ & $22$ & $103$ & $3$ & $ 3$ & \mazethree{1}{0}{3} & \dial{1}\dial{0}\dial{3} \\[-0.2em]
$29$ & $102$ & $1$ & $-1$ & $23$ & $102$ & $1$ & $-1$ & $   $ & $  $ & $ $ & $  $ & \mazethree{1}{0}{2} & \dial{1}\dial{0}\dial{2} \\[-0.2em]
$30$ & $101$ & $1$ & $-1$ & $24$ & $101$ & $1$ & $-1$ & $   $ & $  $ & $ $ & $  $ & \mazethree{1}{0}{1} & \dial{1}\dial{0}\dial{1} \\[-0.2em]
$31$ & $100$ & $3$ & $ 3$ & $25$ & $100$ & $3$ & $ 3$ & $   $ & $  $ & $ $ & $  $ & \mazethree{1}{0}{0} & \dial{1}\dial{0}\dial{0} \\[-0.2em]
$32$ & $200$ & $1$ & $ 1$ & $26$ & $200$ & $1$ & $ 1$ & $   $ & $  $ & $ $ & $  $ & \mazethree{2}{0}{0} & \dial{2}\dial{0}\dial{0} \\[-0.2em]
$33$ & $201$ & $1$ & $ 1$ & $27$ & $201$ & $1$ & $ 1$ & $   $ & $  $ & $ $ & $  $ & \mazethree{2}{0}{1} & \dial{2}\dial{0}\dial{1} \\[-0.2em]
$34$ & $202$ & $1$ & $ 1$ & $28$ & $202$ & $1$ & $ 1$ & $   $ & $  $ & $ $ & $  $ & \mazethree{2}{0}{2} & \dial{2}\dial{0}\dial{2} \\[-0.2em]
$35$ & $203$ & $2$ & $ 2$ & $29$ & $203$ & $2$ & $ 2$ & $23$ & $203$ & $2$ & $ 2$ & \mazethree{2}{0}{3} & \dial{2}\dial{0}\dial{3} \\[-0.2em]
$36$ & $213$ & $1$ & $-1$ & $30$ & $213$ & $1$ & $-1$ & $24$ & $213$ & $1$ & $-1$ & \mazethree{2}{1}{3} & \dial{2}\dial{1}\dial{3} \\[-0.2em]
$37$ & $212$ & $1$ & $-1$ & $31$ & $212$ & $1$ & $-1$ & $25$ & $212$ & $1$ & $-1$ & \mazethree{2}{1}{2} & \dial{2}\dial{1}\dial{2} \\[-0.2em]
$38$ & $211$ & $1$ & $-1$ & $32$ & $211$ & $1$ & $-1$ & $26$ & $211$ & $1$ & $-1$ & \mazethree{2}{1}{1} & \dial{2}\dial{1}\dial{1} \\[-0.2em]
$39$ & $210$ & $2$ & $ 2$ & $33$ & $210$ & $2$ & $ 2$ & $27$ & $210$ & $2$ & $ 2$ & \mazethree{2}{1}{0} & \dial{2}\dial{1}\dial{0} \\[-0.2em]
$40$ & $220$ & $1$ & $ 1$ & $34$ & $220$ & $1$ & $ 1$ & $28$ & $220$ & $1$ & $ 1$ & \mazethree{2}{2}{0} & \dial{2}\dial{2}\dial{0} \\[-0.2em]
$41$ & $221$ & $1$ & $ 1$ & $35$ & $221$ & $1$ & $ 1$ & $29$ & $221$ & $1$ & $ 1$ & \mazethree{2}{2}{1} & \dial{2}\dial{2}\dial{1} \\[-0.2em]
$42$ & $222$ & $1$ & $ 1$ & $36$ & $222$ & $1$ & $ 1$ & $30$ & $222$ & $1$ & $ 1$ & \mazethree{2}{2}{2} & \dial{2}\dial{2}\dial{2} \\[-0.2em]
$43$ & $223$ & $2$ & $ 2$ & $37$ & $223$ & $2$ & $ 2$ & $31$ & $223$ & $2$ & $ 2$ & \mazethree{2}{2}{3} & \dial{2}\dial{2}\dial{3} \\[-0.2em]
$44$ & $233$ & $1$ & $-1$ & $38$ & $233$ & $3$ & $ 3$ & $32$ & $233$ & $3$ & $ 3$ & \mazethree{2}{3}{3} & \dial{2}\dial{3}\dial{3} \\[-0.2em]
$45$ & $232$ & $1$ & $-1$ &&&&&&&&&& \\[-0.2em]
$46$ & $231$ & $1$ & $-1$ &&&&&&&&&& \\[-0.2em]
$47$ & $230$ & $3$ & $ 3$ &&&&&&&&&& \\[-0.2em]
$48$ & $330$ & $1$ & $ 1$ &&&&&&&&&& \\[-0.2em]
$49$ & $331$ & $1$ & $ 1$ &&&&&&&&&& \\[-0.2em]
$50$ & $332$ & $1$ & $ 1$ &&&&&&&&&& \\[-0.2em]
$51$ & $333$ & $2$ & $-2$ & $39$ & $333$ &&& $33$ & $333$ &&& \mazethree{3}{3}{3} & \dial{3}\dial{3}\dial{3} \\[-0.2em]
$52$ & $323$ & $1$ & $-1$ &&&&&&&&&&\\[-0.2em]
$53$ & $322$ & $1$ & $-1$ &&&&&&&&&&\\[-0.2em]
$54$ & $321$ & $1$ & $-1$ &&&&&&&&&&\\[-0.2em]
$55$ & $320$ & $2$ & $-2$ &&&&&&&&&&\\[-0.2em]
$56$ & $310$ & $1$ & $ 1$ &&&&&&&&&&\\[-0.2em]
$57$ & $311$ & $1$ & $ 1$ &&&&&&&&&&\\[-0.2em]
$58$ & $312$ & $1$ & $ 1$ &&&&&&&&&&\\[-0.2em]
$59$ & $313$ & $2$ & $-2$ &&&&&&&&&&\\[-0.2em]
$60$ & $303$ & $1$ & $-1$ &&&&&&&&&&\\[-0.2em]
$61$ & $302$ & $1$ & $-1$ &&&&&&&&&&\\[-0.2em]
$62$ & $301$ & $1$ & $-1$ &&&&&&&&&&\\[-0.2em]
$63$ & $300$ & &&&&&&&&&&&\\[-0.2em]
\end{tabular}}
\caption{\small For $n=2$ digits, we examine the quaternary reflected Gray code $\QUAT{n}$, the longest solution $\LONG{n}$, and the shortest solution $\SHORT{n}$. For each order, we include the ranking, the Gray code, the change sequence, and the unsigned change sequence. The last two columns present the maze representation corresponding to each $2$-digit base $4$ string and the dials representation for the Zigguhooked puzzle.}
\label{tab:orders}

\end{table}

\subsection{Outline}
\label{sec:intro_outline}

Section \ref{sec:solving} gives a non-technical description of how to solve Gray code and Ziggu puzzles.
Section \ref{sec:mazes} explains our Ziggu maze abstraction and how the mazes are formed in Zigguflat and Zigguhooked.
Section \ref{sec:lengths} derives recursive formulae for shortest and longest solution lengths. 
Section \ref{sec:reflected} cover the binary and quaternary reflected Gray codes. 
Section \ref{sec:recursive} provides recursive descriptions of Gray code and Ziggu puzzle solutions.
Sections \ref{sec:sequences}, \ref{sec:successor}, and \ref{sec:position} present the change sequences, successor rules, and comparison and ranking formulae for the solutions to these puzzles. 
Section \ref{sec:loopless} presents loopless algorithms for generating the shortest solution to Ziggu puzzles. 
We present an exciting connection between the shortest solution to Ziggu puzzles and the $2\times n$ grids arising from the pencil-and-paper puzzle Nurikabe in Section \ref{sec:nurikabe}. 
We end with final remarks in Section \ref{sec:final}.
Omitted proofs follow from first principles.

\section{Solving Puzzles from the Start}
\label{sec:solving}

In this section we give simple non-technical descriptions of how to solve Gray code puzzles and Ziggu puzzles from their start states.
The content is tailored to those who wish to solve a puzzle one time as simply as possible.
It also assumes that the solver begins from the start state and knows what the last move was.
In other words, they either solve the puzzle in one shot, or if they take a break, then they take sufficient notes to resume the solution.
Subsequent sections in this paper formalize these results and address more nuanced questions and scenarios.

\subsection{Don't Look Back!}
\label{sec:solving_back}

\begin{figure}[ht]
    \centering
    \begin{subfigure}{0.495\textwidth}
        \centering
        \includegraphics[scale=0.11]{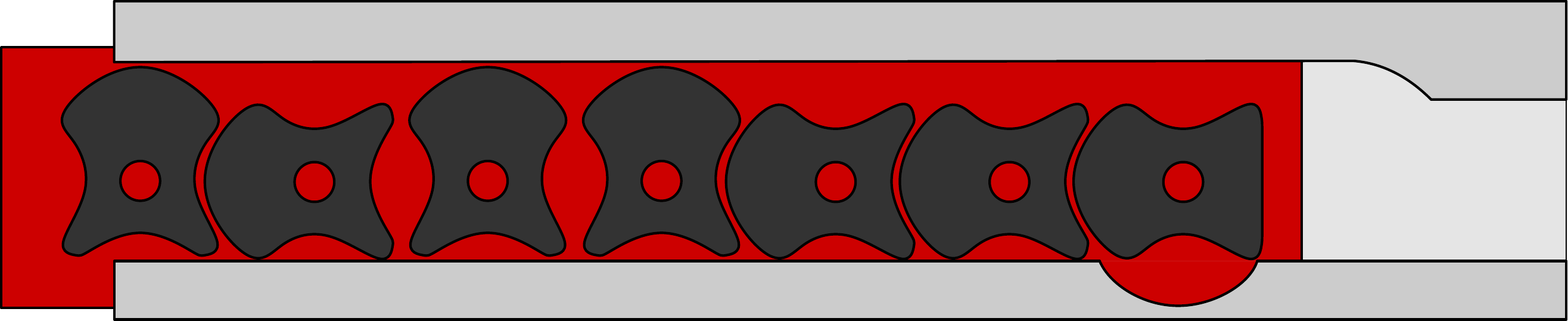} 
        \caption{Rotate dial $1$ for $101100\overline{0} = 1011001$.}
        \label{fig:statesSpinOut_dial7}
    \end{subfigure}
    \hfill
    \begin{subfigure}{0.495\textwidth}
        \centering
        \includegraphics[scale=0.11]{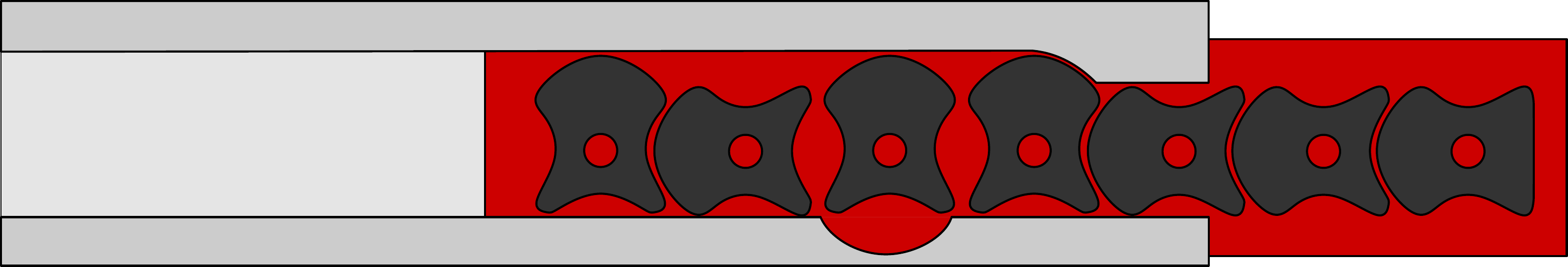} 
        \caption{Rotate dial $3$ for $10\overline{1}1000 = 1001000$.}
        \label{fig:statesSpinOut_dial3}
    \end{subfigure}
    \caption{
    \emph{Spin-Out} states are encoded as a binary string $b_n b_{n-1} \cdots b_1$ based on the orientation of its $n$ dials.
    Horizontal translation is ignored, so (a) and (b) are both encoded as $b_7 b_6 \cdots b_1 = 1011000$.
    The contours of the top and bottom rails ensure that at most two dials can be rotated in any given state.
    These include the rightmost dial $b_1$ as in (a), and the dial to the left of the rightmost vertical dial as in (b).
    In particular, $b_1$ or $b_2$ can be changed from its start state $b_7 b_6 \cdots b_1 = 1111111$, with the former change leading in the wrong direction (cf., Figure \ref{fig:statesBrain4}).
    }
    \label{fig:statesSpinOut}
\end{figure}

Gray code puzzles and Ziggu puzzles can be solved using the same basic principle:
don't look back!
This means that you should never reverse your progress by undoing your most recent move.
This advice is enough to solve many Gray code puzzles due to the following point.
\begin{itemize}
    \item Gray code puzzles have at most one non-reversing move.  
\end{itemize}
In other words, when you move to a new state there is at most one move that doesn't take you back to your previous state.
For many Gray code puzzles (e.g., \emph{Spin-Out}) the only complication is that there are two different moves that can be made from the start state;
the wrong first move leads toward a dead-end from which you need to undo your final move and \emph{then} never look back.
Another way of saying this is that the state space is a path and only one direction takes you to the solution.
For example, the neighbors in Spin-Out%
\footnote{A design flaw in the initial run of \emph{Spin-Out} puzzles introduces more states and shorter solutions by allowing dials on the right to be rotated to the right.
This can be fixed by extending the right end of the frame.
Our discussion of \emph{Spin-Out} ignores the flaw.
} 
are illustrated in Figure \ref{fig:statesSpinOut}.
In \emph{Towers of Hanoi} there can be two non-reversing moves, but this reduces to one when putting additional restrictions on how and when to move the smallest disc.%
\footnote{On every odd numbered turn, move disc $1$ one peg to the left (right) for even (odd) $n$.}

\begin{figure}
    \centering
    \begin{subfigure}{1.0\textwidth}
        \centering
        \includegraphics[width=0.8\linewidth,page=25]{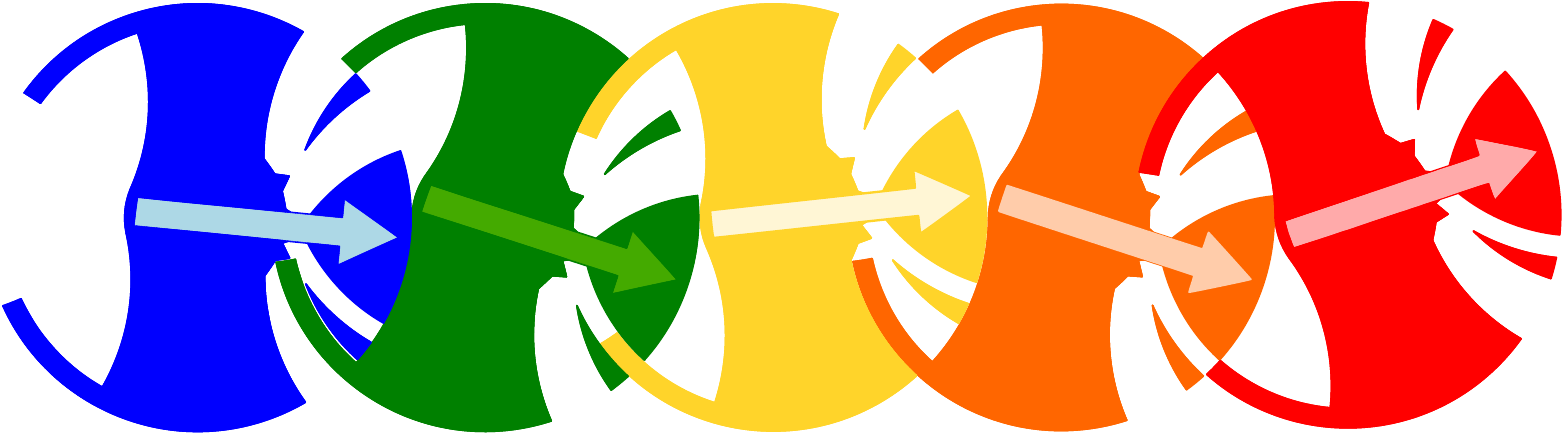}
        \caption{In $\digite{1}\digitd{0}\digitc{2}\digitb{0}\digita{3}$ the digit $\digitd{0}$ is locked by the \digitc{2} digit to its right, but the other four digits can change.
        In the corresponding Zigguhooked state the \digitd{green} dial cannot rotate in either direction due to the orientation of the \digitc{yellow} dial, but all of the other dials can rotate in one direction.
        The number of choices makes this a more difficult state to begin from.
        }
        \label{fig:numMoves_10203}
    \end{subfigure}
    \begin{subfigure}{0.24\textwidth}
        \vspace{0.5em}
        \centering
        \includegraphics[width=1.0\linewidth,page=26]{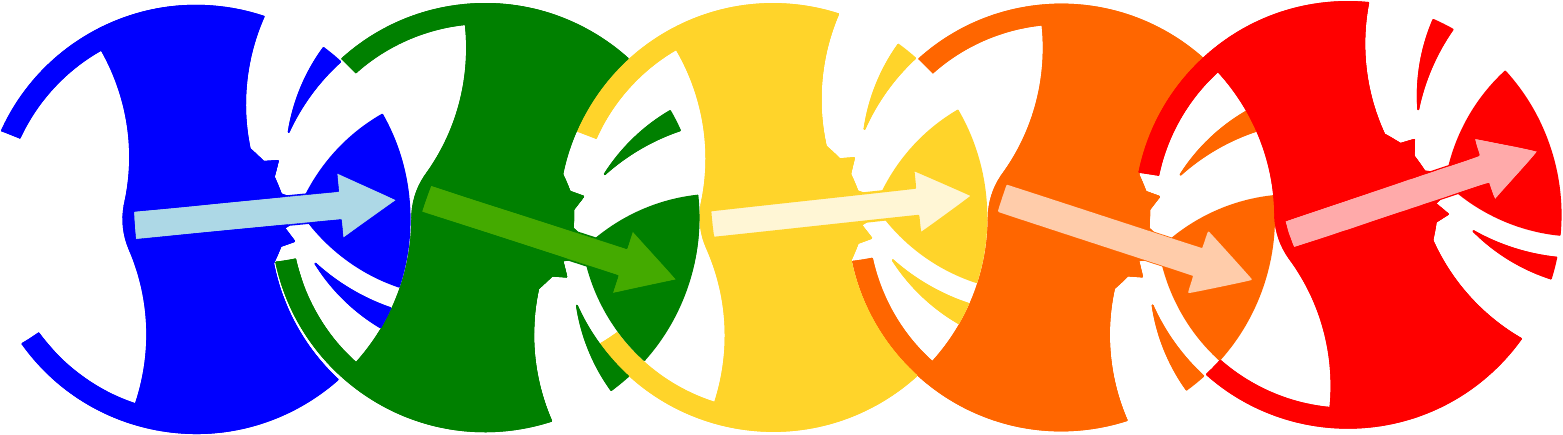}
        \caption{$\overline{\digite{2}}\digitd{0}\digitc{2}\digitb{0}\digita{3}$}
        \label{fig:numMoves_20203}
    \end{subfigure}
    \begin{subfigure}{0.24\textwidth}
        \centering
        \includegraphics[width=1.0\linewidth,page=27]{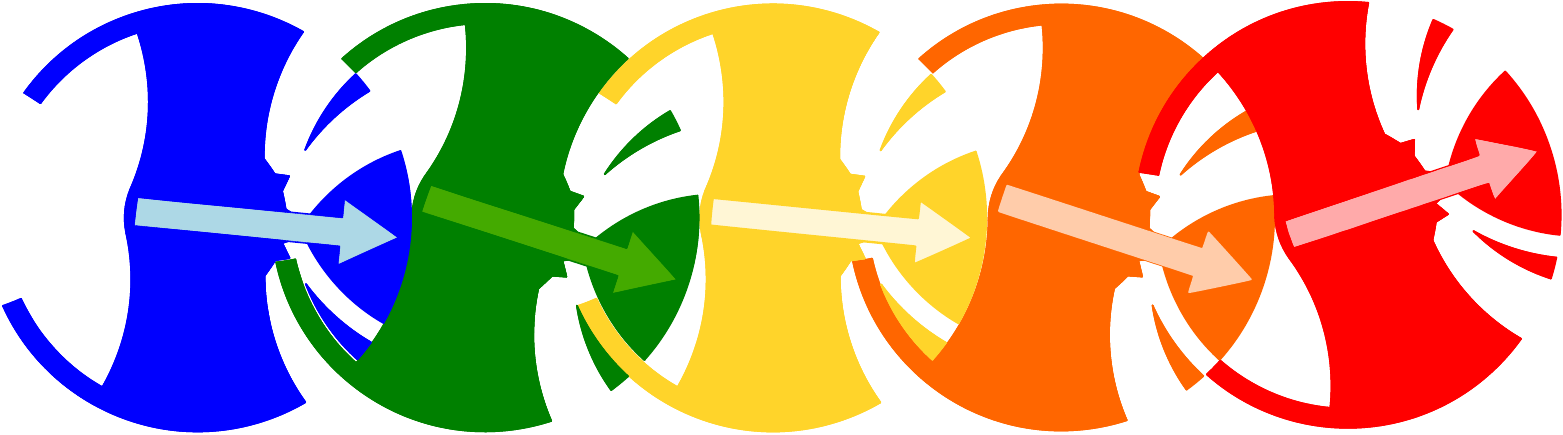}
        \caption{$\digite{1}\digitd{0}\underline{\digitc{1}}\digitb{0}\digita{3}$}
        \label{fig:numMoves_10103}
    \end{subfigure}
    \begin{subfigure}{0.24\textwidth}
        \centering
        \includegraphics[width=1.0\linewidth,page=28]{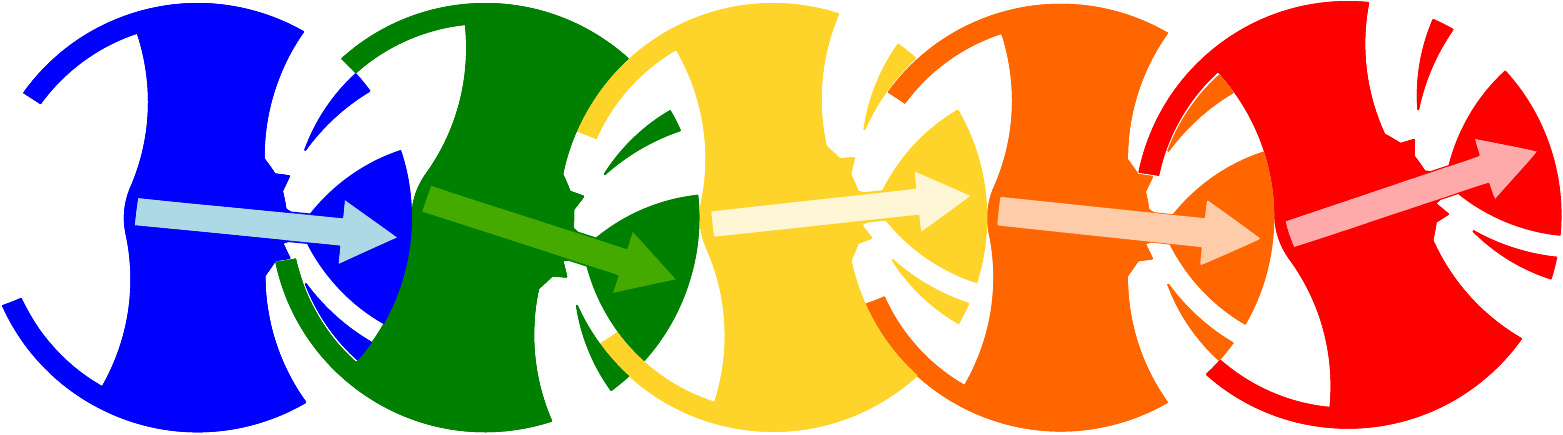}
        \caption{$\digite{1}\digitd{0}\digitc{2}\overline{\digitb{1}}\digita{3}$}
        \label{fig:numMoves_10213}
    \end{subfigure}
    \begin{subfigure}{0.24\textwidth}
        \centering
        \includegraphics[width=1.0\linewidth,page=29]{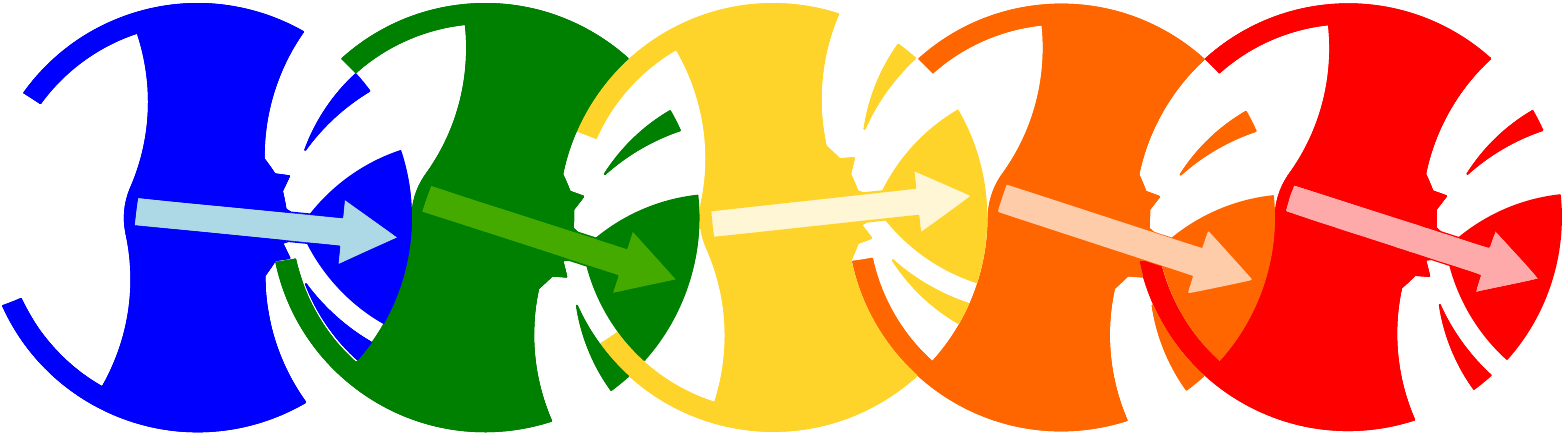}
        \caption{$\digite{1}\digitd{0}\digitc{2}\digitb{0}\underline{\digita{0}}$ UPDATE}
        \label{fig:numMoves_10200}
    \end{subfigure}
    \caption{
    In Ziggu puzzles there are states in which up to $\lceil \frac{n}{2} \rceil + 1$ digits can be changed.
    This is illustrated in (a) by a Zigguhooked state in which $\lceil \frac{5}{2} \rceil + 1 = 4$ of the $n=5$ dials can be changed.
    Digits that can change can only be changed in one way (i.e., incremented or decremented).
    These changes are shown in (b)--(e) with overlines and underlines indicating increments (counterclockwise rotation) and decrements (clockwise rotation), respectively.
    Note that (a) would not be encountered on the shortest solution since states on the shortest solution have at most three digits that can be changed (regardless of $n$). 
    (An exception is the rightmost \digita{red} digit, which can be freely incremented or decremented when it is not $\digita{0}$ or $\digita{3}$.)
    }
    \label{fig:numMoves}
\end{figure}

Ziggu puzzles are more challenging than Gray code puzzles due to the following.
\begin{itemize}
    \item Ziggu puzzles can have more than one non-reversing move.
    \item Ziggu puzzles have solutions of different lengths.
\end{itemize}
A puzzle with $m$ mazes and $n=m+1$ digits can have up to $\lceil\frac{n}{2}\rceil+1$ non-reversing moves%
\footnote{For example, the state $2020 \cdots 202$ (see Section \ref{sec:mazes}) allows changes to each $2$ and the last $0$.}\textsuperscript{,}%
\footnote{The shortest solution only visits states with up to two non-reversing moves.
In this sense, Ziggu puzzles remain relatively simple until you go off course.}.
Figure \ref{fig:numMoves} gives an example of such a state in Zigguhooked.
While Ziggu puzzles have more complicated states, the unique shortest and longest solutions are generated by simple extremal versions of the basic don't look back principle.
\begin{center}
The \textbf{shortest} solution always makes the \textbf{leftmost} non-reversing move. \\
The \textbf{longest} solution always makes the \textbf{rightmost} non-reversing move.
\end{center}
These processes are simplified by the fact that there is only one available move from the start state, so you always begin by moving in the correct direction.
Furthermore, besides the start state the only other ``dead end'' is the solution state.
As a result, you'll never need to backtrack over previous moves.
However, you can return to a previous state by going around in a circle (i.e., through a cycle of states) if you don't follow either of the extremal versions of the basic principle.

\subsection{One Bit to Remember: Resuming a Shortest Solution}
\label{sec:solving_resuming}

By design exponential puzzles take a long time to solve.
So solvers may wish to take a break.
The danger is working backwards when resuming a solution.
In other words, a solver may start by undoing the last move that they made before taking a break.
That can be avoided by making note of the last move that was made.
In fact, less information is needed.
More specifically, only one bit of information is needed.
\underline{One literal bit}.

Gray code puzzles have a ``metronome'' property wherein every second change involves the first piece.
More specifically, the change sequence follows the binary ruler sequence: $1, 2, 1, 3, 1, 2, 1, 4, 1, \ldots$ (\OEIS{A001511}).
Therefore, you only need to remember if your most recent move involved the first piece or not.
The metronome property holds for puzzles that start in an intermediate state (e.g., Spin-Out) although the change sequence is a (reversed) subsequence of \OEIS{A001511}.

The shortest solution to Ziggu puzzles instead have a ``stair-climbing'' property wherein successive changes differ by at most one: $1,1,1,2,1,1,1,2,1,1,1,2,3,2,1,\ldots$. 
Suppose that you always make the $1$ moves in succession (i.e., when changing the first piece you do it three times in a row).
Treating each $1,1,1$ as a single $1$ is natural in most of the Ziggu puzzles (e.g., the first dial in Zigguhooked moves freely between its two extreme states $0$ and $3$).
With this assumption you just remember if your most recent move went up (e.g., $2$ then $3$) or down (e.g., $3$ then $2$).
In the first case you will attempt to move up again (e.g., $2$ then $3$ then $4$);
if the puzzle does not allow this, then you will move down (e.g., $2$ then $3$ then $2$).
In the second case you move down again (e.g., $3$ then $2$ then $1$) or repeatedly move the first piece.

Here is a general tip to avoid undoing your work in a Gray code or Ziggu puzzle.
\begin{center}
Never pause after changing the first piece or going down. 
\end{center}

\section{Ziggu Mazes, Pieces, and Digits}
\label{sec:mazes}

While the Ziggu family consists of a wide family of seemingly distinct puzzles, they all use the same underlying mechanic of the Ziggu maze. 

\begin{figure}
    \centering
    \begin{subfigure}{0.3\linewidth}
        \centering
        \includegraphics[scale=0.4, page=1]{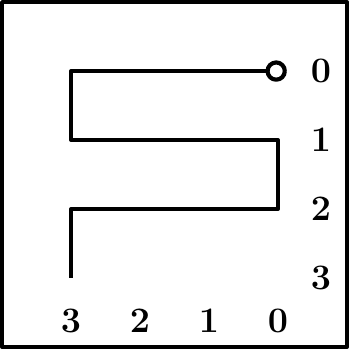}
        \caption{Initial state}
        \label{fig:maze_initial}        
    \end{subfigure}
    \hfill
    \begin{subfigure}{0.3\linewidth}
        \centering
        \includegraphics[scale=0.4, page=2]{mazes.pdf}
        \caption{13 valid locations}
        \label{fig:maze_locations}        
    \end{subfigure}
    \hfill
    \begin{subfigure}{0.3\linewidth}
        \centering
        \includegraphics[scale=0.4, page=3]{mazes.pdf}
        \caption{Maze exit}
        \label{fig:maze_exit}        
    \end{subfigure}
    \caption{(a) The standard Ziggu maze is on a $4$-by-$4$ grid and its position is initially location $(0,0)$.
    (b) There are $13$ valid locations: $(3,i)$ is invalid for $0 \leq i \leq 2$.
    (c) The maze's exit is accessible from location $(3,3)$ (but only if higher mazes have been exited). 
    }
    \label{fig:maze}
\end{figure}

\subsection{Ziggu Mazes}
\label{sec:mazes_mazes}
\begin{definition}
\label{def:maze}
    A \emph{Ziggu maze} is the $\mathsf{S}$-shaped maze in Figure \ref{fig:maze}. It lives on a $4\times 4$ grid, and the 13 valid positions are
    \begin{equation}
        \set{(i,j)\mid i\in\set{0,1,2},j\in\set{0,1,2,3}}\cup\set{(3,3)}.
    \end{equation}
    The rows and columns are labeled $0,1,2,3$ from top to bottom and from right to left, respectively. The 13 valid states are connected as shown in Figure \ref{fig:maze}. The \emph{initial state} of the maze is $(0,0)$, and the maze's \emph{exit} is accessible from $(3,3)$. 
\end{definition}
\begin{definition}
\label{def:m-mazes}
    An \emph{$m$-Ziggu maze} is a sequence of $m$ Ziggu mazes in states 
    \begin{equation}
       (r_1,c_1),(r_2,c_2),\dots,(r_m,c_m)\text{ where } r_i=c_{i+1}\text{ for all } i=1,2,\dots,m-1. 
    \end{equation}
That is, the horizontal position of a maze is the vertical position of the maze to its right. 
\end{definition}
For example, the $3$-Ziggu maze with state $(0,1),(2,0),(3,2)$ has the maze representation \mazefour{1}{0}{2}{3}. While there are $6$ digits in the state description, since $r_1=c_2$ and $r_2=c_3$, we can describe the state of this $3$-Ziggu maze with $4$ digits. 
\begin{definition}
    The \emph{base $4$ representation of an $m$-Ziggu maze} is the length $n=m+1$ base $4$ string $c_1r_1r_2\cdots r_m$.
\end{definition}
Considering the previous example of the $3$-Ziggu maze \mazefour{1}{0}{2}{3}\, in state $(0,1)$, $(2,0)$, $(3,2)$, its base $4$ representation is the base $4$ string $1023$.

Definition \ref{def:maze} immediately gives a bijection between the set of $m$-Ziggu mazes and the set of length $m+1$ base $4$ strings in which $3$ is the only digit appearing to the right of a $3$.

\subsection{Implementation of Ziggu mazes in the puzzles}
\label{sec:mazes_implementations}

The first puzzle in the Ziggu family was Ziggurat. Each \emph{Ziggurat} piece contains a Ziggu maze, and successive pieces have alternating chirality.
For example, the blue pieces in Figure \ref{fig:ziggus_ziggurat} have an \rotatebox{45}{$\mathsf{S}$}-shaped maze and the orange pieces have an \reflectbox{\rotatebox{45}{$\mathsf{S}$}}-shaped maze.
The bottom of each piece contains a nub that navigates the maze two pieces below it.
For example, the nub of the top blue piece passes through the top orange piece and stays within the maze of next blue piece.
The puzzle is solved when every maze is solved (i.e., every nub reaches the end of its respective maze) which allows it to be disassembled.
The challenge comes from the fact that consecutive pieces move in lockstep in one of the two directions of the maze.
This is accomplished by the railing at the bottom of each piece slotting into the channel of the piece below it, as shown on the right side of \ref{fig:ziggus_ziggurat}.
As a result, making progress in one maze can involve moving to the opposite end of the previous maze, and this cascade of back-and-forth movements gives the puzzle its exponential nature.

The Ziggurat puzzle inspired Oskar van Deventer to create the vast family of Ziggu puzzles shown in Figure \ref{fig:ziggus}. These puzzles all use the Ziggu maze, but they feel distinct and provide differing enjoyable puzzling experiences.

\subsubsection{Zigguflat}
\label{sec:mazes_implementations_flat}

\begin{figure}
    \centering
    \includegraphics[width=0.5\linewidth, page=35]{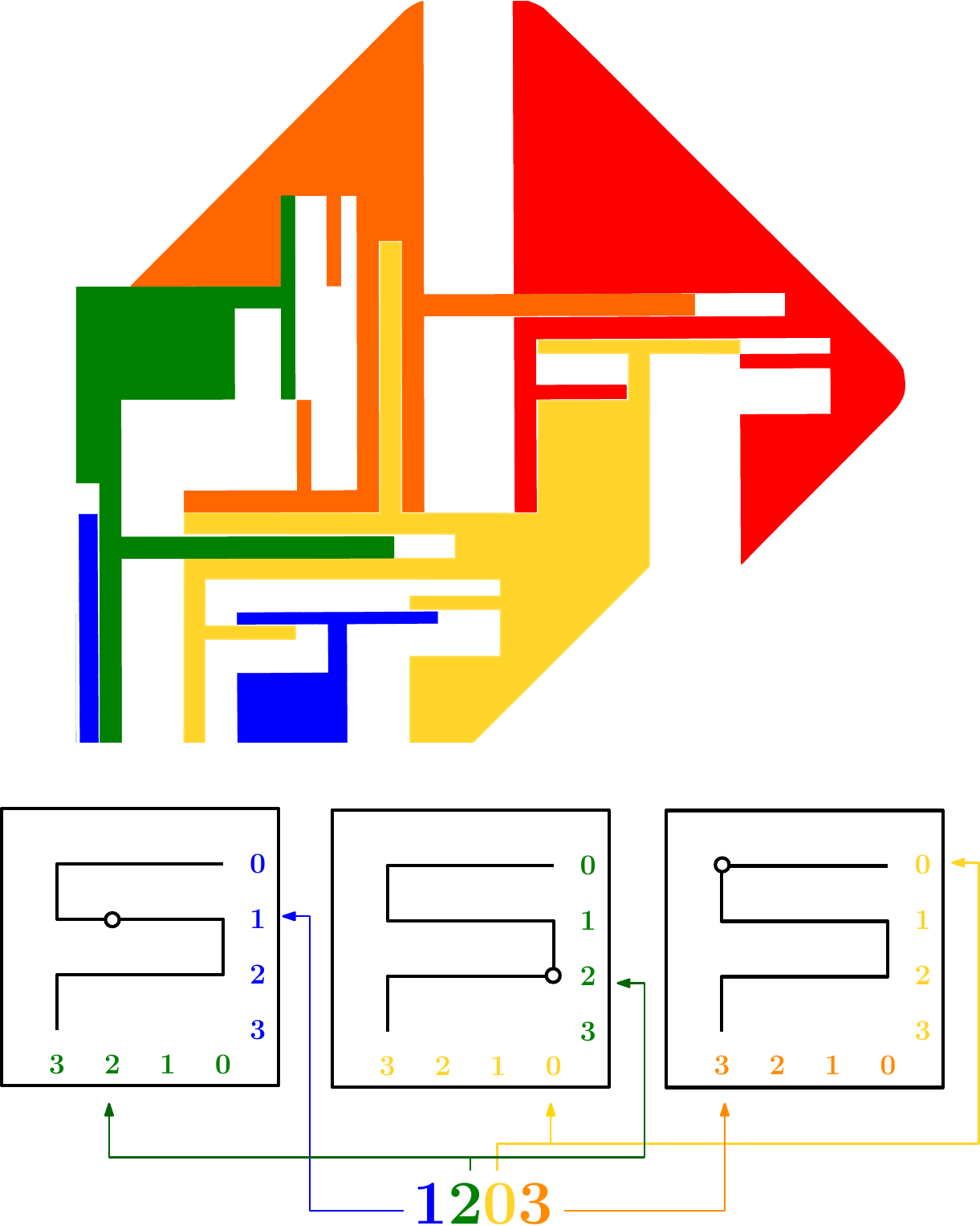}
    \caption{
    Position \digite1\digitd2\digitc0\digitb3 of Zigguflat with three mazes, along with the corresponding maze representation
    \mazefour{1}{2}{0}{3}.  
    The \digite{1} gives the vertical position in the left maze, the \digitd{2} gives the horizontal position in the left maze and the vertical position in the middle maze, the \digitc{0} gives the horizontal position in the middle maze and the vertical position in the right maze, and the \digitb{3} gives the horizontal position in the right maze.}
    \label{fig:zigguflat-maze}
    \includegraphics[width=0.5\linewidth, page=24]{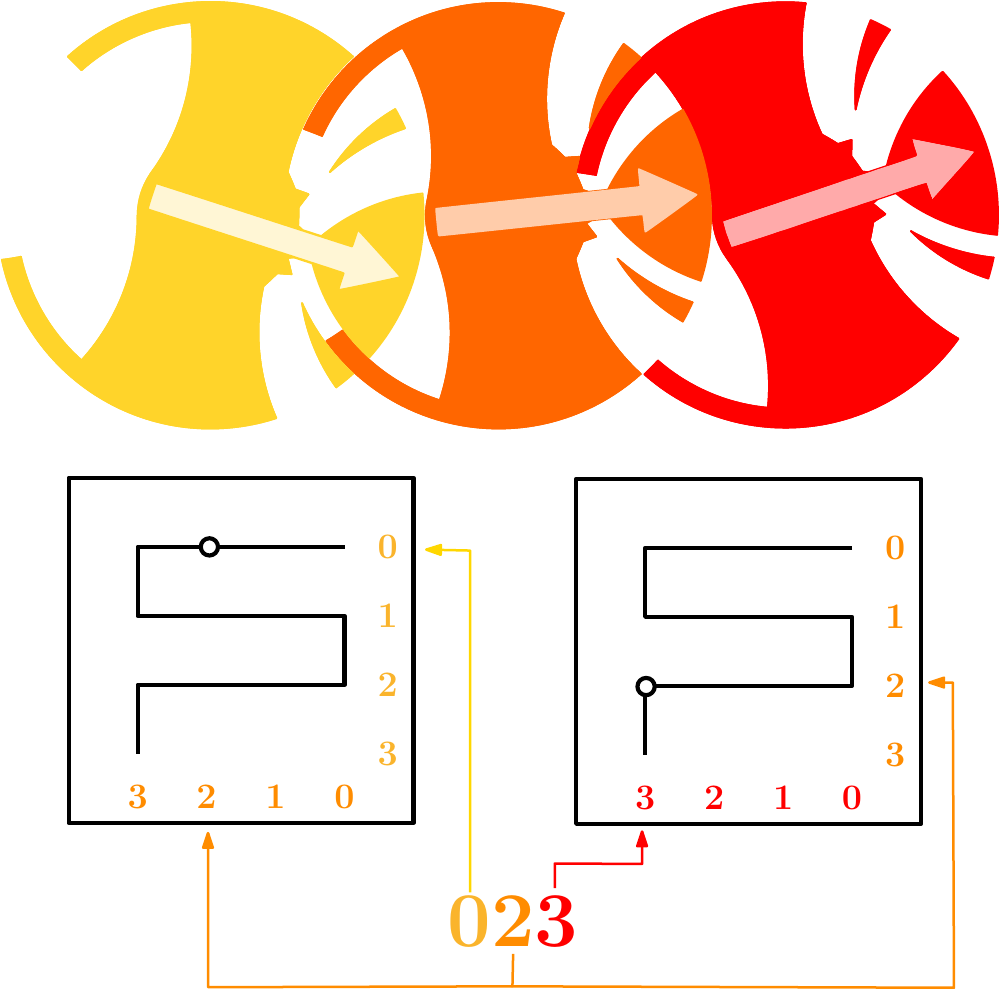}
    \caption{
       Position \digitc0\digitb2\digita3 of Zigguhooked with 3 dials, along with the corresponding maze
    \mazethree{0}{2}{3} and dials \dial{0}\dial{2}\dial{3}. The two mazes exist at the interfaces of adjacent dials.
    The \digitc{0} gives the vertical position in the left maze, the \digitb{2} gives the horizontal position in the left maze and the vertical position in the right maze, and the \digita{3} gives the horizontal position in the right maze.}
    \label{fig:zigguhooked-with-maze}
\end{figure}

The Zigguflat puzzle consists of $p=6$ interlocking pieces as shown in Figure \ref{fig:ziggus_zigguflat}. The mazes are cut into the green, yellow, orange, and red pieces, and those mazes are navigated by $\mathsf T$-shaped nubs of the purple, blue, green, and yellow pieces, respectively. Therefore Zigguflat has $m=p-2=4$ mazes. Notice that the orange and red pieces only do not have $\mathsf T$
-shaped nubs and the purple and blue pieces do not have mazes so that the fully assembled puzzle is a rectangle. As the commercial version of Zigguflat has $m=4$ mazes, it requires $n=m+1=5$ digits to fully describe its state. The correspondence between the physical puzzle state, the digits, and the mazes is shown for a partial version of the Zigguflat puzzle in Figure \ref{fig:zigguflat-maze}. To read off the digits describing the state of the commercial Zigguflat puzzle with 6 pieces, the first four digits describe the vertical positions of the purple, blue, green, and yellow nubs within the mazes they navigate, respectively. The final fifth digit describes the horizontal position of the yellow nub within the red maze.

\subsubsection{Zigguhooked}
\label{sec:mazes_implementations_hooked}

\begin{figure}
    \centering
    \begin{subfigure}{0.495\textwidth}
        \centering
        \includegraphics[width=0.8\textwidth,page=14]{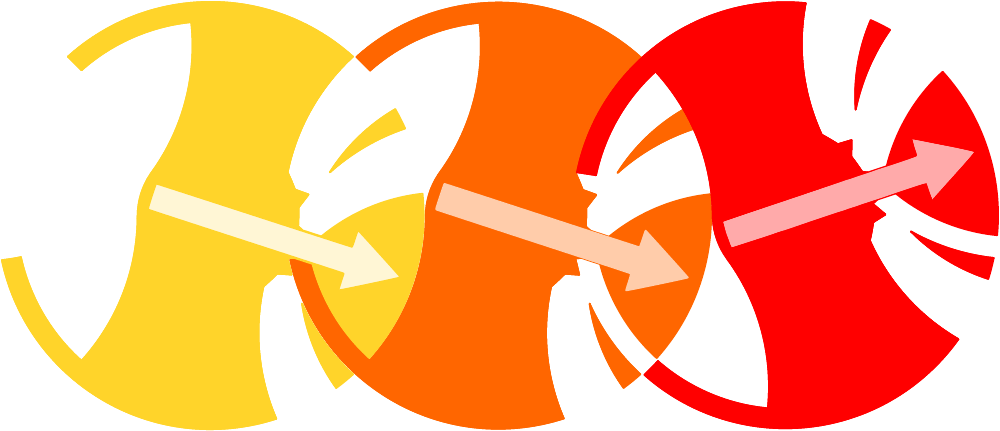} 
        \caption{3-dial Zigguhooked in state \digitc0\digitb0\digita3, \mazethree{0}{0}{3}}
        \label{fig:state-hooked-033}
    \end{subfigure}
    \hfill
    \begin{subfigure}{0.495\textwidth}
        \centering
        \includegraphics[width=0.8\textwidth,page=15]{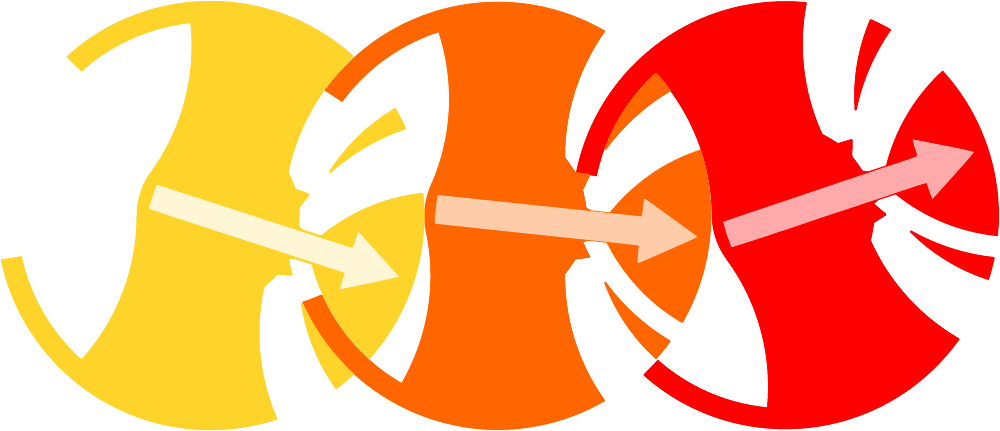} 
        \caption{3-dial Zigguhooked in state  \digitc0\digitb1\digita3, \mazethree{0}{1}{3}}
        \label{fig:state-hooked-133}
    \end{subfigure}
    \caption{Two states of Zigguhooked with three dials. To change between states (a) and (b), rotate the orange dial.}
    \label{fig:state-hooked}
\end{figure}

The Zigguhooked puzzle consists of $p$ interlocking dials with $p=6$ in the commercially available version (see Figure \ref{fig:ziggus_zigguhooked}.) 
Each maze is formed at the interface of two consecutive dials, so there are $m=p-1$ mazes in total.
The rotation of the dials are described by $n=p$ digits. 
Given a pair of adjacent dials, rotating the left dial changes the vertical position within the maze, while rotating the right dial changes the horizontal position within the maze, see Figure \ref{fig:state-hooked}. 
One can also read the digits off of the dials without examining the mazes. 
The radial positions are numbered increasing in a counterclockwise direction, where the states $0$, $1$, $2$, and $3$ correspond to the arrow pointing in the approximately 3:45, 3:15, 2:45, and 2:15 directions of the hour hand on a clock face, respectively. 
So a Zigguhooked puzzle in the state $0213$ would have its dials in position \dial{0}\dial{2}\dial{1}\dial{3}. 
The correspondence between the dial positions, digits, and mazes is shown for a 3-dial Zigguhooked puzzle in Figure \ref{fig:zigguhooked-with-maze}.

\section{Recurrences and Closed Formulae: Solution Lengths}
\label{sec:lengths}

In this section we consider how many moves it takes to solve a puzzle.
We do this using recursive formulae that express the number of moves in terms of the number of moves that solve smaller versions of the puzzle.
In other words, we repeat the classic approach to analyzing the \emph{Towers of Hanoi} (see Figure \ref{fig:recurrenceTowers}) to other Gray code puzzles and Ziggu puzzles.
In each case we consider the number of moves $f(n)$ and the number of states $g(n)$ on a solution separately.
These two quantities are closely related --- the number of states is one more than the number of moves (i.e., $g(n) = f(n)+1$) --- but the recursive formulae arrive at these values in different ways.
We also argue that these recurrences give the unique shortest (or longest) solutions.
The recurrences are summarized in Tables \ref{tab:recurrences_Gray}--\ref{tab:recurrences_Ziggu}.

\subsection{Gray Code Puzzle Solution Lengths}
\label{sec:lengths_Gray}

As was shown in Table \ref{tab:Grays}, most Gray code puzzles run between states that are encoded as $0^n$ and $10^{n-1}$, or between $0^n$ and $1^n$ (or vice versa).

\subsubsection{Between $0^n$ and $10^{n-1}$: Towers of Hanoi, B-Bar, \ldots}
\label{sec:lengths_Gray_Towers}

We previously saw the recurrences for the number of moves $f(n)$ and the number of states $g(n)$ required to solve the Towers of Hanoi in Section \ref{sec:intro}.
The two recurrences are provided in the top of Table \ref{tab:recurrences_Gray}.
Theorem \ref{thm:recurrenceTowersHanoi} substitutes in closed forms of these recurrences.

\begin{theorem} \label{thm:recurrenceTowersHanoi}
    The number of states in the shortest solution to the Towers of Hanoi with $n$ discs is $f(n) = 2^{n}$.
    The number of moves on the shortest solution is $g(n) = 2^{n}-1$.
    These quantities count the number of strings between $0^n$ and $10^{n-1}$ in the binary reflected Gray code.
\end{theorem}

\subsubsection{Between $0^n$ and $1^n$: Baguenaudier, The Brain, Spin-Out, \ldots}
\label{sec:lengths_Gray_Brain}

\begin{figure}
    \centering
    \includegraphics[width=0.8\linewidth]{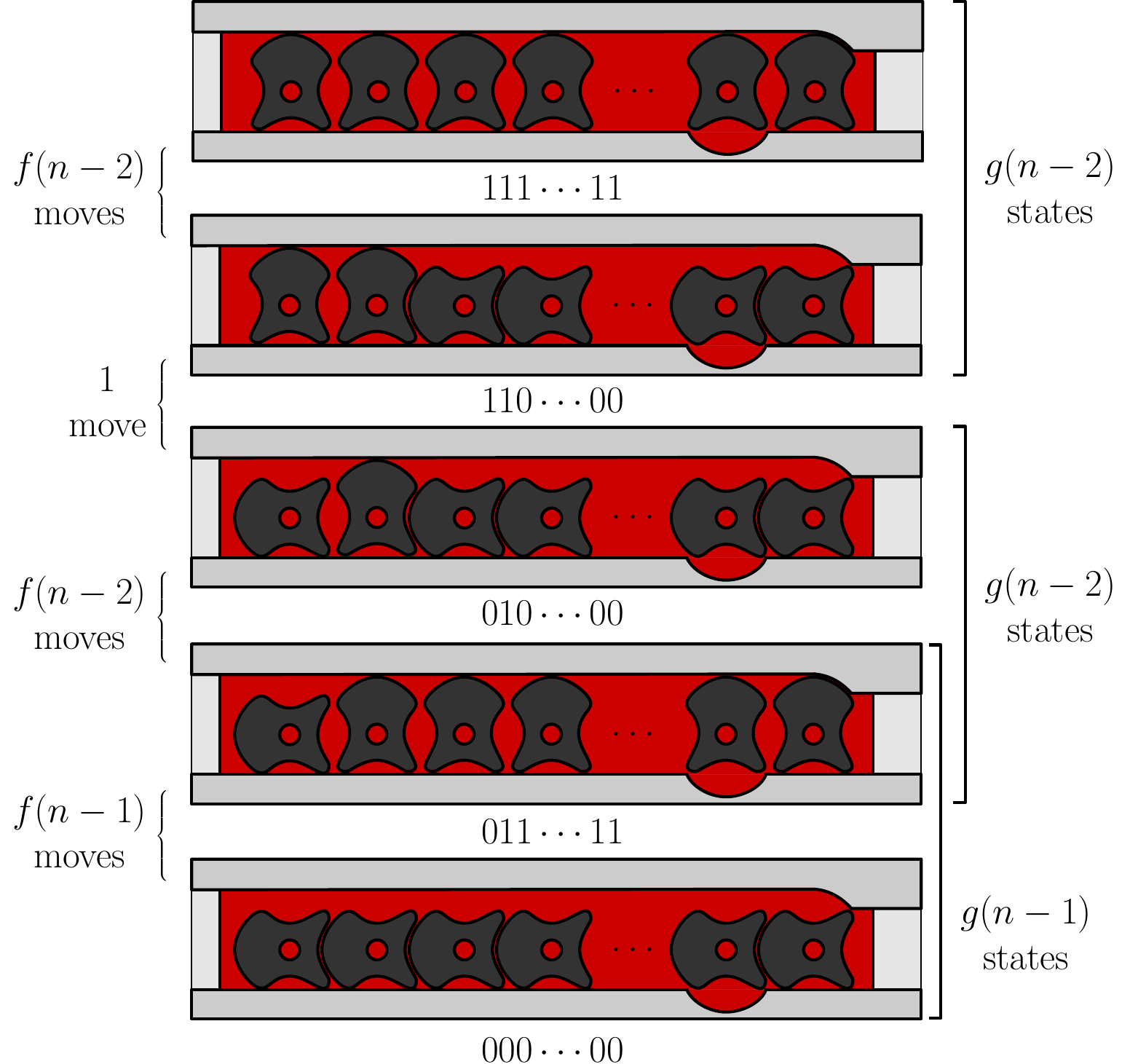}
    \caption{Recurrences for Gray Code puzzles between $0^n$ and $1^n$ (or vice versa).
    \emph{Spin-Out} puzzles with $n$ discs require $f(n) = f(n-1) + 2 f(n-2) + 1$ moves to solve with base cases $f(1) = 1$ and $f(2) = 2$.
    Alternatively, the shortest solution uses $g(n) = g(n-1) + 2g(n-2) - 1$ states (where the $-1$ is from double-counting the state $01^{n-1}$) with $g(1) = 2$ and $g(2) = 3$.
    These recurrences give the number of changes and number of strings in the binary reflected Gray code between $1^n$ and $0^n$.
    They also count the moves and states for additional puzzles including \emph{The Brain}.
    }
    \label{fig:recurrenceSpinOut}
\end{figure}

Recurrences for the number of moves $f(n)$ and the number of states $g(n)$ required to solve Spin-Out appear in Figure \ref{fig:recurrenceSpinOut}.
The two recurrences are provided in the bottom of Table \ref{tab:recurrences_Gray}.
Theorem \ref{thm:recurrenceSpinOut} substitutes in closed forms of these recurrences.

\begin{theorem} \label{thm:recurrenceSpinOut}
    The number of states in the shortest solution to Spin-Out with $n$ dials is $f(n) = \lceil \frac{2^{n+1}}{3} \rceil$.
    The number of moves on the shortest solution is $g(n) = \lfloor \frac{2^{n+1}}{3} \rfloor$
    These quantities count the number of strings between $0^n$ and $1^n$ in the binary reflected Gray code.
\end{theorem}

\begin{table}
\small
\centerline{ 
\begin{tabular}{|c|c|c|c|c|c|c|} 
\hline
\rowcolor{lightgray}
puzzle & quantity & recurrence & base cases & closed form &  sequence \\ \hline
Towers of Hanoi & moves & $2f(n-1) + 1$ & $1$ & $2^n-1$ & \OEIS{A000225} \\ \hline
Towers of Hanoi & states & $2g(n-1)$ & $2$ & $2^n$ & \OEIS{A000079} \\ \hline
Spin-Out & moves & $f(n-1) + 2f(n-2) + 1$ & $1,2$ & $\lfloor 2^{n+1}/3 \rfloor$ & \OEIS{A000975} \\ \hline
Spin-Out & states & $g(n-1) + 2g(n-2) - 1$ & $2,3$ & $\lceil 2^{n+1}/3 \rceil$ & \OEIS{A005578} \\ \hline
\end{tabular}
}
\caption{Formulae for the unique solution for two types of Gray code puzzles with $n$ pieces.  
The base cases are for $n=1,2$.}
\label{tab:recurrences_Gray}
\vspace{1.5em}
\centerline{
\begin{tabular}{|c|c|c|c|c|c|} 
\hline
\rowcolor{lightgray}
quantity & recurrence & base cases & closed form & sequence \\ \hline
most moves & $3f(n-1)+3$ & $3,12$ & $(3^{n+1} - 3)/2$ & \makecell{\OEIS{A029858} \\ \OEIS{A123109}} \\ \hline
most states & $3g(n-1)+1$ & $4,13$ & $(3^{n+1} - 1)/2$ & \OEIS{A003462} \\ \hline
fewest moves & $3f(n-1) - 2f(n-2) + 3$ & $3,12$ & $6 \cdot 2^n - 3 \cdot n - 6$ & N/A \\ \hline
fewest states & $3g(n-1) - 2g(n-2) + 3$ & $4,13$ & $6 \cdot 2^n - 3 \cdot n - 5$ & \OEIS{A101946} \\ \hline
\end{tabular}
}
\caption{Formulae for the unique shortest and longest solutions to Ziggu puzzles with $n$ digits.
The base cases are for $n=1,2$.
}
\label{tab:recurrences_Ziggu}
\end{table}


\subsection{Ziggu Puzzle Solution Lengths}
\label{sec:lengths_ziggu}

\begin{figure}
\centering
\begin{minipage}{0.39\linewidth}
    \begin{subfigure}{1.0\linewidth}
        \centering
        \vspace{2em}
        \includegraphics[scale=1.0,page=7]{mazeMoves.pdf}
        \caption{The top path between $(0,0)$ and $(3,0)$ is highlighted.}
        \label{fig:mazeMoves_top}    
    \end{subfigure}
    \begin{subfigure}{1.0\linewidth}
        \centering
        \includegraphics[scale=1.0,page=5]{mazeMoves.pdf}
        \caption{The shortest solution uses the top path once from $(0,0)$ to $(3,0)$.}
        \label{fig:mazeMoves_short}    
    \end{subfigure}
    \begin{subfigure}{1.0\linewidth}
        \centering
        \includegraphics[scale=1.0,page=6]{mazeMoves.pdf}
        \caption{The longest solution uses the top path every time.}
        \label{fig:mazeMoves_long}    
    \end{subfigure}
\end{minipage}
\begin{minipage}{0.59\linewidth}
    \begin{subfigure}{1.0\linewidth}
        \centering
        \includegraphics[scale=1.0,page=1]{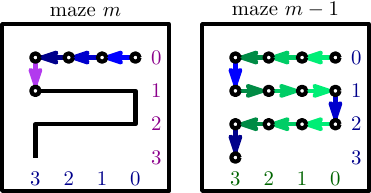}
        \caption{Change the leftmost digit from $0$ to $1$.}
        \label{fig:mazeMoves_1}    
    \end{subfigure}
    \begin{subfigure}{1.0\linewidth}
        \centering
        \includegraphics[scale=1.0,page=3]{mazeMoves.pdf}
        \caption{Change the leftmost digit from $1$ to $2$.}
        \label{fig:mazeMoves_2}    
    \end{subfigure}
    \begin{subfigure}{1.0\linewidth}
        \centering
        \includegraphics[scale=1.0,page=4]{mazeMoves.pdf}
        \caption{Change the leftmost digit from $2$ to $3$.}
        \label{fig:mazeMoves_3}    
    \end{subfigure}
\end{minipage}
\caption{Solving Ziggu puzzles requires moving back-and-forth through the $m$ mazes.
However, the top path shown in (a) only needs to be traversed once.
More specifically, it is traversed only the first time in (b) the shortest solution, and every time in (c) the longest solution.
To solve a Ziggu puzzle we must exit the last maze (i.e., maze $m$).
This involves changing the leftmost digit from (d) $0$ to $1$, then (e) $1$ to $2$, and finally (f) $2$ to $3$.
The shortest solution skips the back-and-forth traversal of the top path in maze $m-1$ in both (e) and (f).
Each omission skips over all but one of the states in maze $m-2$ (not shown).
}
\label{fig:mazeMoves}
\end{figure}

\begin{figure}
    \centering
    \includegraphics[width=1.0\linewidth]{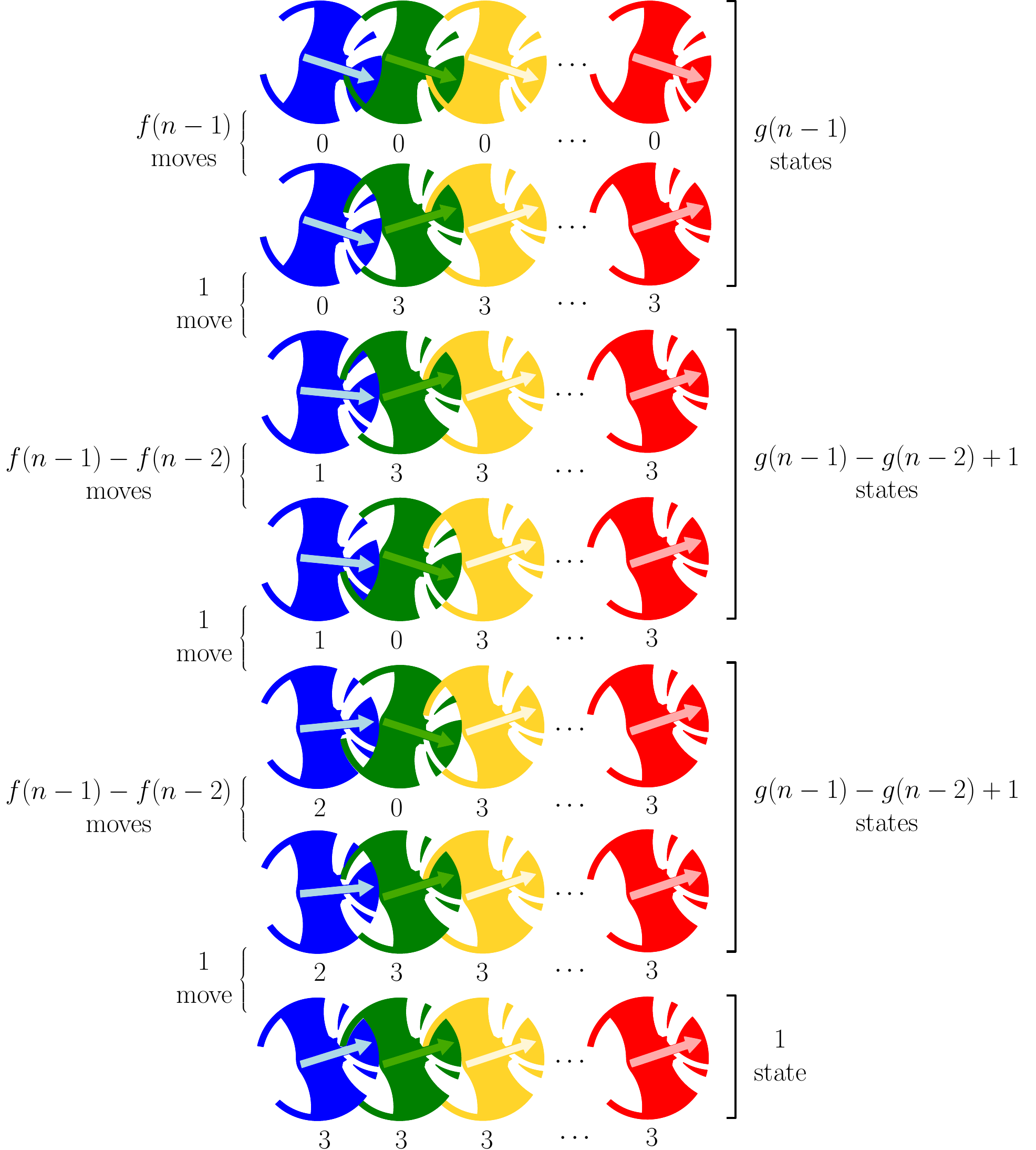}
    \caption{Recurrences for the unique shortest solution to Ziggu puzzles illustrated using Zigguhooked.
    Zigguhooked puzzles with $n$ discs require $f(n) = 3f(n-1) - 2f(n-2) + 3$ moves to solve with base cases $f(1) = 3$ and $f(2) = 12$.
    Alternatively, the shortest solution uses $g(n) = 3g(n-1) - 2g(n-2) + 3$ states with $g(1) = 4$ and $g(2) = 13$.
    In these recurrences, the negative terms come from taking the `shortcut' twice.
    When $n=3$ the shortcut can be seen in the six optional vertices in Figure \ref{fig:statesZiggu} and the additional gap of length six in the $\SHORT{3}$ column Table \ref{tab:orders}.
    }
    \label{fig:recurrenceZigguhookedShort}
\end{figure}

\begin{figure}
    \centering
    \includegraphics[width=0.71\linewidth]{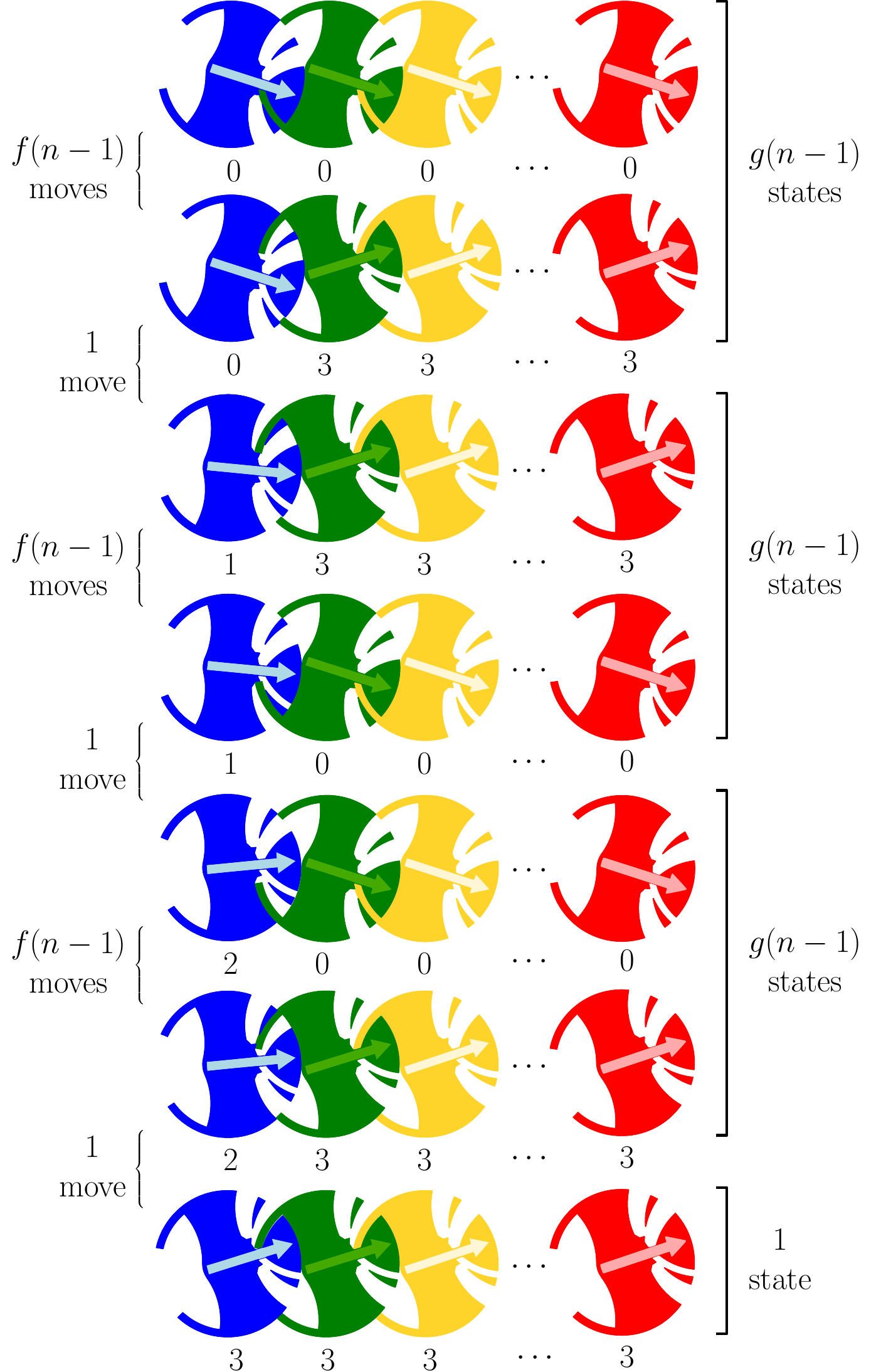}
    \caption{Recurrences for the unique longest solution to Ziggu puzzles illustrated using Zigguhooked.
    Zigguhooked puzzles with $n$ discs can be solved using $f(n) = 3f(n-1) + 3$ moves to solve with base cases $f(1) = 3$ and $f(2) = 12$.
    Alternatively, the shortest solution uses $g(n) = 3g(n-1) + 1$ states with $g(1) = 4$ and $g(2) = 13$.
    In these recurrences, no `shortcuts' are used.
    So the optional six vertices are used in Figure \ref{fig:statesZiggu} and the $\LONG{3}$ column of Table \ref{tab:orders} does not use the optional gap of length six.
    }
    \label{fig:recurrenceZigguhookedLong}
\end{figure}

Ziggu puzzles always have the same start state and target state, so these factors do not result in different recurrences.
However, these puzzles have different shortest and longest solutions, so we will still develop two different pairs of recurrences.
The principle behind our recurrences can be found in Figure \ref{fig:mazeMoves}.
More specifically, this figure shows how the shortest and longest solutions in a Ziggu maze puzzle can be expressed in terms of the solutions to smaller Ziggu maze puzzles.
Figure \ref{fig:recurrenceZigguhookedShort} and \ref{fig:recurrenceZigguhookedLong} then illustrates the recurrences for the shortest and longest solutions in Zigguhooked, respectively.
All four recurrences and summarized in Table \ref{tab:recurrences_Ziggu}.
Theorems \ref{thm:recurrenceZigguShortest}--\ref{thm:recurrenceZigguLongest} substitute in closed forms of the respective formulae.

\begin{theorem} \label{thm:recurrenceZigguShortest}
    The number of states in the shortest solution to Ziggu puzzles with $n$ digits is $f(n) = 6 \cdot 2^n - 3 \cdot n - 5$.
    The number of moves on the shortest solution is $g(n) = 6 \cdot 2^n - 3 \cdot n - 6$.
    These quantities count the number of ziggu strings between $0^n$ and $3^n$ in the quaternary Gray code.
\end{theorem}

\begin{theorem} \label{thm:recurrenceZigguLongest}
    The number of states in the longest solution to Ziggu puzzles with $n$ digits is $f(n) = $.
    The number of moves on the shortest solution is $g(n) = $.
    These quantities count the number of valid strings between $0^n$ and $3^n$ in the quaternary Gray code.
\end{theorem}

\section{Reflected Gray Codes: Binary and Quaternary}
\label{sec:reflected}

In this section we discuss the Gray codes that underlie Gray code puzzles and Ziggu puzzles.
The term \emph{Gray code} can broadly refer to a listing of some set of objects in which consecutive objects differ in a small amount, and in this section we'll consider a Gray code for $n$-bit binary strings, and a Gray code for $n$-digit quaternary strings.
Both Gray codes are based reflection, which reverses a list of strings (but not the individual strings themselves).
For example, $\reverse{123, 456, 789} = 789, 456, 123$.

When ordering numbers or strings, Gray codes are not the same as standard numeric or lexicographic order.
For example, in standard numeric order the $8$-bit binary string $01111111$ is followed by $10000000$ (i.e., all $n$ bits are flipped) whereas it is followed by $01111101$ (i.e., only the second rightmost bit is flipped) in the Gray code in Section \ref{sec:reflected_binary}.

\subsection{Binary Reflected Gray Code}
\label{sec:reflected_binary}

The most well-known Gray code is the \emph{binary reflected Gray code (BRGC)} or simply the \emph{Gray code}.
It is an ordering of all $2^n$ binary strings of length $n$ in which successive strings differ in only a single bit.
For example, the orders for $n=2$ and $n=3$ are provided below, where the overlines mark the bits that change to create the next string.
\begin{align}
    \BRGC{2} &= 0\overline{0},\overline{0}1,1\overline{1},10 \label{eq:BRGC2} \\
    \BRGC{3} &= 00\overline{0},0\overline{0}1,01\overline{1},\overline{0}10,11\overline{0},1\overline{1}1,10\overline{1},100. \label{eq:BRGC3}
\end{align}

The global structure of the binary reflected Gray code is given by the following recursive formulae.
\begin{equation} \label{eq:BRGC}
    \BRGC{1}=0,1\qquad\BRGC{n}=0\cdot\BRGC{n-1},1\cdot\reverse{\BRGC{n-1}}
\end{equation}
In this formula, the $\cdot$ denotes concatenation, the commas denote appending one list to another, and $\reverse{}$ denotes the list in reflected order. 
For example, the first four words of \eqref{eq:BRGC3} are identical to \eqref{eq:BRGC2} except that $0$ is prefixed to them.
Similarly, the last four words of \eqref{eq:BRGC3} are identical to \eqref{eq:BRGC2} written in reflected order and with $1$ prefixed.
The formula \eqref{eq:BRGC} and its resulting order are usually attributed to Frank Gray \cite{Fra1953}, but their histories go back much further \cite{heath1972origins}. 

Its worth noting that the Gray code can also be created using a simple greedy algorithm \cite{Wil2013}.
Start a list with the all-zeros string $0^n$ (where exponentiation denotes repetition) and then repeatedly create the next string in the list from the most recently added string as follows:
complement the rightmost bit that creates a new string.
For example, when $n=3$ the list will begin $000, 001, 011, 010$.
To create the next string we try complementing the bits in $010$ from right-to-left: 
$01\overline{0} = 011$ is already in the list;
$0\overline{1}0 = 000$ is already in the list;
$\overline{0}10 = 110$ is not in the list so it is added to the end of the list.
This notion of not recreating a previous string is similar to the puzzle-solving approach from Section \ref{sec:solving_resuming}.

\subsection{Quaternary Reflected Gray Code}
\label{sec:reflected_quaternary}

We define the quaternary reflected Gray code recursively as follows. Let $\QUAT{1} = 0,1,2,3$, and for $n\geq 2$, let
\begin{align}
   \QUAT{n} &= \begin{cases}
       0 \cdot \QUAT{n-1}, \\
       1 \cdot \reverse{\QUAT{n-1}}, \\
       2 \cdot \QUAT{n-1}, \\
       3 \cdot \reverse{\QUAT{n-1}} \\
   \end{cases}\label{eq:rec_quat} 
\end{align}

Note that the global structure of the quaternary reflected Gray code is very similar to the structure of the binary reflected Gray code, as both formulas involve prepending digits to the list of one size smaller, with every second copy in reverse order.

\section{Recursive Formulae: The Global Structure of Solutions}
\label{sec:recursive}


Understanding how the solution of a puzzle changes as the number of mazes in the puzzle changes requires knowledge of the global structure of the solution. We describe the global structure of a sequence using a \textit{recursive formula}, where the formula for a given $n$ is built up from the sequences for smaller $n$.

\subsection{Recursive Structure of Gray Code Puzzle Solutions}
\label{sec:recursive_BRGC}

\begin{figure}
    \centering
    \begin{subfigure}{\textwidth}
        \centering
        \includegraphics[width=0.97\textwidth,page=2]{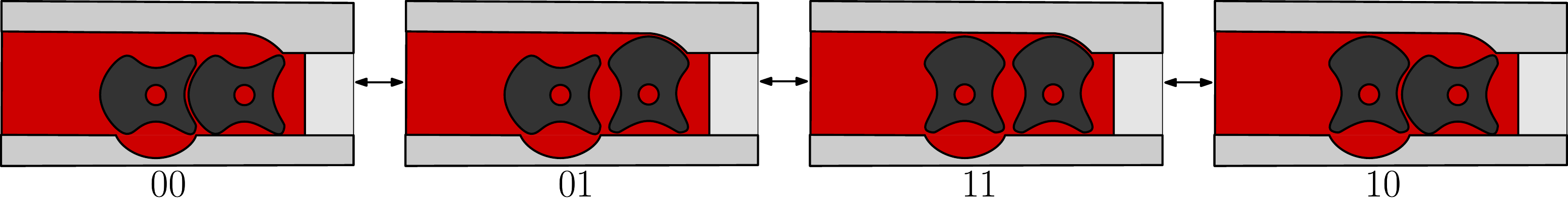} 
        \caption{The states of \textit{Spin-Out} with two dials.}
        \label{fig:spinout-rec-2}
    \end{subfigure}
    \vfill
    \begin{subfigure}{\textwidth}
        \centering
        \includegraphics[width=0.97\textwidth,page=3]{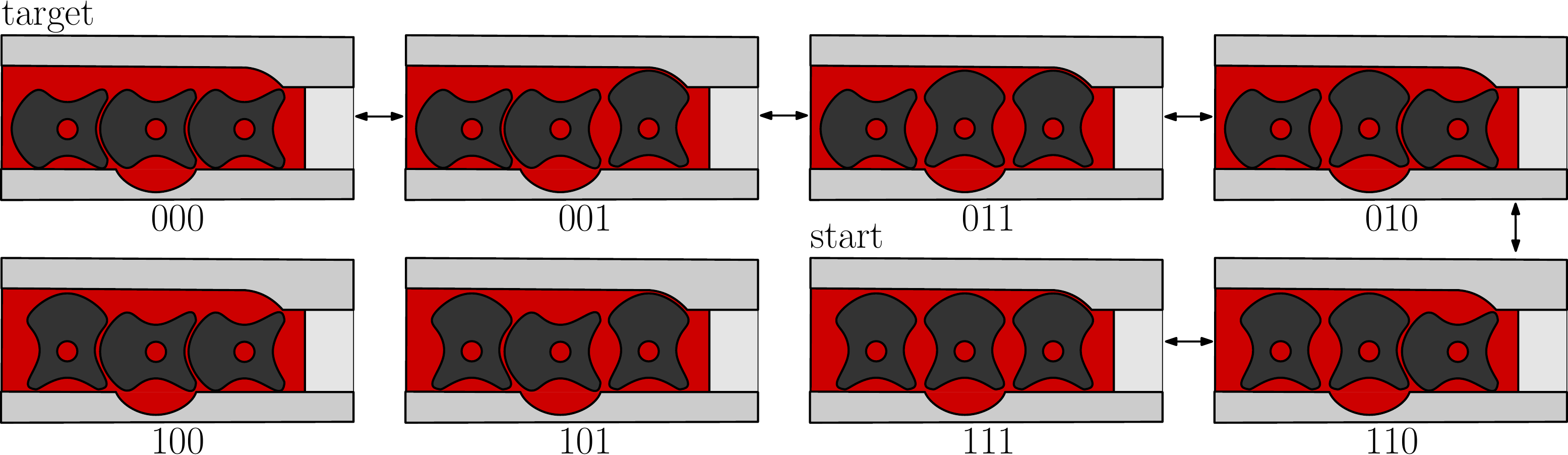} 
        \caption{The states of \textit{Spin-Out} with three dials.}
        \label{fig:spinout-rec-3}
    \end{subfigure}
    \caption{This figure illustrates the recursive structure of a modified version of \textit{Spin-Out} that has only three dials. Figure \ref{fig:spinout-rec-2} illustrates $\BRGC{2}$. Below each $2$-dial state in Figure \ref{fig:spinout-rec-3} is the corresponding $3$-dial state with either a $0$ or $1$ prepended, depending on whether the $3$-dial state is in the second or third row of the figure, respectively. }
    \label{fig:SpinOut-recurisve}
\end{figure}

The state space of many Gray code puzzles is precisely the binary reflected Gray code, whose recursive description is given by \eqref{eq:BRGC}. \textit{Spin-Out} is one such Gray code puzzle, and its recursive structure is illustrated in Figure \ref{fig:SpinOut-recurisve}.

\subsection{Recursive Structure of the Longest Ziggu Puzzle Solution}
\label{sec:recursive_long}

Let $\LONG{1} = 0,1,2,3$, and for $n\geq 2$, let
\begin{align}
   \LONG{n} &= \begin{cases}
       0 \cdot \LONG{n-1}, \\
       1 \cdot \reverse{\LONG{n-1}}, \\
       2 \cdot \LONG{n-1}, \\
       3^n.
   \end{cases}\label{eq:rec_long}
\end{align}

This recursive formula is derived from the recursive formula for the quaternary reflected Gray code, with the modification that instead of recursively including the entire list $3\cdot\reverse{\QUAT{n-1}}$, we append the state $3^n$, as there is no need to continue beyond the solved state $3^n$.

\begin{figure}
    \centering
    \begin{subfigure}{\textwidth}
        \centering
        \includegraphics[width=0.97\textwidth,page=36]{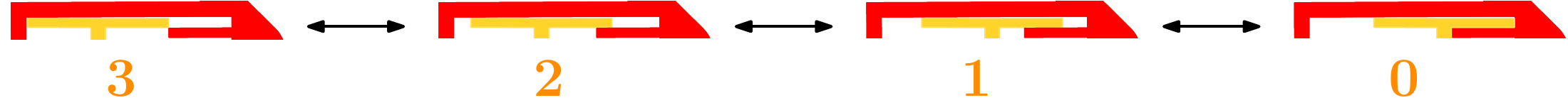} 
        \caption{The states of Zigguflat with only the horizontal axis of a single maze.}
        \label{fig:zf-rec-1}
    \end{subfigure}
    \vfill
    \begin{subfigure}{\textwidth}
        \centering
        \includegraphics[width=0.97\textwidth,page=37]{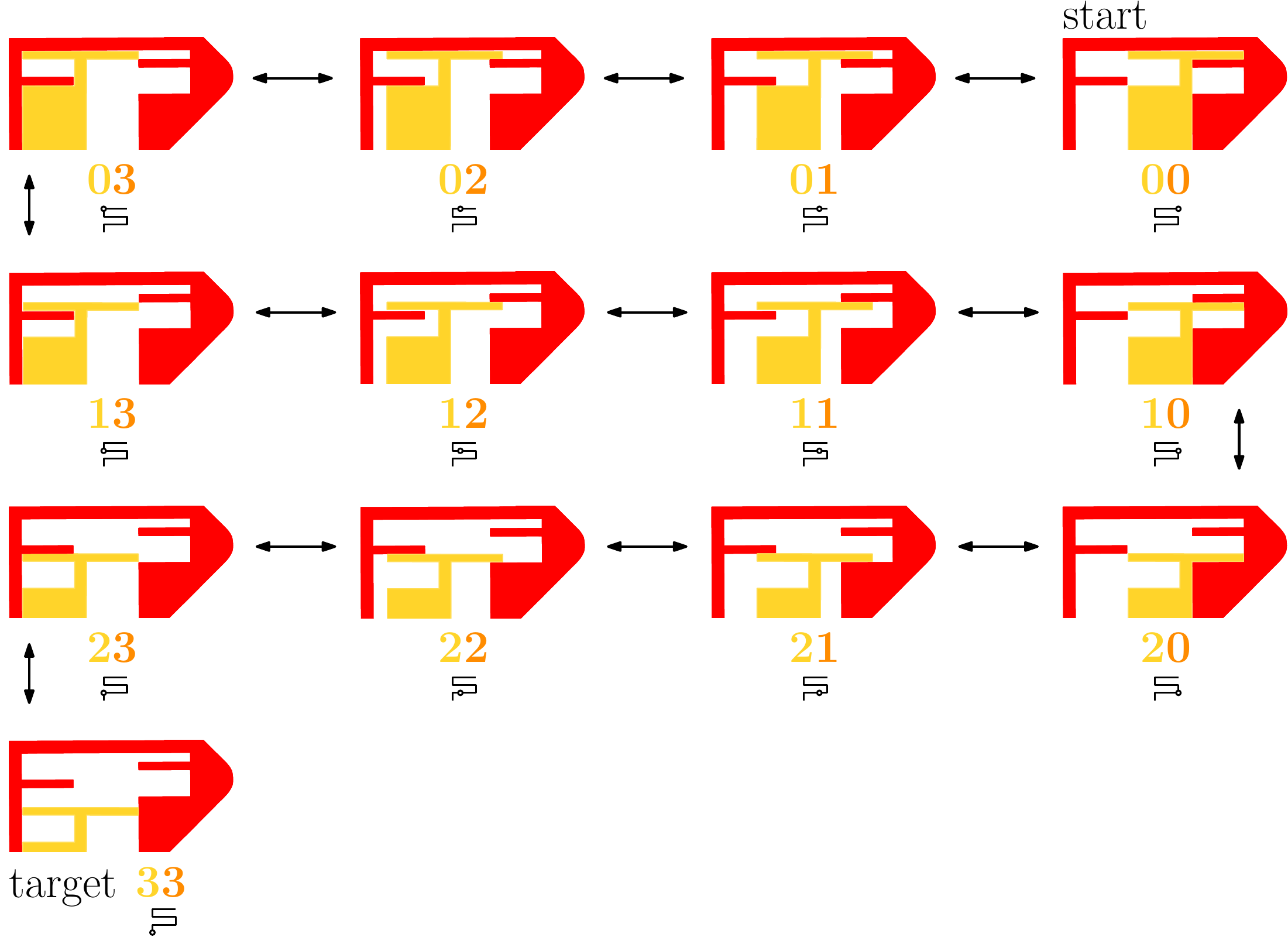} 
        \caption{The states of Zigguflat with one maze.}
        \label{fig:zf-rec-2}
    \end{subfigure}
    \caption{This figure illustrates the recursive structure of both the longest and shortest solution to the Ziggu puzzles (they agree when $n=1,2$) using a smaller version of Zigguflat with only one maze. Figure \ref{fig:zf-rec-1} is the solution with only the horizontal part of the maze, which has 4 states. The successive rows in Figure \ref{fig:zf-rec-2} are created by prepending either a $0$, $1$, or $2$ in each successive row to the state given by the column.}
    \label{fig:zf-recursive}
\end{figure}

\begin{figure}
    \centering
    \begin{subfigure}{\textwidth}
        \centering
        \includegraphics[width=0.97\textwidth,page=22]{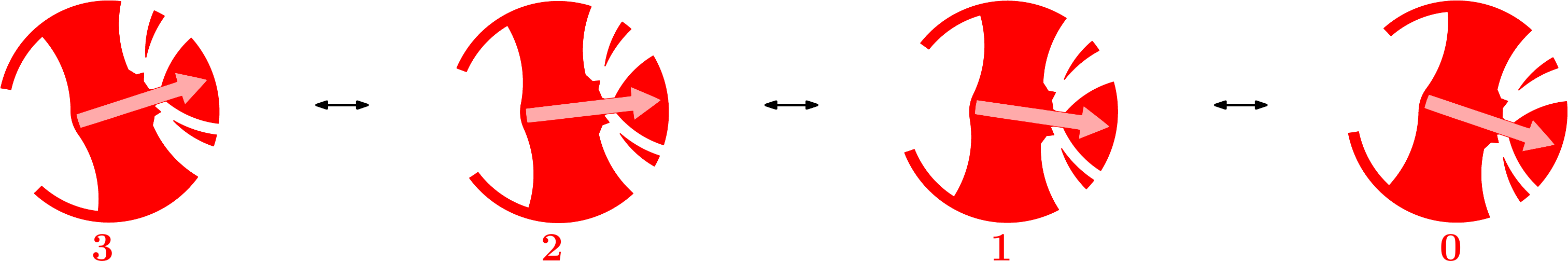} 
        \caption{The states of Zigguhooked with one dial.}
        \label{fig:zigguhooked-rec-1}
    \end{subfigure}
    \vfill
    \begin{subfigure}{\textwidth}
        \centering
        \includegraphics[width=0.97\textwidth,page=23]{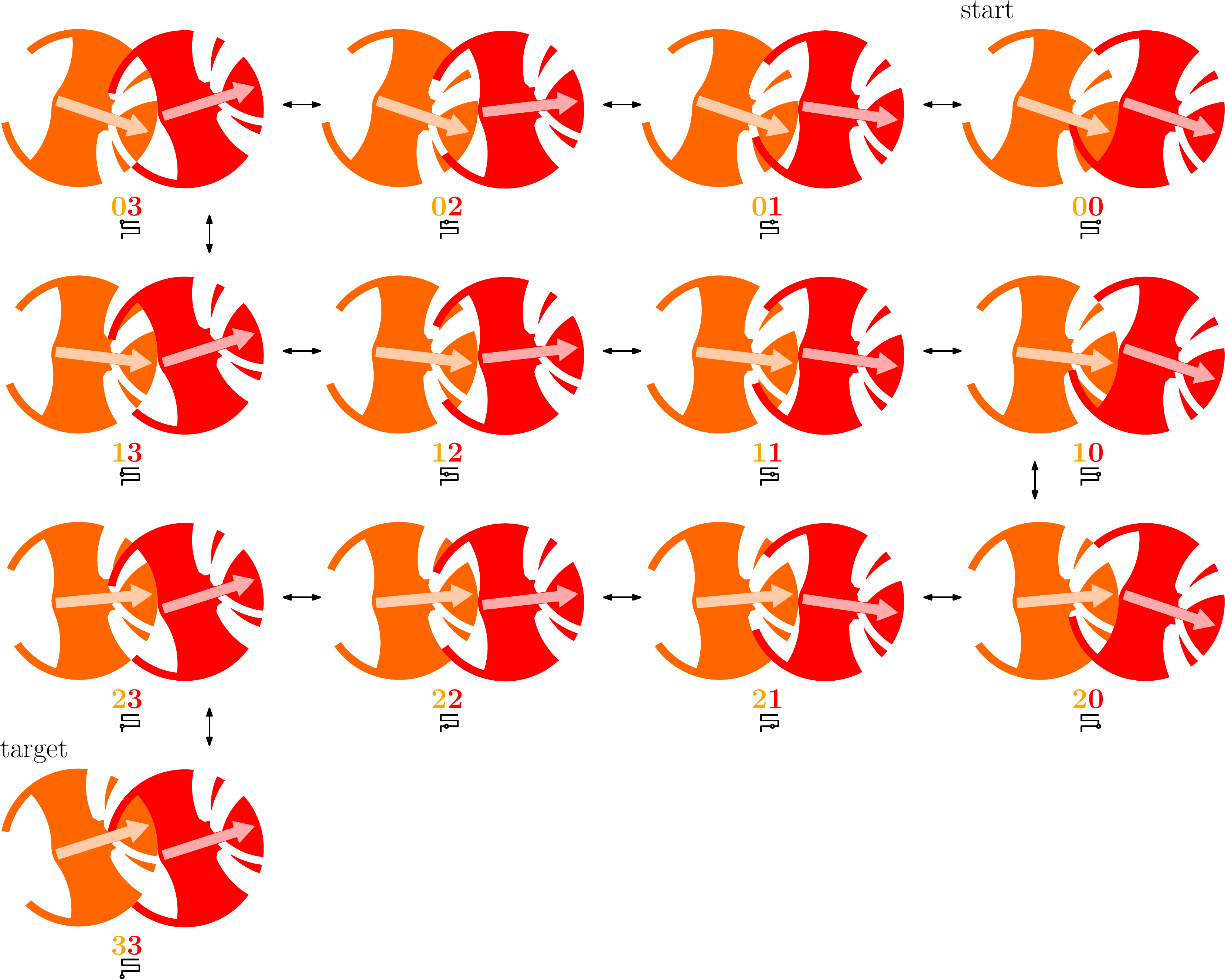} 
        \caption{The states of Zigguhooked with two dials.}
        \label{fig:zigguhooked-rec-2}
    \end{subfigure}
    \caption{This figure illustrates the recursive structure of both the longest and shortest solution to the Ziggu puzzles (they agree when $n=1,2$) using a smaller version of Zigguhooked with only 2 dials. Figure \ref{fig:zigguhooked-rec-1} is the solution with only 1 dial, which has 4 states. The successive rows in Figure \ref{fig:zigguhooked-rec-2} are created by prepending either a $0$, $1$, or $2$ in each successive row to the state given by the column.}
    \label{fig:hooked-recursive}
\end{figure}

The recursive structures of Zigguflat and Zigguhooked are illustrated in \ref{fig:zf-recursive} and \ref{fig:hooked-recursive}, respectively.

\subsection{Recursive Structure of the Shortest Ziggu Puzzle Solution}
\label{sec:recursive_short}

Let $\SHORT{1}=\Core(1)=0,1,2,3$, and for $n\geq2$, let
\begin{align}
    \SHORT{n}    & = \begin{cases}  0\cdot\SHORT{n-1}, \\
                     \Core(n)[1:]
                 \end{cases}\label{eq:rec_short}             \\[-0.2em]
    \Core(n) & = \begin{cases} 03^{n-1},                   \\[-0.2em]
                     1\cdot\reverse{\Core(n-1)}, \\
                     2\cdot\Core(n-1),           \\[-0.2em]
                     3^n.
                 \end{cases}\nonumber
\end{align}
The $\Core(n)[1:]$ denotes all but the first element of $\Core(n)$. This ``list slicing'' is needed because the word $03^{n-1}$ appears as both the last word of $0\cdot\SHORT{n-1}$ and the first word of $\Core(n).$ $\SHORT{n}$ is the sublist of $\LONG{n}$ obtained by omitting all states where $03$ is followed by a digit that is not $3$.

\section{Integer Sequences: The Changes}
\label{sec:sequences}

Experienced sequential puzzle solvers do not focus on successive puzzle states, but rather the changes to perform to create the states.
For Gray code puzzles, these changes fall into a familiar rhythmic pattern related to counting in binary.
\begin{equation} \label{eq:rulerBinary5}
    \rulerBinary{} = 1,2,1,3,1,2,1,4,1,2,1,3,1,2,1,5,1,2,1,3,1,2,1,4,1,2,1,3,1,2,1,{.}{.}{.}
\end{equation}
In this section we review the change sequence for Gray code puzzles, and then consider analogous change sequences for longest and shortest solutions to Ziggu puzzles.
We'll see that Ziggu puzzle solutions follow a pattern for counting in quaternary (i.e., base $4$), except that portions of the sequence are skipped over.
\begin{equation} \label{eq:rulerQuat3}
    \small
    \rulerQuat{} = 
    1,1,1,2,1,1,1,2,1,1,1,2,1,1,1,3,1,1,1,2,1,1,1,2,1,1,1,2,\hcancel[red]{1,1,1,}3,\hcancel[red]{1,1,1,} 
    {.}{.}{.}
\end{equation}

\begin{figure}
    \makebox[\textwidth][c]{\includegraphics[width=1.2\textwidth]{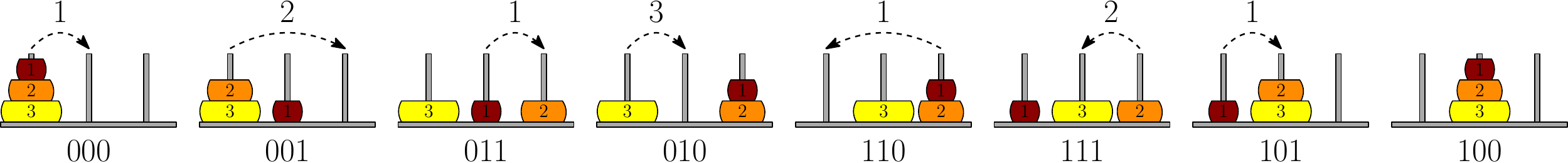}}
    \caption{Solving the \emph{Towers of Hanoi} with $n=3$ discs.
    The solution moves the discs using binary ruler sequence $\rulerBinary{3} = 1, 2, 1, 3, 1, 2, 1$.
    Each disc has a unique peg onto which it can be moved, except for disc $1$ which always moves one peg to the right (cyclically).
    The signed sequence $\srulerBinary{3} = 1, 2, -1, 3, 1, -2, -1$ also specifies whether the disc $i$ moves off of or onto disc $i+1$ (or the first peg for the largest disc).}
    \label{fig:hanoi-3}
\end{figure}

\subsection{Sequences for Gray Code Puzzles}
\label{sec:sequences_BRGC}


The change sequence for Gray code puzzles is both well-known and intimately related to other concepts including binary counting and the binary reflected Gray code.
Here we review its unsigned and signed versions, and relate them to the process of solving Gray code puzzles.

\subsubsection{Binary Ruler Sequence}
\label{sec:sequences_BRGC_ruler}

The integer sequence in \eqref{eq:rulerBinary5} is known as the \emph{binary ruler sequence} \OEIS{A001511}.
Its values give the sequence of piece numbers to modify in a Gray code puzzle%
\footnote{Gray code puzzles that don't run from $0^n$ to $10^{n-1}$ will use a subsequence of some $\rulerBinary{n}$ (or its reversal). For example, Spin-Out's solution begins $3,1,2,1,4,1,2,1,3,\ldots$.}.
For example, a \emph{Towers of Hanoi} speedrunner will move disc $1$, then disc $2$, disc $1$, disc $3$, and so on%
\footnote{There is always a unique peg onto which a disc can be moved, except for disc $1$.
The solver should always move this smallest disc in a consistent direction (i.e., to the right with wrap-around).}  
according to \eqref{eq:rulerBinary5}.
The puzzle with $n$ discs is solved immediately before the first copy of $n+1$ in the sequence. 
There are $2^{n}-1$ values in the sequence before the first copy of $n+1$, and this initial part of the sequence is denoted $\rulerBinary{n}$.

\begin{figure}
    \centering
    \includegraphics[width=1.0\linewidth]{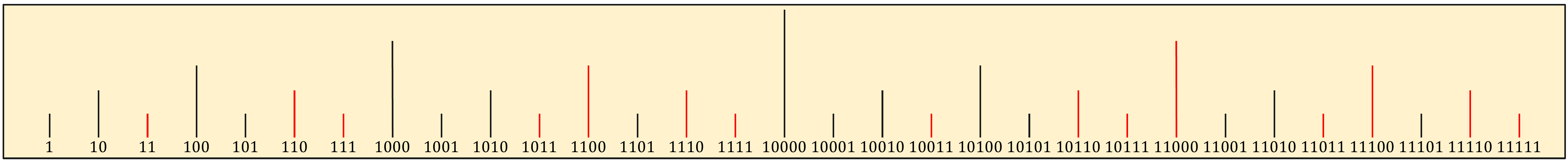}
    \caption{A binary ruler with $n=5$ bits (leading $0$s omitted). 
    The height of each tick is one more than the number of $0$s at the end of the string (e.g., $101\underline{00}$ has height $3$) or equivalently, one more than the largest power of two that divides it (e.g., $10100$ is divisible by $2^2 = 4$). 
    The sequence of heights is the binary ruler sequence $\rulerBinary{5}$ seen in \eqref{eq:rulerBinary5}.
    The signed ruler sequence $\srulerBinary{5}$ \eqref{eq:srulerBinary5} negates the red ticks, which occur when the aforementioned quotient is $3 \mod 4$, or equivalently when the zeroes are preceded by two $1$s (e.g., $10\underline{11}0$ is red but $10\underline{01}0$ is not).
    The sequences also specify the bit changes in the binary reflected Gray code and the pieces to change when solving Gray code puzzles.
    }
    \label{fig:srulerBinary5}
\end{figure}

The name of the binary ruler sequence comes from the tick marks on a base-two ruler (see Figure \ref{fig:srulerBinary5}).
Each successive value gives the number of bits that change in lexicographic order, or equivalently,%
\footnote{The role of the ruler sequence in these two orders can be seen in Oskar van Deventer's \emph{Gray Code Counter}.
} the index of the single bit that changes in $b_n b_{n-1} \cdots b_1$ in the binary reflected Gray Code, as seen below for $n=5$.
\begin{align}
    \lexBinary{5} &= 0000\overline{0}, 000\overline{01}, 0001\overline{0}, 00\overline{011}, 0010\overline{0}, 001\overline{01}, \ldots \\
    \BRGC{5} &= 0000\overline{0}, 000\overline{0}1, 0001\overline{1}, 00\overline{0}10, 0011\overline{0}, 001\overline{1}1, \ldots
\end{align}

The binary ruler sequence can be defined in several ways.
Most directly, the $j^{\text{th}}$ entry is one more than the largest power of two that divides $j$, or equivalently, one more than the number of consecutive $0$s that appears at the end of its binary representation.
For example, the $152^{\text{nd}}$ entry of the sequence is $4$ due to the fact that $152$ is divisible by $2^3$ (but not $2^4$), or equivalently, $10011000$ ends with $3$ copies of $0$.  
The sequence can also be defined recursively with base case $\rulerBinary{1}=1$ and either inductive case below.
\begin{align} 
    \rulerBinary{n} &= \rulerBinary{n-1}, n, \reverse{\rulerBinary{n-1}} \label{eq:rulerRec_reverse}\\
    &= \rulerBinary{n-1}, n, \rulerBinary{n-1}\label{eq:rulerRec_norev}
\end{align}
In other words, the sequence for $n$-bits is comprised of two copies of the sequence for $(n-1)$-bits with a single copy of $n$ between them.
For example, \eqref{eq:rulerRec_norev} gives $\rulerBinary{3} = \rulerBinary{2}, 3, \rulerBinary{2} = 1,2,1,3,1,2,1$.
In \eqref{eq:rulerRec_reverse} the $\reverse{}$ denotes reversing the second subsequence, and this definition aligns with the standard recursive definition of the binary \underline{reflected} Gray code.
But the reversal is redundant here as $\rulerBinary{n}$ is always a palindrome.
Nevertheless, \eqref{eq:rulerRec_reverse} is helpful for understanding related sequences.

\subsubsection{Signed Binary Ruler Sequence}
\label{sec:sequences_BRGC_sruler}

The \emph{signed binary ruler sequence} \OEIS{A164677} is the same as the binary ruler sequence except that some symbols are negated. 
\begin{equation} \label{eq:srulerBinary5}
    \srulerBinary{} = 1, 2, -1, 3, 1, -2, -1, 4, 1, 2, -1, -3, 1, -2, -1, 5, 1, 2, -1, \ldots
\end{equation}
The signs provide the direction in which the bits change in the binary reflected Gray code.
More specifically, if the $j$th entry is $i$, then bit $b_i$ is changed from $0$ to $1$, whereas an entry of ${-}i$ means that bit $b_i$ is changed from $1$ to $0$.
The sequence is defined directly in the same way, except that the sign is positive or negative when the quotient is $1$ or $3$ modulo $4$, respectively, and this is equivalent to whether the trailing $0$s in the binary representation are preceded by $01$ or $11$, respectively.
For example, the $152^{\text{nd}}$ entry of the signed sequence is $-4$ since $152 / 2^3 = 19 \equiv 3 \mod{4}$, or equivalently, the trailing $0$s in $100\underline{11}000$ are preceded by $11$.

For Gray code puzzle solvers, the additional information provided in the signed sequence is usually not helpful.
This is because the way in which an individual piece can be changed is readily apparent.
For example, in \emph{The Brain} there is no question about whether a wedge can be pulled out or pushed in.
However, the additional information can be useful in certain puzzle solving situations (e.g., blindfolded speedruns).

When using the binary reflected Gray code as inspiration, we can construct the following recurrence using $\srulerBinary{1}=1$ and
\begin{equation} \label{eq:srulerRec_reverse}
    \srulerBinary{n} = \srulerBinary{n-1}, n, \reverse{\complement{\srulerBinary{n-1}}} 
\end{equation}
where the overline represents complementation (i.e., negation of every value in the subsequence).
For example, $\srulerBinary{3} = 1,2,-1,3,1,-2,-1$ where the first subsequence $\srulerBinary{2} = 1,2,-1$ is `undone' by reversing and complementing it to create the second subsequence $\reverse{\complement{1,2,-1}} = 1,-2,-1$.
At first glance, it may seem difficult to make the type of simplification that we saw for the unsigned sequence in \eqref{eq:rulerRec_reverse}--\eqref{eq:rulerRec_norev} because of the pattern of negations.
Fortunately, each value in the binary ruler sequence $\rulerBinary{n}$ appears an even number of times, except for the single value of $n$.
As a result, the two subsequences in \eqref{eq:srulerRec_reverse} will have identical negation patterns except that $n-1$ will be positive in the first subsequence and negative in the second subsequence.
Since the largest value in a binary ruler sequence exactly once in the middle of it, we can construct the following alternate recursive formula
\begin{equation} \label{eq:srulerRec_norev}
    \srulerBinary{n} = \srulerBinary{n-1}, n, \complementMid{\srulerBinary{n-1}} 
\end{equation}
where $\complementMid{}$ complements the middle value of an odd-length sequence.

There is also a simpler way of viewing the signed binary ruler sequence: negate every second copy of each value in the unsigned binary ruler sequence.
In other words, puzzle solvers simply need to remember the parity of how many times each value has been seen.


\subsection{Sequence for Quaternary Gray Code}
\label{sec:ruler_GC}

\begin{figure}
    \centering
    \makebox[\textwidth][c]{\includegraphics[width=1.2\textwidth]{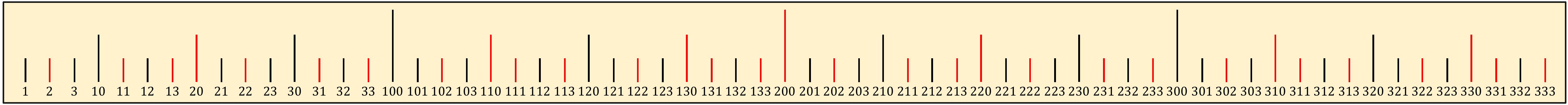}}
    \makebox[\textwidth][c]{\includegraphics[width=1.2\textwidth]{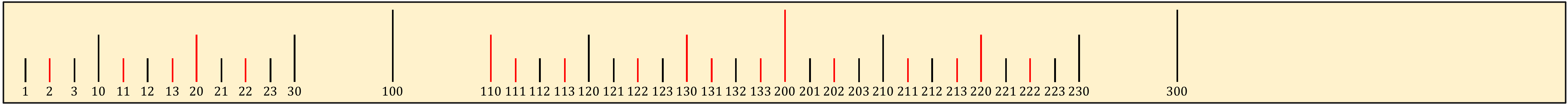}}
    \makebox[\textwidth][c]{\includegraphics[width=1.2\textwidth]{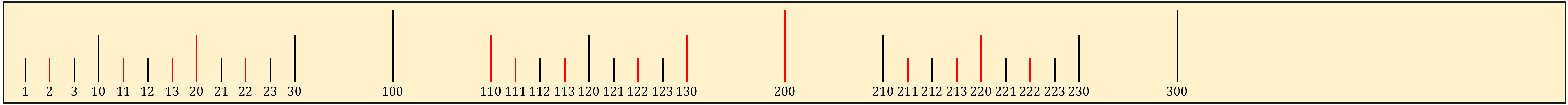}}
    \caption{The quaternary rulers with $n=3$ digits (leading $0$s omitted). 
    The height of each tick is one more than the number of $0$s at the end of the string (e.g., $13\underline{0}$ has height $2$) or one more than the largest power of four that divides it (e.g., $130$ is $28$ in decimal which is divisible by $4^1 = 4$).
    The sequence of heights is the quaternary ruler sequence $\rulerQuat{3}$ seen in \eqref{eq:rulerQuat3}.
    The signed version $\srulerQuat{3}$ \eqref{eq:rulerQuat_norev} negates the red ticks, which occur when the aforementioned quotient is $3 \mod 4$, or equivalently, when the zeroes are preceded by two $1$s (e.g., $10\underline{11}0$ is red but $10\underline{01}0$ is not).
    The sequences also specify the digit changes in the quaternary reflected Gray code, and subsequences specify the pieces to change when solving Ziggu puzzles using both the longest and shortest solutions.
    }
    \label{fig:srulerQuat3}
\end{figure}

\subsubsection{Quaternary Ruler Sequence}
\label{sec:sequences_QRGC_ruler}

We can define the quaternary ruler sequence $\rulerQuat{n}$ recursively, with $\rulerQuat{1}=1,1,1$ and
\begin{align}
\rulerQuat{n} 
&= \rulerQuat{n-1}, n, \reverse{\rulerQuat{n-1}},  n, \rulerQuat{n-1}, n, \reverse{\rulerQuat{n-1}} \label{eq:rulerQuat_reverse}  \\
&= \rulerQuat{n-1}, n, \rulerQuat{n-1},  n, \rulerQuat{n-1}, n, \rulerQuat{n-1} \label{eq:rulerQuat_norev}
\end{align}
for all $n>1$. We have two equivalent definitions for $\rulerQuat{n}$ since it is a palindrome, implying that $\reverse{\rulerQuat{n-1}}=\rulerQuat{n}$. Note that \eqref{eq:rulerQuat_norev} matches the recursive definition of $\QUAT{n}$ given in \eqref{eq:rec_quat}. 

As an example, \eqref{eq:rulerQuat_norev} gives that 
\begin{align*}
\small
    \rulerQuat{3} =&\ \rulerQuat{2},3,\rulerQuat{2},3,\rulerQuat{2},3,\rulerQuat{2} \\
    =&\ 1,1,1,2,1,1,1,2,1,1,1,2,1,1,1,3,1,1,1,2,1,1,1,2,1,1,1,2,1,1,1,3, \\
    &\ 1,1,1,2,1,1,1,2,1,1,1,2,1,1,1,3,1,1,1,2,1,1,1,2,1,1,1,2,1,1,1,3
\end{align*}
This change sequence is shown in the $\rulerQuat{3}$ column of Table \ref{tab:orders}, where it can be compared with the states in $\QUAT{3}$.

\subsubsection{Signed Quaternary Ruler Sequence}
\label{sec:sequences_QRGC_sruler}

Just as for the binary reflected Gray code, we can define a signed quaternary ruler sequence $\srulerQuat{n}$, where $i$ and $-i$ indicate that digit $q_i$ should be incremented or decremented, respectively. Recall that we label the digits in a quaternary word $q=q_n\cdots q_2q_1$. We can define $\srulerQuat{n}$ recursively by letting $\srulerQuat{1}=1,1,1$ and defining
\begin{align}
\srulerQuat{n} 
&= \srulerQuat{n-1}, n, \reverse{\complement{\srulerQuat{n-1}}},  n, \srulerQuat{n-1}, n, \reverse{\complement{\srulerQuat{n-1}}} \label{eq:srulerQuat_reverse}  \\
&= \srulerQuat{n-1}, n, \complementMax{\srulerQuat{n-1}},  n, \srulerQuat{n-1}, n, \complementMax{\srulerQuat{n-1}} \label{eq:srulerQuat_norev}
\end{align}
for all $n>1$. Here the function $\complementMax{}$ complements the sign of the largest indices in the sequence. As an example, 
\begin{align*}
    \complementMax{\srulerQuat{2}} &= \complementMax{1,1,1,2,-1,-1,-1,2,1,1,1,2,-1,-1,-1} \\
    &= 1,1,1,-2,-1,-1,-1,-2,1,1,1,-2,-1,-1,-1, 
\end{align*}
so $\complementMax{}$ performed $\complement{2}=-2$ since $2$ is the maximum entry found in $\srulerQuat{}.$ Therefore 
\begin{align*}
\small
    \srulerQuat{3} =&\ \srulerQuat{2},3,\complementMax{\srulerQuat{2}},3,\rulerQuat{2},3,\complementMax{\srulerQuat{2}} \\
    =&\ 1,1,1,2,-1,-1,-1,2,1,1,1,2,-1,-1,-1,3,1,1,1,-2,-1,-1,-1, \\
    &\ -2,1,1,1,-2,-1,-1,-1,3,1,1,1,2,-1,-1,-1,2,1,1,1,2,-1,-1, \\
    &\ -1,3,1,1,1,-2,-1,-1,-1,-2,1,1,1,-2,-1,-1,-1
\end{align*}

Observe that the $\complementMid{}$ function in \eqref{eq:srulerRec_norev} could be replaced with $\complementMax{}$ since the maximum index appears only once in $\srulerBinary{}$ at precisely the middle of the sequence. This highlights the similarity between \eqref{eq:srulerQuat_norev} and \eqref{eq:srulerRec_norev}. The signed sequence $\srulerQuat{3}$ is also included as a column in Table \ref{tab:orders}.

\subsection{Sequence for the Longest Solution to Ziggu Puzzles}
\label{sec:ruler_long}

In this section, we construct signed and unsigned change sequences for the longest solution to the Ziggu puzzles. As before, we represent a state in the longest solution to a Ziggu puzzle by $q=q_n\cdots q_2q_1\in\LONG{n}$. The change sequence for the longest solution is a sequence of the indices to change; i.e., if $i$ appears in the change sequence, we should change digit $q_i$ to obtain the next state in the longest solution. 

\subsubsection{Change Sequence for the Longest Solution}
\label{sec:sequences_long_ruler}

We can define the (unsigned) change sequence for the longest solution to the Ziggu puzzles $\rulerLong{n}$ recursively, with $\rulerLong{1}=1,1,1$ and
\begin{equation}\label{eq:rulerLong} 
\rulerLong{n} = \rulerLong{n-1}, n, \reverse{\rulerLong{n-1}},  n, \rulerLong{n-1}, n 
\end{equation}

for all $n>1$. Recall that the longest solution is the sublist of the quaternary Gray code obtained by omitting all states where a $3$ is followed by a digit that is not a $3$. This motivates the definition of $\rulerLong{n}$, as we obtain \eqref{eq:rulerLong} by omitting the last term $\reverse{\rulerQuat{n-1}}$ from \eqref{eq:rulerQuat_reverse}. Since this final term was omitted, $\rulerLong{n}$ is not a palindrome, and so we require the reversal of the middle $\rulerLong{n-1}$ term. The structure of $\rulerLong{n}$ matches the recursive definition of the Gray code $\LONG{n}$ given in \eqref{eq:rec_long}.

As an example, 
\begin{align}
\small
    \rulerLong{3} =&\ \rulerLong{2},3,\reverse{\rulerLong{2}},3,\rulerLong{2},3 \nonumber\\
    =&\  1,1,1,2,1,1,1,2,1,1,1,2,3,2,1,1,1,2,1,1,1,\nonumber\\
    &\ 2,1,1,1,3,1,1,1,2,1,1,1,2,1,1,1,2,3\nonumber\\
    =&\ 1,1,1,2,1,1,1,2,1,1,1,2,\hcancel[red]{1,1,1,}3,\hcancel[red]{1,1,1,}2,1,1,1,2,1,1,1,2,1,1,1,3, \label{eq:rulerLong_strike}\\
    &\ 1,1,1,2,1,1,1,2,1,1,1,2,\hcancel[red]{1,1,1,}3,\hcancel[red]{1,1,1,2,1,1,1,2,1,1,1,2,1,1,1,3}\nonumber
\end{align}
where \eqref{eq:rulerLong_strike} is $\rulerQuat{3}$ with certain entries omitted. The omitted entries correspond to the states in Table \ref{tab:orders} that are present in the $\QUAT{3}$ column but omitted in the $\LONG{3}$ column.

\begin{remark}\label{rem:no_signs_needed}
At first, this sequence seems to be omitting the information about whether to increment or decrement digit $q_i$. If one is solving a physical Ziggu puzzle, there is no ambiguity, as the pieces of the puzzle can only be manipulated in one of those two suggested ways. If we are working with the base $4$ representations, we still see that there is no ambiguity. If digit $q_i$ is a $0$ or $3$, then it can only be incremented or decremented, respectively. If $q_i$ is $1$ or $2$, then we increment $q_i$ if $\sum_{j=i+1}^n q_j$, the sum of the digits to the left of $q_i$, is odd, and decrement if the sum is even. If $q_i$ is not the rightmost digit, then we can equivalently consider the parity of $q_i+q_{i-1}$, the sum of $q_i$ and the digit immediately to its right. This fact is exploited in Section \ref{sec:loopless_short} to generate the shortest solution to the Ziggu puzzles in Algorithm \ref{alg:fast-ziggu}. In summary, one can always determine whether to increment or decrement a digit, even if we aren't working with the physical puzzle.
\end{remark}

Nevertheless, we can still construct a change sequence with added information about whether to increment or decrement a given digit. 

\subsubsection{Signed Change Sequence for the Longest Solution}
\label{sec:sequences_long_sruler}

Recall that given state $q=q_n\cdots q_2q_1\in\LONG{n}$, an entry $i$ or $-i$ in a change sequence instructs one to increment or digit $q_i$, respectively. We define the signed change sequence for the longest solution recursively by letting $\srulerLong{1}=1,1,1$ and letting 
\begin{equation}\label{eq:srulerLong}
    \srulerLong{n} = \srulerLong{n-1}, n, \reverse{\complement{\srulerLong{n-1}}},  n, \srulerLong{n-1}, n
\end{equation}
As in Section \ref{sec:sequences_QRGC_sruler}, we obtain the formula for $\srulerLong{n}$ in \eqref{eq:srulerLong} by omitting the last term $\reverse{\complement{\srulerQuat{n-1}}}$ from \eqref{eq:srulerQuat_reverse}.

For example, 
\begin{align*}
    \srulerLong{3} =&\ \srulerLong{2},3,\reverse{\complement{\srulerLong{2}}},3,\srulerLong{2},3 \\
    =&\  1,1,1,2,-1,-1,-1,2,1,1,1,2,3,-2,-1,-1,-1,-2,1,1,1,-2, \\
    =&\ -1,-1,-1,3,1,1,1,2,-1,-1,-1,2,1,1,1,2,3
\end{align*}
Since we omitted the final term from \eqref{eq:srulerQuat_reverse} when constructing \eqref{eq:srulerLong}, the unsigned sequence in \eqref{eq:rulerLong} is not a palindrome, and so we require the reversal of the middle $\srulerLong{n-1}$ term in \eqref{eq:srulerLong}. We have the same relationship seen before between $\rulerLong{n}$ and $\rulerQuat{n}$, as $\srulerLong{n}$ is a sublist of $\srulerQuat{n}$, using the same strikeouts as indicated in \eqref{eq:rulerLong_strike}. This relationship is seen in Table \ref{tab:orders}.

\subsection{Sequence for the Shortest Solution to Ziggu Puzzles}
\label{sec:ruler_short}

In this section, we construct change sequences for the shortest solution to Ziggu puzzles. We represent states in the shortest solution as $q=q_n\cdots q_2q_1\in\SHORT{n}$, and the presence of $i$ in a change sequence indicates that digit $q_i$ should be changed. 

\subsubsection{Change Sequence for the Shortest Solution}
\label{sec:sequences_short_ruler}

We can define an (unsigned) change sequence for the shortest solution recursively as follows. Let $\rulerShort{n}=\rulerShort[1]{n}=1,1,1$, and for all $n>1$, let $\rulerShort{n}=\rulerShort[1]{n},$ where

\begin{align}
    \rulerShort[f]{n} 
    &= \begin{cases}
        \rulerShort[f]{n-1},n,\reverse{\left(\rulerShort[0]{n-1}\right)},n,\rulerShort[0]{n-1},n &\text{if }f=1 \\
         \phantom{\rulerShort[f]{n-1},}\ n,\reverse{\left(\rulerShort[0]{n-1}\right)},n,\rulerShort[0]{n-1},n&\text{if }f=0
    \end{cases}\label{eq:rulerShort_reverse} \\
    &= \begin{cases}
        \rulerShort[f]{n-1},n,\rulerShort[0]{n-1},n,\rulerShort[0]{n-1},n &\text{if }f=1 \\
         \phantom{\rulerShort[f]{n-1},}\ n,\rulerShort[0]{n-1},n,\rulerShort[0]{n-1},n&\text{if }f=0
    \end{cases}\label{eq:rulerShort_norev}
\end{align}

We have introduced a flag variable $f$, which indicates whether or not to recursively include the prefix $\rulerShort[f]{n-1}$. This has the advantage of being efficient to implement computationally. Since $\rulerShort[0]{n-1}$ is a palindrome, we can omit its reversal, and thus obtain \eqref{eq:rulerShort_norev} from \eqref{eq:rulerShort_reverse}.

For example, we see that 
\begin{align}
\small
    \rulerShort{3} =&\ \rulerShort[1]{3} \nonumber\\
    =&\ \rulerShort[1]{2},3,\rulerShort[0]{2},3,\rulerShort[0]{2},3 \nonumber\\
    =&\ 1,1,1,2,1,1,1,2,1,1,1,2,3,2,1,1,1,2,1,1,1,2,3,2,1,1,1,2,1,1,1,2,3 \nonumber\\
    =&\  1,1,1,2,1,1,1,2,1,1,1,2,3,2,1,1,1,2,1,1,1,\label{eq:rulerShort_strike}\\
    &\ 2,\hcancel[red]{1,1,1,}3,\hcancel[red]{1,1,1,}2,1,1,1,2,1,1,1,2,3\nonumber
\end{align}
where \eqref{eq:rulerShort_strike} is $\rulerLong{3}$ with the states removed to created $\rulerShort{3}$. How $\rulerShort{3}$ is a sublist of $\rulerLong{3}$ is further illustrated in Table \ref{tab:orders}, where one can see how the change sequence $\rulerShort{3}$ acts on the states to produce $\SHORT{3}$.

Unlike the recursive definition of the shortest solution $\SHORT{n}$ given in \eqref{eq:rec_short}, which requires a two-part recursion, the definition of the change sequence $\rulerShort{n}$ requires only one recursion and uses a flag variable $f$. The slice in \eqref{eq:rec_short} is required since the state $03^{n-1}$ is repeated between the two recursive parts, but we don't see this repetition in the change sequence, which is why there is no slice and only a flagged single-part recursion in \eqref{eq:rulerShort_norev}.

\subsubsection{Signed Change Sequence for the Shortest Solution}
\label{sec:sequences_short_sruler}

By Remark \ref{rem:no_signs_needed}, the unsigned change sequence $\rulerShort{n}$ is sufficient to produce the shortest solution to the Ziggu puzzles, as each digit can only be incremented or decremented, depending on the other digits. For completeness purposes, we construct a signed change sequence that, give a state $q=q_n\cdots q_2q_1\in\SHORT{n}$ explicitly includes the information about whether to increment ($i$) or decrement ($-i$) digit $q_i$.

Let $\srulerShort{n}=\srulerShort[1]{n}=1,1,1$, and for all $n>1$, let $\srulerShort{n}=\srulerShort[1]{n},$ where
\begin{align}
    \srulerShort[f]{n} 
    &= \begin{cases}
        \srulerShort[f]{n-1},n,\reverse{\complement{\srulerShort[0]{n-1}}},n,\srulerShort[0]{n-1},n &\text{if }f=1 \\
        \phantom{\srulerShort[f]{n-1},}\ n,\reverse{\complement{\srulerShort[0]{n-1}}},n,\srulerShort[0]{n-1},n&\text{if }f=0
    \end{cases}\label{eq:srulerShort_reverse} \\
    &= \begin{cases}
        \srulerShort[f]{n-1},n,\complementMax{\srulerShort[0]{n-1}},n,\srulerShort[0]{n-1},n &\text{if }f=1 \\
        \phantom{\srulerShort[f]{n-1},}\ n,\complementMax{\srulerShort[0]{n-1}},n,\srulerShort[0]{n-1},n&\text{if }f=0
    \end{cases}\label{eq:srulerShort_norev}
\end{align}
where we use the same flag variable concept as in \eqref{eq:rulerShort_reverse}.
Since $\rulerShort[0]{n}$ is a palindrome, we see that $\reverse{\complement{\srulerShort[0]{n-1}}}=\complementMax{\srulerShort[0]{n-1}},$ where $\complementMax{}$ complements the sign of the maximum entries in a sequence. For example, 
\begin{align*}
\small
    \srulerShort{3} =&\ \srulerShort[1]{3} \\
    =&\ \srulerShort[1]{2},3,\srulerShort[0]{2},3,\srulerShort[0]{2},3 \\
    =&\ 1,1,1,2,-1,-1,-1,2,1,1,1,2,3,-2,-1,-1,-1,\\
    &\ -2,1,1,1,-2,3,2,-1,-1,-1,2,1,1,1,2,3
\end{align*}

We have the same relationship seen before between $\rulerShort{n}$ and $\rulerLong{n}$, as $\srulerShort{n}$ is a sublist of $\srulerLong{n}$, using the same strikeouts as indicated in \eqref{eq:rulerShort_strike}. This relationship is seen in Table \ref{tab:orders}.

\section{Successor and Predecessor Rules: The Next and Previous States}
\label{sec:successor}

Puzzles such as the Gray code and Ziggu puzzles can take a long time to solve, and so it is conceivable that a puzzle solver might set a puzzle down during a solving session to say, make a cup of tea. If the puzzle solver noted their last move or could recall it, they may be able to resume their solution using the approach outlined in Section \ref{sec:solving_resuming}. However, if a puzzle solver has no knowledge as to the previous moves, they would use a \emph{successor rule} to determine the next state to which they should advance the puzzle.

Give a string in a Gray code order, a successor rule computes the next string in the Gray code.


\subsection{Next String in the Binary Reflected Gray Code}
\label{sec:successor_Gray}

next state in BRGC.
this is the same as next state in a solution for some Gray code puzzles but not all.
for example, Hanoi last BRGC state is solved state, so always successor.
in contrast, in Spin-Out solved is $0^n$ which is first in BRGC so always make progress to solution by applying predecessor.
however, need to apply successor or predecessor when returning to start state $1^n$ because it is not at extreme in path; 
see Section \ref{sec:position_comparison} and Figure [ref] for further discussion of this.

\subsection{Next String in the Quaternary Gray Code}
\label{sec:successor_GC}

\begin{definition} \label{def:successorGC}
Let $q = q_n \cdots q_2 q_1 \in \QUAT{n}$ be a quaternary string of length $n$.
Its \emph{quaternary successor} is obtained by incrementing or decrementing digit $q_i$ as follows
\begin{equation} \label{eq:successorGC}
    \nextGC{q} = 
    \begin{cases}
        q_n \cdots q_{i+1} \; q_i{+}1 \; q_{i-1} \cdots q_1 & \text{if $q_i < 3$ and $\displaystyle \sum_{j=i+1}^{n} q_j$ is even} \\
        q_n \cdots q_{i-1} \; q_i{-}1 \; q_{i-1} \cdots q_1 & \text{if $q_i > 0$ and $\displaystyle \sum_{j=i+1}^{n} q_j$ is odd} \\
    \end{cases}
\end{equation}
where $i$ is the minimum index such that one of the two cases holds.
If there is no such value of $i$, then $q$'s quaternary successor $\nextGC{q}$ is undefined.
\end{definition}

\begin{example}
    If $q=012310$, we see that $i=2$ since $q_2=1<3$ and $\sum_{j=i+1}^nq_j=3+2+1+0=6$ is even. \eqref{eq:successorGC} gives that $\nextGC{q}=012320$, so digit $q_2$ was incremented.
\end{example}

\begin{remark} \label{rem:successorGC}
    If $\nextGC{q}$ is undefined for $q \in \QUAT{n}$, then $q = 30^{n-1}$, which is the last string in the quaternary Gray code 
\end{remark}

\subsection{Next State in the Longest Solution to Ziggu Puzzles}
\label{sec:successor_long}

\begin{definition} \label{def:successorLong}
Let $q= q_n \cdots q_2 q_1 \in \LONG{n}$ be a long quaternary string of length $n$.
Its \emph{long successor} is obtained by incrementing or decrementing digit $q_i$ as follows
\begin{equation} \label{eq:successorLong}
    \nextLong{q} = 
    \begin{cases}
        q_n \cdots q_{i+1} \; q_i{+}1 \; q_{i-1} \cdots q_1 & \text{if $\displaystyle \sum_{j=i+1}^{n} q_j$ is even and $q_i < 3$} \\ 
        q_n \cdots q_{i-1} \; q_i{-}1 \; q_{i-1} \cdots q_1 & \text{if $\displaystyle \sum_{j=i+1}^{n} q_j$ is odd, $q_i > 0$ and $q_{i-1} \neq 3$} \\
    \end{cases}
\end{equation}
where $i$ is the minimum index such that one of the two cases holds.
If there is no such value of $i$, then $q$'s long successor $\nextLong{q}$ is undefined.
\end{definition}

\begin{example}\label{ex:successor_long}
If $q=20103$, then $i=1$, since $q_1=3>0$ and $\sum_{j=i+1}^nq_j=0+1+0+2=3$ is odd. Therefore $\nextLong{q}=20102$, as digit $q_1$ was decremented.
\end{example}

Comparing \eqref{eq:successorGC} and \eqref{eq:successorLong}, we see that the longest solution omits the strings in the quaternary reflected Gray code where a $3$ is followed by digits that are not $3$.

\begin{remark} \label{rem:successorLong}
    If $\nextGC{q}$ is undefined for $q \in \LONG{n}$, then $q = 3^{n}$, which is the last string in the longest solution. 
\end{remark}

\subsection{Next State in the Shortest Solution to Ziggu Puzzles}
\label{sec:successor_short}

\begin{definition} \label{def:successorShort}
Let $q = q_n \cdots q_2 q_1 \in \SHORT{n}$ be a short quaternary string of length~$n$.
Its \emph{short successor} is obtained by incrementing or decrementing digit $q_i$ as follows
\begin{equation} \label{eq:successorShort}
    \nextShort{q} = 
    \begin{cases}
        q_n \cdotsSquish q_{i+1} \; q_i{+}1 \; q_{i-1} \cdotsSquish q_1 & \text{if $\displaystyle\sum_{j=i+1}^nq_j$ is even, $q_i< 3$} \\
        q_n \cdotsSquish q_{i-1} \; q_i{-}1 \; q_{i-1} \cdotsSquish q_1 & \text{if $\displaystyle\sum_{j=i+1}^nq_j$ is odd, $q_i > 0$, $q_{i-1} \neq 3$, $q_{i-1} q_i \neq 03$} \\
    \end{cases}
\end{equation}
where $i$ is the minimum index such that one of the two cases holds.
If there is no such value of $i$, then $q$'s short successor $\nextShort{q}$ is undefined.
\end{definition}

\begin{example}
    Considering the same state as in Example \ref{ex:successor_long}, let $q=20103$. We see that now $i\neq1$ since $q_2q_1=03$. Instead, $i=3$ since $q_3=1<3$ and $\sum_{j=i+1}^n=0+2=2$ is even. Therefore $\nextShort{q}=20203$, where we incremented digit $q_3$.
\end{example}

Comparing \eqref{eq:successorLong} and \eqref{eq:successorShort}, we see that the shortest solution omits the strings in the longest solution where a $03$ is followed by digits that are not $3$. Remark \ref{rem:successorLong} holds for \eqref{eq:successorShort} as well.

\section{Position: Ordering States Using Comparison \& Ranking}
\label{sec:position}

Determining which moves will advance a puzzle in an arbitrary state towards the initial or solved state can seem quite daunting. In Gray code puzzles where the initial state is $1^n$ but the puzzle can have different valid states that also have first digit $1$, one may have to ``backtrack" in the Gray code to find the shortest path to the official starting state of the puzzle. A comparison formula allows a puzzle solver to determine the relative order of the current state and the solved state, thus indicating whether or not to proceed forwards or backwards in the Gray code. 

For the Ziggu puzzles where the initial state $0^n$ and the solved state $3^n$ are the first and final valid states in both the longest and shortest solutions $\LONG{n}$ and $\SHORT{n}$, a comparison formula is more generally useful when one is curious about the relative order of two puzzle states.

We present a comparison formula that can compute the relative order of two states in all four Gray codes we're considering: the binary and quaternary reflected Gray codes, and the longest and shortest solutions to the Ziggu puzzles.

\subsection{A Comparison Formula}
\label{sec:position_comparison}

\begin{figure}
    \centering
    \begin{subfigure}{0.32\textwidth}
        \centering
        \includegraphics[scale=0.0975]{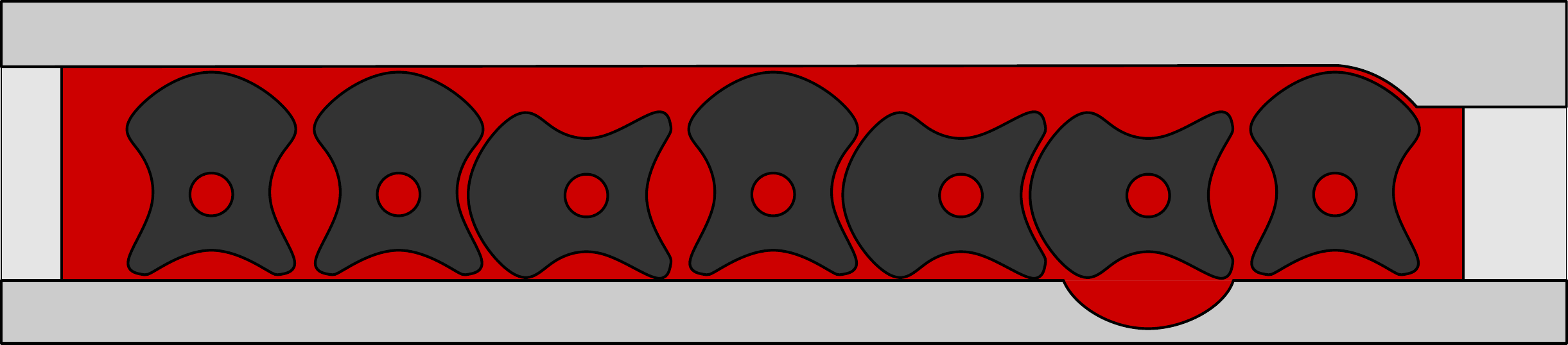} 
        \caption{$t=1101001$.}
        \label{fig:compareSpinOut_1101001}
    \end{subfigure}
    \hfill
    \begin{subfigure}{0.32\textwidth}
        \centering
        \includegraphics[scale=0.0975]{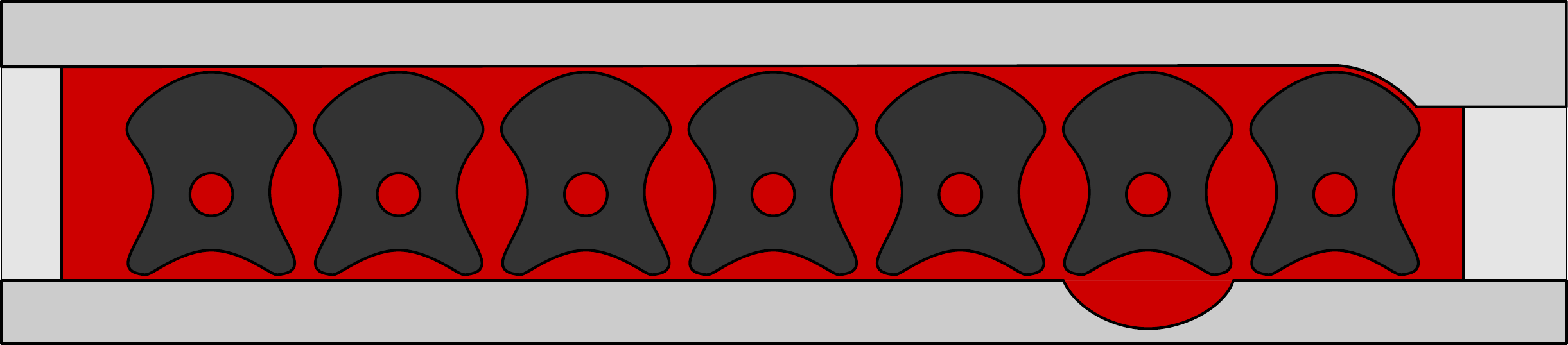} 
        \caption{$u=1111111$.}
        \label{fig:compareSpinOut_1111111}
    \end{subfigure}
    \hfill
    \begin{subfigure}{0.32\textwidth}
        \centering
        \includegraphics[scale=0.0975]{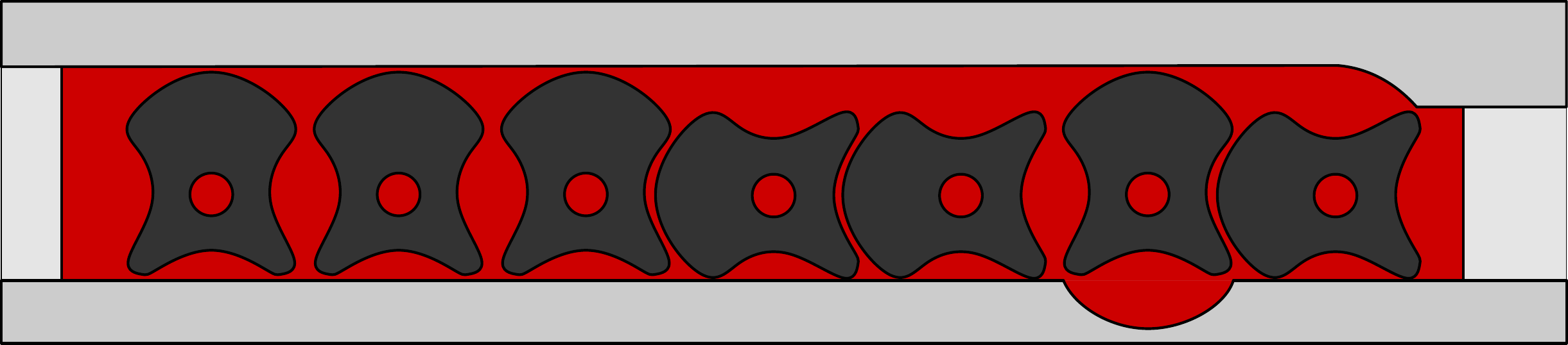} 
        \caption{$v=1110010$.}
        \label{fig:compareSpinOut_1110010}
    \end{subfigure}\caption{
    The initial state of \textit{Spin-Out} is shown in (b) $u=1111111$. If handed the puzzle in state (a) $t=1101001$ or in state (c) $v=1110010$ and asked to return the puzzle to the initial state (b), one would use the comparison formula in \eqref{eqn:comparison} to determine whether to use the successor or predecessor rule to advance the puzzle towards the initial state. 
    The leftmost digit where $t$ and $u$ differ is the fifth digit (read right-to-left), since $t_7t_6t_5=110$ and $u_7u_6u_5=111$. Since the sum of the digits to the left of the fifth digit is even ($t_7+t_6=2$) and $t_5<u_5$, \eqref{eqn:comparison} implies that $t<u$, so $t$ occurs before $u$ in $\BRGC{7}$. Thus, repeatedly apply the successor rule to advance state (a) to state (b).
    Notice that the leftmost digit where $u$ and $v$ differ is the fourth digit since $u_7u_6u_5u_4=1111$ and $v_7v_6v_5v_4=1110$. The sum of the digits to the left of the fourth digit is odd ($u_7+u_6+u_5=3$), and since $v_4<u_4$, \eqref{eqn:comparison} determines that $v>u$. Thus $v$ occurs after $u$ in $\BRGC{7}$, so repeatedly apply the predecessor rule to (c) to return the puzzle to the initial state in (b).}
\label{fig:compareSpinOut}
\end{figure}

Given two distinct states $w,v\in\BRGC{n}$ with $w=w_n\cdots w_2w_1$ and $v=v_n\cdots v_2v_1$, we can determine if $w$ occurs before $v$ in the binary reflected Gray code as follows. Let $w_k$ be the leftmost digit where $w$ and $v$ differ, i.e, $k$ is the index such that $w_i=v_i$ for all $k<i\leq n$ and $w_k\neq v_k.$ We compute that $w$ occurs before $v$ as follows.
\begin{equation}\label{eqn:comparison}
    w < v\text{ if } \begin{cases}
        w_k<v_k &\text{if }\sum_{i=k+1}^n w_i\text{ is even} \\
        w_k>v_k &\text{if }\sum_{i=k+1}^n w_i\text{ is odd.} \\
    \end{cases}
\end{equation}
\begin{example}
    We see that $101<100$ in $\BRGC{3}$, since $k=1$, $w_3+w_2=1+0$ is odd, and $1>0$.
\end{example}
The parity condition on the sum of the digits to the left of $w_k$ captures the number of times that the recursive sublists are written in reverse order. Figure \ref{fig:compareSpinOut} illustrates the use of this comparison formula for the Spin Out puzzle.

Equation \ref{eqn:comparison} also gives the relative order for distinct states $w$ and $v$ in the quaternary reflect Gray code $\QUAT{n}$, as the  parity condition still captures this same sublist reversal behavior. 
\begin{example}
    We see that $11013<11023$ in $\QUAT{5}$ since $k=2$, $w_3+w_4+w_5=2$ is even, and $1<2$. 
\end{example}

Since both the longest and shortest solutions are sublists of $\QUAT{n}$, we see immediately that if we are given distinct states $w,v\in\LONG{n}$ or $w,v\in\SHORT{n},$ then \eqref{eqn:comparison} provides their relative order within these respective lists. 

\begin{figure}
    \centering
    \begin{subfigure}{0.32\textwidth}
        \centering
        \includegraphics[scale=0.0975]{SpinOut1111111.pdf} 
        \caption{$1111111$ has rank $85$.}
        \label{fig:rankSpinOut_1111111}
    \end{subfigure}
    \hfill
    \begin{subfigure}{0.32\textwidth}
        \centering
        \includegraphics[scale=0.0975]{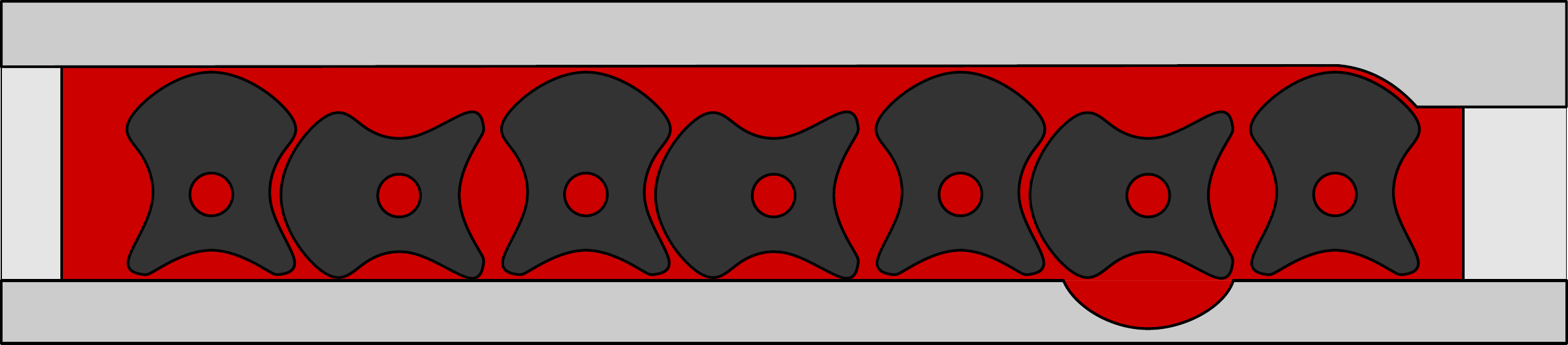} 
        \caption{$1010101$ has rank $102$.}
        \label{fig:rankSpinOut_1010101}
    \end{subfigure}
    \hfill
    \begin{subfigure}{0.32\textwidth}
        \centering
        \includegraphics[scale=0.0975]{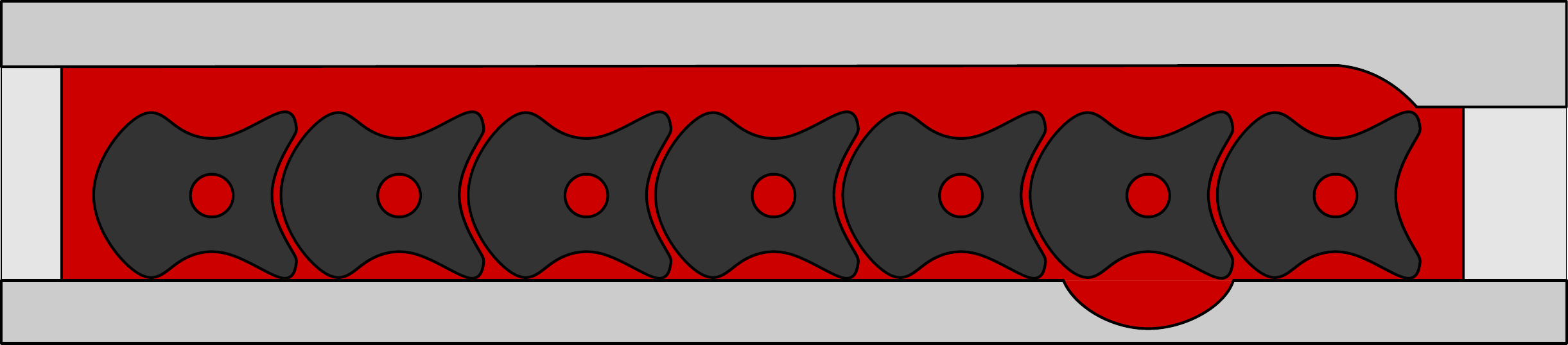} 
        \caption{$0000000$ has rank $0$.}
        \label{fig:rankSpinOut_0000000}
    \end{subfigure}
    \caption{
    The commercial version of \emph{Spin-Out} has $n=7$ dials with start (a) $1111111$ and solved (c) $0000000$ states.
    How many moves are needed to solve the puzzle from (b) $1010101$?
    The answer depends on the ranks of the binary encodings in the binary reflected Gray code for $n=7$.
    The absolute difference between (b) and (c) is $|102-0|=102$, so eighty-five moves are required.
    Similarly, $|102-85|=17$ moves convert it to the start state.
    Since the former quantity is larger than the distance between the start and final states $|85-0| = 85$, we know that the solver has gone the wrong way!
    That is, state (b) is not on the shortest solution.}
    \label{fig:rankSpinOut}
\end{figure}

Beyond understanding the relative order of states within solution sets, a puzzle solver may want to know exactly how many states remain to solve a puzzle. This problem is known as \textit{ranking}, where given a state $w$, one computes the order in the given list, where states are numbered $0,1,2$, and so on. 

\subsection{A ranking formula for Gray Code Puzzles}

To rank the binary reflected Gray code, complement each digit $b$ if the number of $1$s to the left of $b$ is odd, and then convert the resulting string to base $10$. This is computed by the following formula, where $w=w_nw_{n-1}\cdots w_1$ is the given word.

\begin{equation}\label{eqn:rank-brgc}
    \rankBRGC{w} = \sum_{k=1}^n 2^{k-1}a_k\text{ where } a_k=\begin{cases}
        w_k & \text{if }\sum_{j=k+1}^n w_j\text{ is even} \\
        1-w_k & \text{if }\sum_{j=k+1}^n w_j\text{ is odd.} \\
    \end{cases}
\end{equation}

\begin{example}
    The binary word $10110010$ has rank $220$ since $1\overline{01}1\overline{001}0=11011100$, which is $220$ is base $10$. 
\end{example}
An example of ranking Spin-Out states is given in Figure \ref{fig:rankSpinOut}. One can implement \eqref{eqn:rank-brgc} as an efficient algorithm using bitwise operations. 

\begin{figure}[ht]
    \centering
    \captionsetup[subfigure]{justification=centering}
    \begin{subfigure}{0.3\textwidth}
        \includegraphics[width=\textwidth,page=16]{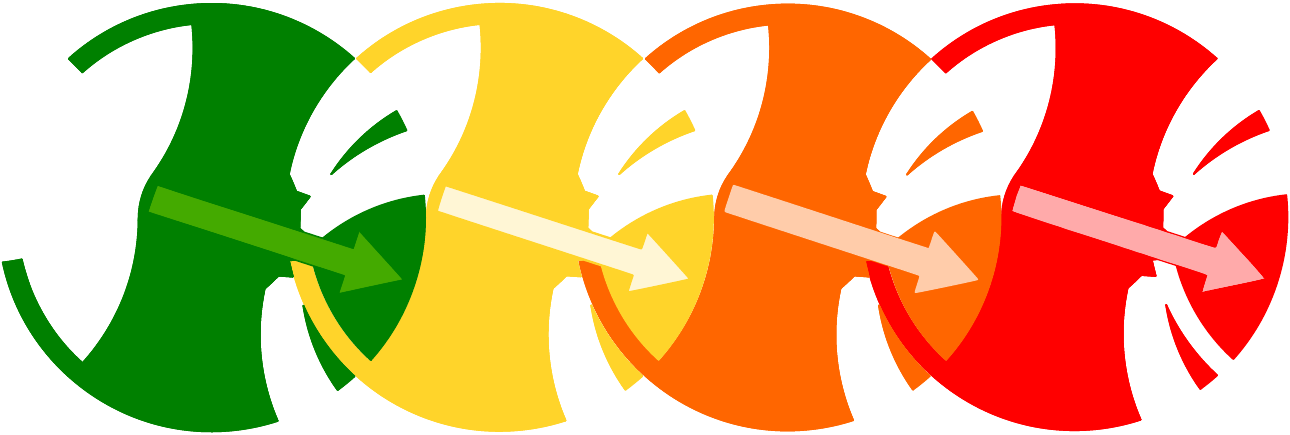}
        \caption{$q=\digitd0\digitc0\digitb0\digita0$\\
        \dial0\dial0\dial0\dial0\\
        \mazefour0000\\
        $\rankShort{q}=0$ \\
        $\rankLong{q}=0$ \\
        $\rankQuat{q}=0$}
        \label{fig:zf-0000}
    \end{subfigure}
    \hfill
    \begin{subfigure}{0.3\textwidth}
        \includegraphics[width=\textwidth,page=17]{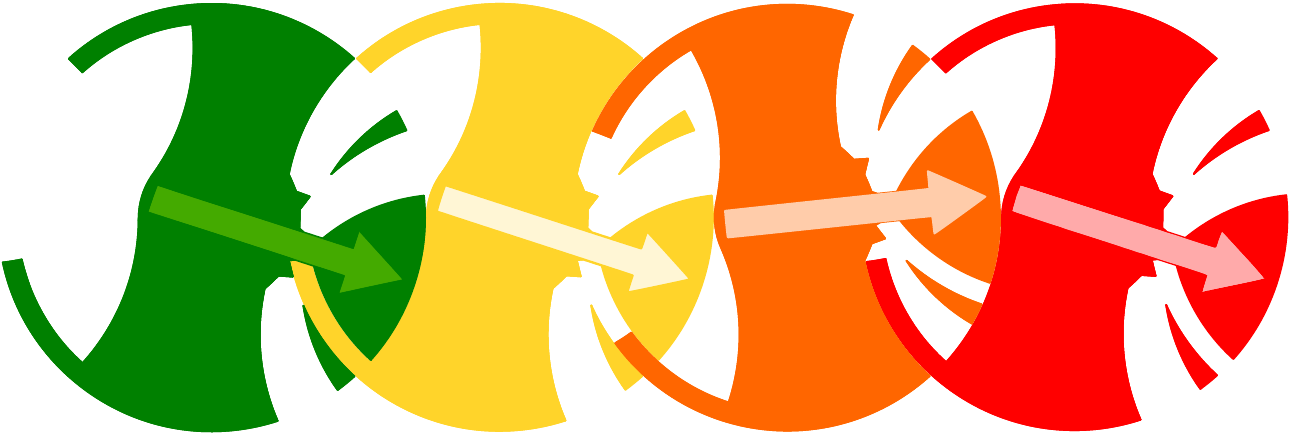}
        \caption{$q=\digitd0\digitc0\digitb2\digita0$\\
        \dial0\dial0\dial2\dial0 \\
        \mazefour0020 \\
        $\rankShort{q}=8$ \\
        $\rankLong{q}=8$ \\
        $\rankQuat{q}=8$}
        \label{fig:zf-0020}
    \end{subfigure}
    \hfill
    \begin{subfigure}{0.3\textwidth}
        \includegraphics[width=\textwidth,page=18]{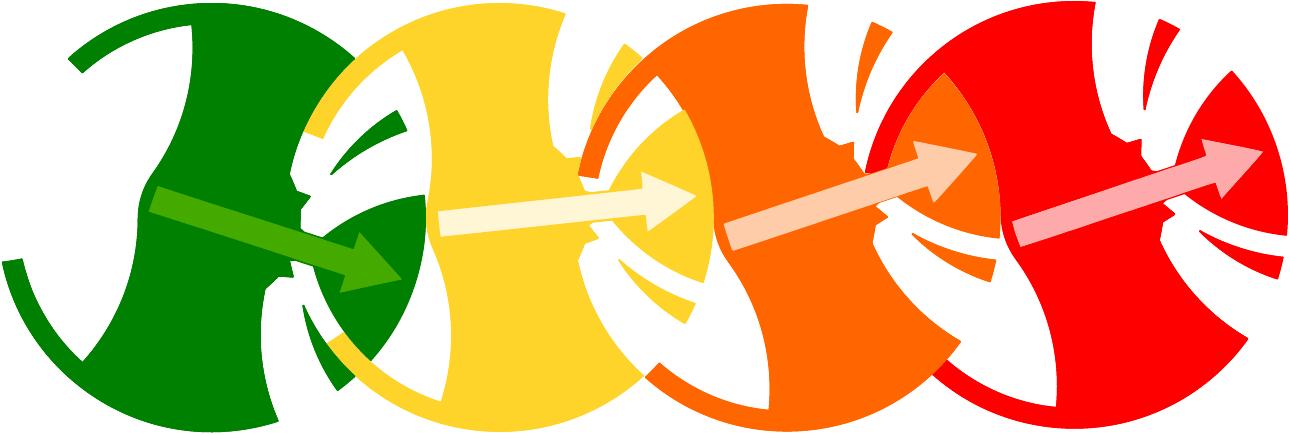}
        \caption{$q=\digitd0\digitc2\digitb3\digita3$\\
        \dial0\dial2\dial3\dial3\\
        \mazefour0233\\
        $\rankShort{q}=32$ \\
        $\rankLong{q}=38$ \\
        $\rankQuat{q}=44$}
        \label{fig:zf-0233}
    \end{subfigure}
    \vfill
    \begin{subfigure}{0.3\textwidth}
        \includegraphics[width=\textwidth,page=19]{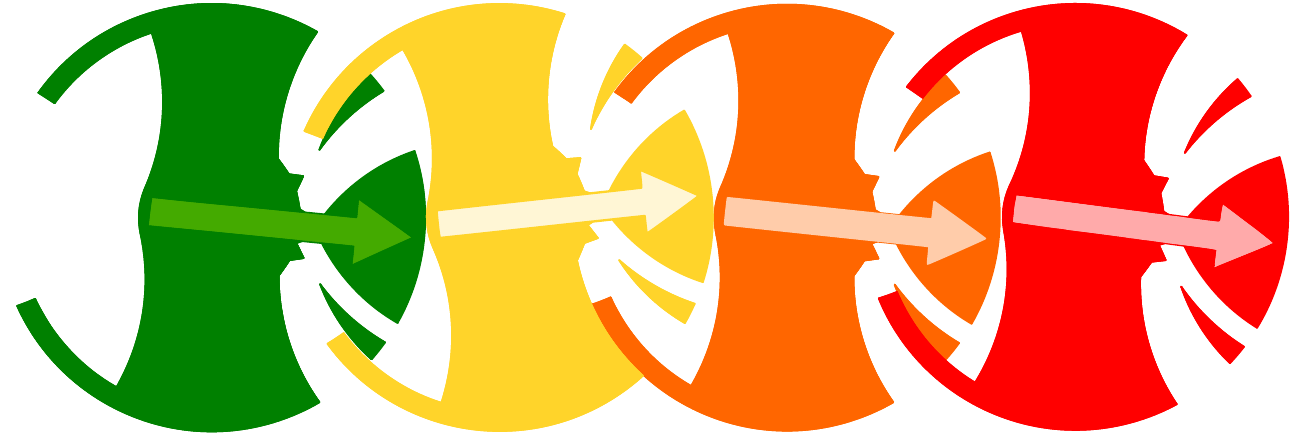}
        \caption{$q=\digitd1\digitc2\digitb1\digita1$\\
        \dial1\dial2\dial1\dial1\\
        \mazefour1211\\
        $\rankShort{q}=41$ \\
        $\rankLong{q}=47$ \\
        $\rankQuat{q}=89$}
        \label{fig:zf-1211}
    \end{subfigure}
    \hfill
    \begin{subfigure}{0.3\textwidth}
        \includegraphics[width=\textwidth,page=20]{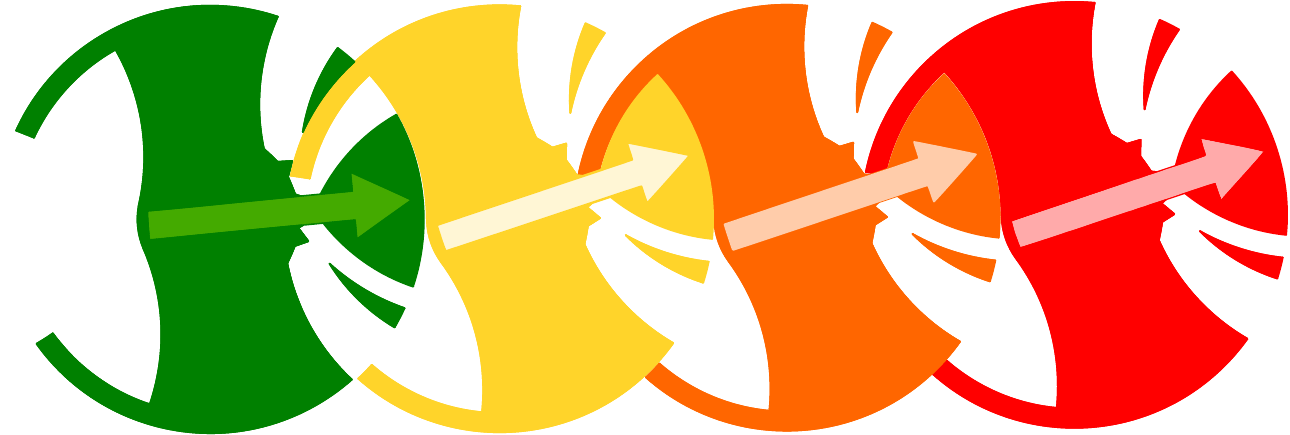}
        \caption{$q=\digitd2\digitc3\digitb3\digita3$\\
        \dial2\dial3\dial3\dial3\\
        \mazefour2333\\
        $\rankShort{q}=77$ \\
        $\rankLong{q}=119$ \\
        $\rankQuat{q}=179$}
        \label{fig:zf-2333}
    \end{subfigure}
    \hfill
    \begin{subfigure}{0.3\textwidth}
        \includegraphics[width=\textwidth,page=21]{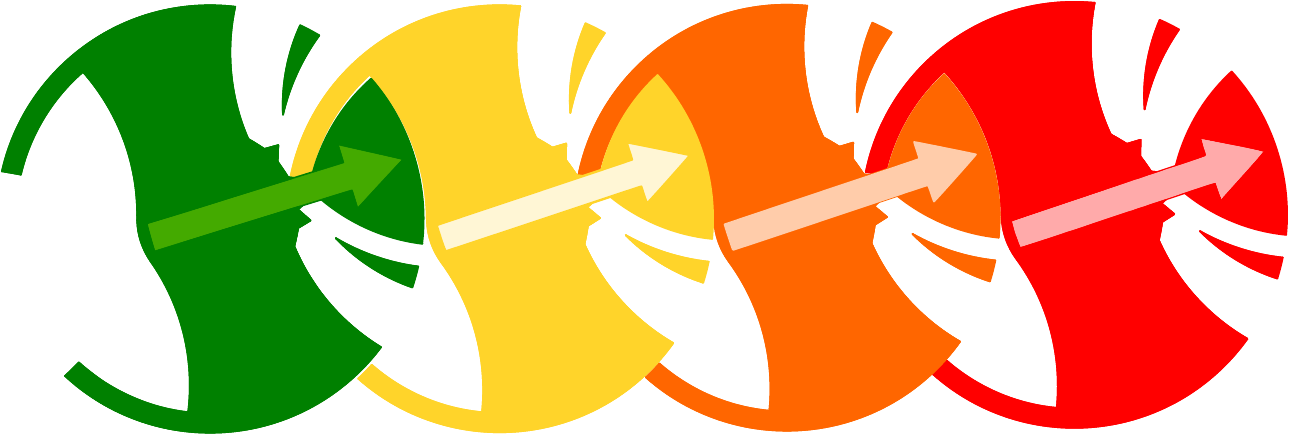}
        \caption{$q=\digitd3\digitc3\digitb3\digita3$\\
        \dial3\dial3\dial3\dial3\\
        \mazefour3333\\
        $\rankShort{q}=78$ \\
        $\rankLong{q}=120$ \\
        $\rankQuat{q}=204$}
        \label{fig:zf-3333}
    \end{subfigure}
    \caption{Six states of Zigguhooked with $4$ dials, with their corresponding base 4 state, their maze representation, and their rank in the shortest solution, the longest solution, and the quaternary grade code, respectively.}
    \label{fig:rank-ziggu}
\end{figure}

In the subsequent sections, we present ranking formulas and algorithms for the Ziggu puzzles. For an example of ranking Zigguhooked, see Figure \ref{fig:rank-ziggu}.

\subsection{A ranking formula for the Quaternary reflected Gray code}

Given a word $w=w_n\cdots w_2w_1$ of length $n$ in the Quaternary reflected Gray code, we can compute the rank of $w$ using the following formula. The ranking is performed by first replacing each digit $w_k$ with $3-w_k$ if the sum of the digits to the left of $w_k$ is odd, and then converting the resulting string from base $4$ to base $10$.
\begin{equation}\label{eqn:rank-quat}
    \rankQuat{w} = \sum_{k=1}^n 4^{k-1}a_k\text{ where } a_k=\begin{cases}
        w_k & \text{if }\sum_{j=k+1}^n w_j\text{ is even} \\
        3-w_k & \text{if }\sum_{j=k+1}^n w_j\text{ is odd.} \\
    \end{cases}
\end{equation}
Equation \ref{eqn:rank-quat} can be implemented by an algorithm that has $O(n)$ time complexity. 

The rank of several states in the length $4$ Zigguhooked puzzle as part of the quaternary reflected Gray code is given in Figure \ref{fig:rank-ziggu}. The ranks of all length $3$ base $4$ strings in shown in Table \ref{tab:orders}.

\subsection{A ranking algorithm for the Longest solution to Ziggu Puzzles}

Given a word $w=w_0w_1\cdots w_{n-1}$ of length $n$ in the shortest solution, Algorithm \ref{alg:rank-long} will output the rank of $w$ in $\LONG n$. This algorithm computes the rank of a given state in $O(n)$ time. It is dependent on $\rankQuat w$, as first the rank of $w$ as a word in the Quaternary Gray code is computed, and then the number of words that the longest solution skips that are included in $\QUAT n$ are subtracted.

Figure \ref{fig:rank-ziggu} also exhibits the ranks of several states of Zigguhooked as part of the longest solution, and the ranking for all states in $\LONG{3}$ are shown in Table \ref{tab:orders}.

\begin{algorithm}[ht]
    \caption{A ranking algorithm for the longest solution.}
    \label{alg:rank-long}
    \begin{algorithmic}[1]
        \Procedure{$\rankLong w$}{}
            \Comment{$w=w_0w_1\dots w_{n-1}$ is the word to rank.}
            \If{$w=0^n$}
            \Comment{All zero word is always first.}
                \State \Output $0$
            \EndIf
            \State{$b\gets0$}
            \State{$l\gets0$}
            \State{$m\gets 0$}
            \State{$s\gets0$}
            \While{$w_m=0$}
                \State{$m\gets m+1$}
            \EndWhile
            \LComment{$m$ is the index of the leftmost nonzero letter.}
            \For{$i\in\set{m,m+1,\dots,n-3}$}
                \If{$l\equiv1\bmod2\lor(i>m\land w_{i-1}=3)$}
                    \State{$m\gets\lfloor(4-w_i)/2\rfloor$}
                \Else
                    \State{$m\gets\lfloor(w_i+1)/2\rfloor$}
                \EndIf
                \State{$l\gets l+w_i$}
                \If{$i=m$}
                    \State{$s\gets s + m\cdot A033114(n-2-i)$}
                \ElsIf{$i=m+1$}
                    \State{$b\gets 2w_m$}
                    \State{$s\gets s + (m+b)\cdot A033114(n-2-i)$}
                \Else
                    \If{$l-w_i-w_{i-1}\equiv1\bmod2\lor(i>1\land w_{i-2}=3)$}
                        \State{$p\gets 2\cdot\max(2-w_{i-1},0)$}
                    \Else
                        \State{$p\gets 2w_{i-1}$}
                    \EndIf
                    \State{$b\gets 3b+p$}
                    \State{$s\gets s + (m+b)\cdot A033114(n-2-i)$}
                \EndIf
            \EndFor
            \State \Output $\rankQuat{w}-6s$
        \EndProcedure
        \Procedure{A033114}{$n$}
            \State \Output $\lfloor4^{n+1}/15\rfloor$
        \EndProcedure
    \end{algorithmic}
\end{algorithm}

\subsection{A ranking algorithm for the Shortest solution to Ziggu Puzzles}

Given a word $w=w_0w_1\cdots w_{n-1}$ of length $n$ in the shortest solution, Algorithm \ref{alg:rank-short} will return the rank of $w$ in $\SHORT n$. This algorithm computes the rank of a given state in $O(n)$ time. It is dependent on $\rankLong w$, as first the rank of $w$ as a word in the longest solution is computed, and then the number of words that the shortest solution skips that are included in $\LONG w$ are subtracted.

To give a complete overview of the ranking, Figure \ref{fig:rank-ziggu} also includes the rank of several Zigguhooked as part of the shortest solution. Table \ref{tab:orders} includes the ranking of all states in $\SHORT{3}$ as well.

\begin{algorithm}[ht]
    \caption{A ranking algorithm for the shortest solution.}
    \label{alg:rank-short}
    \begin{algorithmic}[1]
        \Procedure{$\rankShort w$}{}
            \Comment{$w=w_0w_1\dots w_{n-1}$ is the word to rank.}
            \If{$w=0^n$}
            \Comment{All zero word is always first.}
                \State \Output $0$
            \EndIf
            \State{$a\gets0$}
            \State{$b\gets0$}
            \State{$l\gets0$}
            \State{$m\gets 0$}
            \State{$s\gets0$}
            \While{$w_m=0$}
                \State{$m\gets m+1$}
            \EndWhile
            \LComment{$m$ is the index of the leftmost nonzero letter.}
            \For{$i\in\set{m,m+1,\dots,n-3}$}
                \If{$l\equiv1\bmod2\lor(i>m\land w_{i-1}=3)\lor(i>1\land w_{i-2}\neq0\land w_{i-1}w_i=03)$}
                    \State{$m\gets\lfloor(3-w_i)/2\rfloor$}
                \Else
                    \State{$m\gets\lfloor w_i/2\rfloor$}
                \EndIf
                \State{$l\gets l+w_i$}
                \If{$i=m$}
                    \State{$s\gets s + m\cdot A003462(n-2-i)$}
                \ElsIf{$i=m+1$}
                    \State{$b\gets w_m$}
                    \State{$a\gets m+b$}
                    \State{$s\gets s + a\cdot A003462(n-2-i)$}
                \Else
                    \If{$w_i=3\land w_{i-1}\in\set{0,3}$}
                        \State{$a\gets 2a+1$}
                        \State{$s\gets s + a\cdot A003462(n-2-i)$}
                    \Else
                        \If{$(l-w_i-w_{i-1}\equiv1\bmod2\land w_{i-1}=1)\lor(l-w_i-w_{i-1}\equiv0\bmod2\land w_{i-1}=2)$}
                            \State{$p\gets 2$}
                        \Else
                            \State{$p\gets 1$}
                        \EndIf
                        \State{$b\gets 2b+p$}
                        \State{$a\gets b+m$}
                        \State{$s\gets s + a\cdot A003462(n-2-i)$}
                    \EndIf
                \EndIf
            \EndFor
            \State \Output $\rankLong{w}-6s$
        \EndProcedure
        \Procedure{A003462}{$n$}
            \State \Output $(3^n-1)/2$
        \EndProcedure
    \end{algorithmic}
\end{algorithm}

\section{Loopless Algorithms for Ziggu Puzzles}
\label{sec:loopless}



Combinatorial generation refers to the efficient generation of every object of a particular type \cite{ruskey2003combinatorial}.
In this context, the best type of algorithm is a \emph{loopless algorithm}, which creates each successive object in worst-case $O(1)$-time.
More specifically, a single object is stored in memory, and the generation algorithm makes a constant-sized modification to obtain the next object.
Loopless algorithms were pioneered by Ehrlich \cite{Ehr1973} and exist for most common types of combinatorial objects (e.g., all $n$-bit binary strings, all permutations of $\{1,2,\ldots,n\}$, etc.).
In this section we provide loopless algorithms for the shortest solutions to Ziggu puzzles.

\subsection{Loopless Algorithms for the Shortest Solution}
\label{sec:loopless_short}

Algorithm \ref{alg:fast-ziggu} requires one index variable and can generate each word in constant time. It stores the current word $w$ and current right-to-left $1$-index $i$. 

If $i=1$, so the rightmost digit is to be changed, it is either increment thrice from $0\to1\to2\to3$ or decremented thrice from $3\to2\to1\to0$. Once these three states have been output, we increment $i$, and so consider the second-to-the-right digit in the next iteration of the loop (lines 6-15). 

Algorithm \ref{alg:fast-ziggu} begins with the starting word of $n$ zeroes and initial index $i=1$, and so the algorithm will first change the rightmost digit. 

If $i>1$, the algorithm considers the parity of $w_i$ and $w_{i-1}$, the digit to its right. If $w_i$ and $w_{i-1}$ have opposite parity, we increment $w_i$. If they have the same parity, we decrement $w_i$. Finally, if this change results in $w_i$ being extremal, equal to $0$ or $3$, we increment $i$, and so the next loop iteration will modify the digit to the left. If $w_i$ is not extremal, we decrement $i$, and thus next consider the digit to the right of $w_i$ (lines 16-25).

If the puzzle solver is willing to sacrifice an additional $n$ bits of auxiliary memory, this algorithm can be simplified considerably. Algorithm \ref{alg:fast-ziggu} also generates each subsequent word in worst-case constant time. Just like Algorithm \ref{alg:fast-ziggu-dir}, Algorithm \ref{alg:fast-ziggu} stores the current word $w$ and the current right-to-left $1$-index. The additional $n$ bits of memory are required for the array of directions $d$, where each entry is either $1$ or $-1$. 

Algorithm \ref{alg:fast-ziggu-dir} has the same initial state as Algorithm \ref{alg:fast-ziggu}, with initial word $w=0^n$ and starting index $i=1$. The array of direction $d$ is initialized to have every entry be a $1$, indicating that each digit should be increased (lines 2-4).

In every loop iteration, Algorithm \ref{alg:fast-ziggu-dir} changes $w_i$ by the value $d[i]$, the $i$-th right-to-left $1$-indexed value in the array $d$ (line 7). If this change results in $w_i$ being extremal, so $w_i$ equals $0$ or $3$, $d[i]$ is set from increasing $(1)$ to decreasing $(-1)$ or vice versa (line 9). In the next loop iteration, the digit to the left of $w_i$ will be changed (line 10). If $w_i$ was instead changed to a nonextremal value ($1$ or $2$), then the next loop iteration will consider the digit to the right of $w_i$, if $w_i$ is not the rightmost digit (lines 11-12). If $w_i$ is the rightmost digit, $w_i$ will be changed again in the next loop iteration. 
\begin{algorithm}[!ht]
    \caption{An algorithm for generating the shortest solution to the Ziggu puzzles. Each successive state is generated in worst-case $O(1)$ time, making it a loopless algorithm. It uses one index variable, which requires $O(\log(n))$ memory.}
    \label{alg:fast-ziggu}
    \begin{algorithmic}[1]
        \Procedure{\FastZiggu}{$n$}
        \State $w\gets 0^n$
        \Comment{The start word $0^n$, $n$ copies of $0$.}
        \State $i\gets 1$
        \Comment{$1$-based (right-to-left) index for changing $w$.}
        \State \Visit $w$, None
        \While{$i\leq n$}
        \If{$i=1$}
        \Comment{If we are changing the rightmost digit.}
        \If{$w_i=0$}
        \For{$m\in\set{1,2,3}$}
        \State{$w_i \gets w_i+1$}
        \Comment{Increment $w_i$ three times if $w_i=0$.}
        \State \Visit $w$, $i$
        \EndFor
        \Else
        \For{$m\in\set{1,2,3}$}
        \State{$w_i \gets w_i-1$}
        \Comment{Decrement $w_i$ three times if $w_i=3$.}
        \State \Visit $w$, $i$
        \EndFor
        \EndIf
        \State{$i\gets i+1$}
        \Comment{Consider the digit to the left.}
        \Else
        \If{$(w_i+w_{i-1})\equiv1\bmod2$}
        \Comment{If $w_i$ and $w_{i-1}$ have opposite parity}
        \State{$w_i\gets w_i+1$}
        \Comment{Increment $w_i$}
        \Else
        \Comment{If $w_i$ and $w_{i+1}$ have the same parity}
        \State{$w_i\gets w_i-1$}
        \Comment{Decrement $w_i$}
        \EndIf
        \State \Visit $w$, $i$
        \If{$w_i\in \set{0,3}$}
        \Comment{If $w_i$ is extremal}
        \State{$i\gets i+1$}
        \Comment{Consider the digit to the left}
        \ElsIf{$i>1$}
        \State{$i\gets i-1$}
        \Comment{Consider the digit to the right}
        \EndIf
        \EndIf
        \EndWhile
        \EndProcedure
    \end{algorithmic}
\end{algorithm}

\begin{algorithm}[!ht]
    \caption{An algorithm for generating the shortest solution to the Ziggu puzzles. Each successive state is generated in worst-case $O(1)$ time, making it a loopless algorithm. The algorithm uses one index variable, which requires $O(\log(n))$ memory, and $n$ bits of additional memory to store directions.}
    \label{alg:fast-ziggu-dir}
    \begin{algorithmic}[1]
        \Procedure{\FastZigguDir}{$n$}
        \State $d\gets [1]*n$
        \Comment{The starting directions, a length $n$ array of $1$s.}
        \State $w\gets 0^n$
        \Comment{The starting word of $n$ zeroes.}
        \State $i\gets 1$
        \Comment{$1$-based (right-to-left) index for changing $w$.}
        \State \Visit $w$, None
        \While{$i\leq n$}
        \State $w_i \gets w_i+d[i]$
        \Comment{Change $w_i$ by $\pm1$, depending on value of $d[i]$.}
        \If{$w_i\in \set{0,3}$}
        \Comment{Check if $w_i$ is now extremal.}
        \State $d[i]\gets -d[i]$
        \Comment{Change $d[i]$ from $\pm1$ to $\mp1$.}
        \State $i \gets i+1$
        \Comment{Consider the letter to the left.}
        \ElsIf{$i<n-1$}
        \State $i \gets i-1$
        \Comment{Consider the letter to the right.}
        \EndIf
        \State \Visit $w$, $i$
        \EndWhile
        \EndProcedure
    \end{algorithmic}
\end{algorithm}



\section{A Bijection with Nurikabe puzzles}
\label{sec:nurikabe}

Since the discovery of the binary reflected Gray code, many exciting bijections between the Gray code and other combinatorial objects have been uncovered. There are many such sequences whose rich connections to other mathematical objects have yet to be discovered. One such sequence is OEIS A101946, the sequence whose $n$-th term is the length of the shortest solution to an $n$-Ziggu puzzle. In this section, we present an explicit bijection between the set of $2\times n$ Nurikabe grids and the set of states in the shortest solution to an $n$-Ziggu puzzle. 

Nurikabe is a pencil and paper puzzle by Japanese puzzle publisher Nikoli \cite{Nik2021}. The starting state is a grid of white squares, some of which contain positive integers. The goal is to color in the squares without numbers in either black or white such that each number $k$ lives in exactly one edge-connected ``island'' of exactly $k$ white squares, the black cells are all edge-connected, and no $2\times2$ square is colored entirely in black.

We consider a counting problem inspired by Nurikabe.

\begin{figure}[ht]
    \centering
    \includegraphics[width=\textwidth]{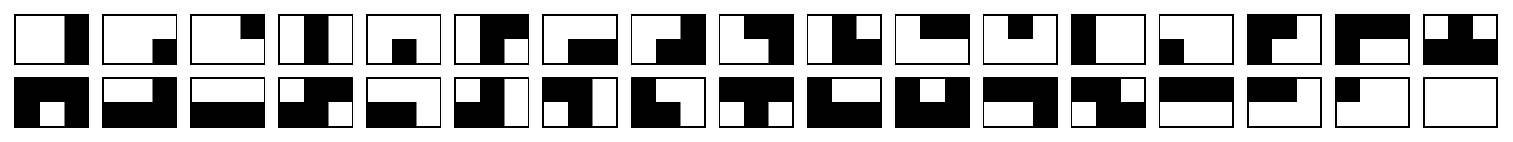}
    \caption{The $34$ distinct $2\times 3$ Nurikabe grids.}
    \label{fig:nurikabe4}
\end{figure}

\begin{definition}\label{def:nurikabe}
    A \emph{$2\times n$ Nurikabe grid} is a coloring of a $2\times n$ grid by black and white such that the following conditions are satisfied.
    \begin{enumerate}
        \item The black squares are edge connected; i.e., they form a polyomino.
        \item No $2\times2$ subgrid of the grid consists entirely of black squares.
    \end{enumerate}
    Let $\Nurikabe{n}$ denote the set of all $2\times n$ Nurikabe grids.
\end{definition}
The $2\times3$ Nurikabe grids are shown in Figure \ref{fig:nurikabe4}.

\begin{theorem}\label{thm:nurikabe-count}
    Let $a_n=|\Nurikabe{n}|.$ Then
    \begin{equation*}
        a_n = 6\cdot2^n-3n-5 \quad\text{for all }n\geq0.
    \end{equation*}
\end{theorem}
\begin{proof}[Proof Sketch]
    Let $b_n$ be the number of $2\times n$ Nurikabe grids that have at least one black cell in all $n$ columns, and let $c_n$ be the number of $2\times n$ Nurikabe grids with at least one column consisting only of white squares. Thus $a_n=b_n+c_n$. One can show using counting arguments that $b_0=0$, $c_0=1$, $b_n = 3\cdot 2^{n-1}$, and $c_n=9\cdot 2^{n-1}-3n-5$ for all $n\geq1$.
\end{proof}

This sequence appears in the Online Encyclopedia of Integer Sequences as sequence A101946. A comment stating Theorem \ref{thm:nurikabe-count} without proof was contributed to the encyclopedia by Hugo van der Sanden on February 22, 2024.

\begin{theorem}\label{thm:bij}
    For all $n\geq1$, there exists a bijection between $\Nurikabe{n}$, the set of $2\times n$ Nurikabe grids, and $\SHORT{n}$, the sets of states in the shortest solution to an $n$-Ziggu puzzle.
\end{theorem}
A constructive proof amounts to the bijection in Definition \ref{def:bij} below, which one can verify is in fact a bijection.

\begin{definition}\label{def:bij}
    Given a word $w=w_1w_2\dots w_n\in\SHORT{n}$, we will construct a $2\times n$ Nurikabe grid $g_1g_2\dots g_n\in\Nurikabe{n}$, where $g_i$ is the $i$-th column of the grid. Let $s:\SHORT{1}\to\Nurikabe{1}$ be defined by
    \begin{equation*}
        s(w_i) = \begin{cases}
            \colbb, & \text{if }w_i=0  \\[-0.2em]
            \colwb, & \text{if }w_i=1  \\[-0.2em]
            \colbw, & \text{if }w_i=2  \\[-0.2em]
            \colww, & \text{if }w_i=3.
        \end{cases}
    \end{equation*}
    Let $p:\Nurikabe{1}\times\SHORT{1}\to\Nurikabe{1}$ be given by
    \begin{equation*}
        p(c,w_i) =
        \begin{cases}
            \colbb, & \text{if }(c,w_i)\in\set{\Big(\colwb,2\Big),\Big(\colbw,1\Big)} \\[0.4em]
            \colwb, & \text{if }w_i=1\text{ and }c\in\Big\{\colbb,\colwb\Big\}                    \\[0.4em]
            \colbw, & \text{if }w_i=2\text{ and }c\in\Big\{\colbb,\colbw\Big\}                    \\[0.4em]
        \end{cases}
    \end{equation*}

    Let $k$ be such that $w_kw_{k+1}\dots w_n$ is a maximal suffix of the form $03^{n-k}$ or $3^{n-k+1}$. If no such suffix exists, let $k=n+1$. Let $l$ be such that $w_1w_2\dots w_l$ is a maximal prefix of the form $0^l$, if $w$ has a prefix consisting of zeroes. If no such prefix exists, let $l=0$.

    Define $\sigma:\SHORT{n}\to\Nurikabe{n}$ to be
    \begin{equation*}
        \sigma(w) = g_1g_2\dots g_n,
    \end{equation*}
    where working right to left, we have
    \begin{equation*}
        g_i = \begin{cases}
            s(w_i),         & k\leq i\leq n, \\
            p(g_{i+1},w_i), & l < i < k,     \\[-0.2em]
            \colww,         & 1\leq i\leq l.
        \end{cases}
    \end{equation*}
\end{definition}

\begin{figure}[ht]
    \centering
    \includegraphics[width=0.8\textwidth]{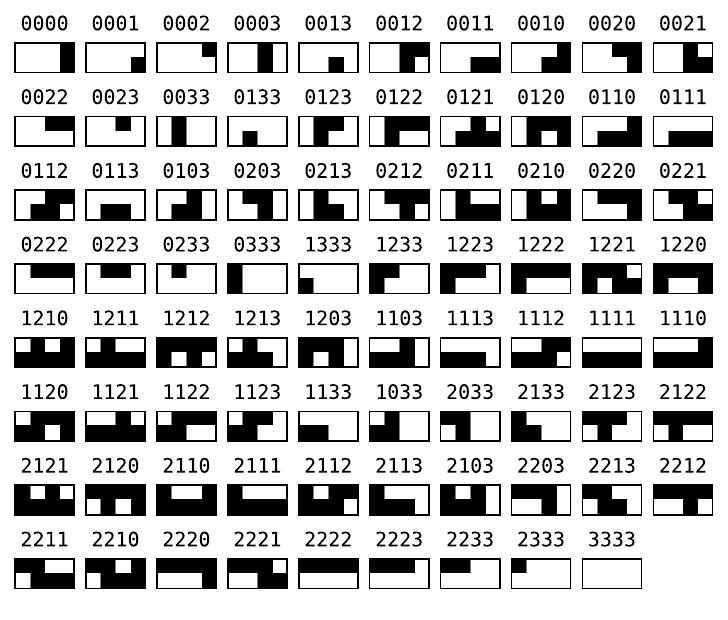}
    \caption{The bijection between $\Nurikabe{4}$, the set of $2\times 4$ Nurikabe grids, and $\SHORT{4}$, the set of states in the shortest solution to a $4$-Ziggu puzzle, from Theorem \ref{thm:bij}.}
    \label{fig:bij}
\end{figure}

Figure \ref{fig:bij} exhibits the bijection for $n=4$.

\begin{remark}
    Let $w=w_nw_{n-1}\dots w_1\in\SHORT{n}$, and let $\rho(w)=r(w_n)r(w_{n-1})\dots r(w_1)$, where
    \begin{equation*}
        r(w_i) = \begin{cases}
            0 & \text{if }w_i = 0 \\[-0.2em]
            1 & \text{if }w_i =2  \\[-0.2em]
            2 & \text{if }w_i =1  \\[-0.2em]
            3 & \text{if }w_i =3.
        \end{cases}
    \end{equation*}
    Then reflecting the $2\times n$ Nurikabe grid $\sigma(r(w))$ about its horizontal axis corresponds to the grid $\sigma(w)$.
\end{remark}

\section{Final Remarks}
\label{sec:final}

Preliminary results in this paper were shared by the second author at a \emph{Gathering for Gardner --- Celebration of Mind} presentation.
This led to several enjoyable communications with the community including Bram Cohen, Oskar van Deventer, and Goetz Schwandtner.
Along with Hirokazu Iwasawa they have been making parallel investigations into Ziggu puzzles and their generalizations.
In particular, they have considered variations to the basic {$\mathsf{S}$}-shaped maze and how this changes the solution lengths.
This research has resulted in new puzzle designs, including van Deventer's recent \emph{Zigguphi} puzzle. 
The authors look forward to playing these new puzzles, and hope that some of the results presented here will help in their analysis.
The authors of this paper have noticed that our results generalize to mazes that continue the back-and-forth pattern (i.e., $2g+1$ rows with $g \geq 1$).


\FloatBarrier
\bibliographystyle{amsplain}
\bibliography{small-math-docs.bib}




\end{document}